\documentclass[10pt,a4paper, oneside]{article}

\usepackage[T1]{fontenc}
\usepackage{authblk}
\usepackage[a4paper,left=3cm,right=2cm,top=2.5cm,bottom=2.5cm]{geometry}

\usepackage{wrapfig}
\usepackage{dcolumn}  
\usepackage{bm}       
\usepackage{amssymb}   
\usepackage{algorithm}
\usepackage{amsmath}
\usepackage{algorithm,algpseudocode}
\usepackage{algorithm}
\usepackage{enumerate}
 \usepackage{booktabs}
\usepackage{multirow}
\usepackage{makecell}
\usepackage{pgfplots}
\usepackage{tikz}
\usetikzlibrary{patterns}
\usepackage{empheq}
\usepackage{calrsfs}
\usepackage[colorinlistoftodos,prependcaption]{todonotes}
\usepackage{xargs}
\usetikzlibrary{arrows}
\usepackage{subfig}
\usepackage{enumitem}
\usepackage[english]{babel}
\usepackage[utf8x]{inputenc}
\usepackage[english]{babel}
\usepackage[utf8x]{inputenc}
\usepackage[cal=boondox]{mathalfa}

\newcolumntype{"}{@{\hskip\tabcolsep\vrule width 1pt\hskip\tabcolsep}}

\newtheorem{remark}{Remark}

\newcommandx{\info}[2][1=]{\todo[linecolor=blue,backgroundcolor=blue!25,bordercolor=black,#1]{#2}}

\definecolor{myblack}{RGB}{53, 53, 53}
\definecolor{myblue}{RGB}{40, 75, 99}
\definecolor{myred}{RGB}{192, 50, 33}
\definecolor{myyellow}{RGB}{255, 166, 48}
\definecolor{mywhite}{RGB}{240, 237, 238}
\definecolor{mygreen}{RGB}{0, 102, 0}

\definecolor{green1}{RGB}{9, 82, 86}
\definecolor{green2}{RGB}{8, 127, 140}
\definecolor{green3}{RGB}{6, 167, 125}
\definecolor{green4}{RGB}{79, 109, 122}
\definecolor{green5}{RGB}{192, 214, 223}
\definecolor{violet}{RGB}{26,69,131}
\def \eps{\bm{\varepsilon}}		

\def \si{{\mathbf{p}_i}}			
\def \Hi{{\mathbf{H}_i}}			
\def \gi{{\mathbf{g}_i}}			
\def \xxi{{\mathbf{x}_i}}			
\def \xxip{{\mathbf{x}_{i+1}}}			

\def \pi{{\mathbf{p}_{i}}}			
\def \0{\mathbf{0}}			%

\def \x{{\mathbf{x}}}			
\def \s{{\mathbf{p}}}			
\def \m{{\mathbf{m}}}			
\def \R{{\mathbb{R}}}			
\def \N {\mathbb{N}}			
\def \g{{\mathbf{g}}}			%
\def \H{{\mathbf{H}}}			%
\def \lb{{\mathbf{l}}}			%
\def \tu{{\mathbf{tu}}}			%
\def \tl{{\mathbf{tl}}}			%
\def \e{{\mathbf{e}}}
\def \l{{\mathbf{l}}}			
\def \u{{\mathbf{u}}}			
\def \y{{\mathbf{y}}}			

\def \h{\mathcal{h}}

\def \Rev{{\mathbf{R}}}			
\def \P{{\mathbf{P}}}			
\def \I{{\mathbf{I}}}			

\def \B{{\mathbf{B}}}			
\def \D{{\mathbf{D}}}			
\def \PI{{\mathbf{\Pi}}}

\newcommand{\utopia}{{\sc Utopia}}
\newcommand{\petsc}{{\sc PETSc}}

\DeclareMathAlphabet{\pazocal}{OMS}{zplm}{m}{n}

\algnewcommand\algorithmiconput{\textbf{Constants:}}
\algnewcommand\algorithmicinput{\textbf{Input:}}
\algnewcommand\algorithmicoutput{\textbf{Output:}}
\algnewcommand{\algorithmicgoto}{\textbf{go to}}%

\algnewcommand\Constants{\item[\algorithmiconput]}
\algnewcommand\Input{\item[\algorithmicinput]}%
\algnewcommand\Output{\item[\algorithmicoutput]}%
\algnewcommand{\Goto}[1]{\algorithmicgoto~\ref{#1}}%

\title{Recursive multilevel trust region method with application to fully monolithic phase-field models of brittle fracture}
\providecommand{\keywords}[1]{\textbf{{Keywords:}} #1}

\author[1]{Alena Kopani\v{c}\'akov\'a\thanks{alena.kopanicakova@usi.ch}}
\author[1]{Rolf Krause}
\affil[1]{Institute of Computational Science, Universit\`a della Svizzera italiana, Lugano 6900, Switzerland}

\begin{document}

%

\thispagestyle{empty}
\newpage

\maketitle
\setcounter{page}{1} 

\begin{abstract}
The simulation of crack initiation and propagation in an elastic material is difficult, as crack paths with complex topologies have to be resolved. 
Phase-field approach allows to simulate crack behavior by circumventing the need to explicitly model crack paths. 
However, the underlying mathematical model gives rise to a non-convex constrained minimization problem.
In this work, we propose a recursive multilevel trust region (RMTR) method to efficiently solve such a minimization problem.
The RMTR method combines the global convergence property of the trust region method and the optimality of the multilevel method. 
The solution process is accelerated  by employing level dependent objective functions, minimization of which provides correction to the original/fine-level problem.
In the context of the phase-field fracture approach, it is challenging to design efficient level dependent objective functions as the underlying mathematical model relies on the mesh dependent parameters. 
We introduce level dependent objective functions that combine fine level description of the crack path with the coarse level discretization. 
The overall performance and the convergence properties of the proposed RMTR method are investigated by means of several numerical examples in three dimensions.
\end{abstract}



\section{Introduction}
Predicting damage, crack and fragmentation patterns is a long-lasting challenge in computational mechanics.
The phase-field approach to fracture addresses this task in an elegant way and has therefore become very popular. 
The basic idea behind the phase-field approach is to regularize sharp crack interfaces in the structure by means of smeared, volumetric crack zones.
This removes the need to represent discontinuities and therefore avoids possible remeshing. 
The volumetric approximation of a crack surface is based on the Ambrosio-Torelli functional \cite{ambrosio1990approximation}, which employs an auxiliary damage variable. 
This damage variable acts as an indicator, which marks the state of the material from intact to fully broken. 
The transition between the broken and intact part of the domain is modeled in diffused manner. 
The size of the transition zone is controlled by the mesh dependent length-scale parameter.
  
Several approaches to phase-field fracture can be found in the literature. 
In the physics community, evolution equations are derived by adapting the phase transition formalism of Landau and Ginzburg \cite{landau1980statistical}. 
In the mechanics community, phase-field fracture approach is based on the variational formulation of the brittle fracture developed by  Francfort and Marigo \cite{francfort1998revisiting}, which provides an extension to classical Griffiths theory of fracture \cite{griffith1921phenomena}.  
The first numerical implementation of the variational phase-field approach is presented by Bourdin et al. in \cite{bourdin2000numerical}. 
Miehe et al.  \cite{miehe2010thermodynamically, miehe2010phase} modify the underlying mathematical model by introducing rate-independent formulation that ensures the local growth of the damage variable. 
Their formulation also introduces anisotropic energy functional, which allows for the degradation of the  tensile energy with increased damage. 
Those modifications result into a thermodynamically consistent phase field formulation of the brittle fracture. 
Since then, the phase-field approach has become popular and was extended in various directions, such as dynamics  \cite{bourdin2011time} \cite{Kuhn_dynamic}, large deformations \cite{del2007variational} \cite{hesch2017framework}  \cite{hesch2014thermodynamically}, cohesive fracture \cite{Vignollet2014} \cite{conti2016phase}, shells and plates \cite{amiri2014phase} \cite{ambati2016phase},  and used within several multiphysics applications  \cite{bilgen2017_pneumatic} \cite {Wick15},  \cite{Heider2016} \cite{abdollahi2012phase}. 

An important feature of any phase-field fracture simulator is to ensure crack irreversibility \cite{miehe2010thermodynamically}, which introduces a lower bound on the phase field variable. 
The resulting constraint can be treated directly in the solution strategy, for example by using active set \cite{facchinei1998active}, or semi-smooth Newton's methods \cite{ulbrich2011semismooth}.
In order to simplify the solution process, alternatives, built on the modified energy, or weak formulations, are widely considered in the literature. 
For example, Bourdin \cite{bourdin2000numerical} transforms inequality constraints into equality by using penalization term from \cite{artina2015}. 
A strategy based on augmented Lagrangian penalization is introduced by Wheeler et. al in \cite{WHEELER201469}. 
Miehe et. al \cite{miehe2010phase} handle the irreversibility by introducing so-called history variable, which replaces maximal tensile energy accumulated within the loading process. 
Disposing the inequality constraints usually causes deterioration of the conditioning of the Hessian matrix \cite{nocedal2006numerical}, or leads to inconsistencies between original energy functional and the modified weak form.  

In this work, the phase-field formulation of fracture is based on the variational framework. 
The fracture problem is modeled by using an energy functional, whose minimizer corresponds to the associated solution. 
Algorithmically, the minimization process can be realized  by using monolithic or staggered solution scheme.  
In the monolithic approach, the energy functional is minimized simultaneously for displacement and phase field \cite{Vignollet2014}. 
This is numerically challenging, as the energy functional of the coupled problem is non-convex.  
In contrast, using a staggered solution scheme,  the energy functional is minimized individually for the displacement and phase field, which gives rise to two convex minimization problems.
This makes staggered solution scheme  very robust and popular in the engineering community. 
However, Gerasimov et. al show in \cite{GERASIMOV2016276}, that if convergent, the monolithic scheme can be faster than the staggered scheme since the staggered scheme tends to underestimate crack speed.  
Therefore, it is an important and challenging task to design efficient and reliable solution strategies for the fully monolithic scheme. 

To our knowledge, there are only a few published results concerning the monolithic framework. 
These advances were accomplished either by modifying the underlying weak formulations or by applying some globalization strategy during minimization process.
Vignollet et. al \cite{vignollet2014phase} augment the weak form by using a dissipation-based arc-length procedure \cite{verhoosel2009dissipation}. 
Heister et. al apply "convexification" trick, where the phase field variable in the critical terms is replaced by its extrapolation in time \cite{heister2015primal}. 
Gerasimov et. al employ a line-search based Newton's solver, which takes into account negative curvature \cite{GERASIMOV2016276}. 
A globalization of Newton's method, based on error oriented approach from \cite{deuflhard2011newton},  is successfully applied by Wick in  \cite{wick2017error}.  
In \cite{wick2017modified}, Wick employs so-called modified Newton's method, which is based on dynamical switching between full Newton's and modified Newton's steps. 
The modified Newton's steps are based on an altered Hessian matrix, where coupling terms are scaled in an appropriate manner. 

All aforementioned globalizations of Newton's method require the solution of a system of linear equations in order to determine step direction. 
Solving those linear systems can become computationally demanding as they are usually very large and ill-conditioned, see Section \ref{section:pf}. 
An ideal choice of the linear solution strategy is a multilevel method due to its optimal complexity and scalability. 
Linear multilevel methods have been developed for more than 50 years and are used extensively in many areas. 
An introduction to multigrid methods can be found for example in \cite{Briggs2000multigrid, hackbusch2013multi}. 
In the context of the phase field fracture models, an algebraic multilevel method is applied by Farrell and Maurini in \cite{farrell2017linear}. 
The authors of  \cite{farrell2017linear} employ smoothed aggregation in order to precondition systems related to the displacement and the phase-field subproblems. 
Also within the staggered solution scheme, the geometric multigrid method is applied by Bilgen et. al in \cite{bilgen2017phase}. 
To the best of our knowledge, the multilevel solution strategy has not been applied yet to the coupled systems arising in the monolithic framework. 

In this work, we present a globally convergent and fully-nonlinear multilevel solution method, for the phase field fracture problems originating in the monolithic framework.
We employ a recursive multilevel trust region method (RMTR) \cite{Gratton2008recursive, Gross2009}, which combines the trust region globalization strategy with the nonlinear multilevel framework.
As other nonlinear multilevel methods, for instance FAS~\cite{brandt1977multi}, or MG/OPT~\cite{Nash2000multigrid}, the RMTR accelerates the solution process by employing level dependent objective functions.
The minimization of these level dependent objective functions provides a coarse level corrections to the original/fine-level problem.
Quality of the coarse level corrections is determined by the characteristics of the level dependent objective functions. 
A construction of an efficient level dependent objective functions for the phase-field fracture problems is a challenging task.
This is due to the fact that the description of the fracture zones depends on the mesh-dependent length-scale parameter.
The main contribution of this paper is a design of novel level dependent objective functions for the phase field fracture problems, see Section \ref{section:coarse_level_models}.

Employing recursive multilevel trust region method, with trust region radius defined by the $\ell_{ \infty}$ norm,  allows for a natural treatment of the variable bounds \cite{conn2000trust}.
For this reason, the presented RMTR method  is set-up in such a way, that all corrections satisfy the irreversibility condition. 
This is achieved by performing coarse level minimization with respect to a restricted irreversibility condition, see Section \ref{section:rmtr}.

This paper is organized as follows: 
We start with a short review of phase-field fracture model, Section \ref{section:pf}. 
In  Section \ref{section:rmtr}, we provide an introduction to the trust region and multilevel trust region methods. 
 Section \ref{section:transfer}  discusses our choice of the transfer operators. 
 Section \ref{section:coarse_level_models} comments on how to build efficient level-dependent objective functions. 
The multilevel treatment of the irreversibility condition is addressed in Section \ref{section:constrains}. 
Numerical examples used to test the RMTR Algorithm will be introduced in Section \ref{section:numerical_examples}. 
Finally, Section \ref{section:numerical_results} demonstrates overall performance of the RMTR solution strategy. 
The presented work is summarized in Section \ref{section:conclusion}.

\section{Variational formulation of phase-field fracture}
\label{section:pf}
\begin{figure}
\subfloat[\label{fig:1a}]{
    \begin{tikzpicture}[scale=0.55]
      \draw[color=black,fill=gray!8,very thick] (0,0) -- (5,0) -- (5,5) -- (0,5) --cycle;
      \draw[color=black,ultra thick] (0,2.5)--(2.5,2.5);
      \draw [thick] (1.5,2) node[below] {$\Gamma$}; 
      \draw [thick] (5.1,2.5) node[right] {$\partial \Omega$} ; 
      \draw(4.5,4.5) node {$\Omega$}; 
    \end{tikzpicture}
}
    ~\hfill
\subfloat[\label{fig:1b}]{
    \begin{tikzpicture}[scale=0.55]    
      \draw[color=black,fill=gray!8,very thick] (0,0) -- (5,0) -- (5,5) -- (0,5) --cycle;

     	\draw[color=myred, <->] (0.95,2.9)--(0.95,2.5);
	\draw[color=myred] (0.85, 2.74) node[right] { \tiny ${l_s}$}; 
      
      \draw[color=black,ultra thick] (0,2.5)--(2.5,2.5);
      \draw[color=myred,dashed,thick] (0,2.1)--(2.5,2.1);
      \draw[color=myred,dashed,thick] (0,2.9)--(2.5,2.9);
      \draw[color=myred,dashed,thick] (2.5,2.1) arc (-90:90:0.4);
      \draw[thick,color=myred] (2.5,2) node[below] {$\Gamma_{l_s}$}; 
      \draw [thick] (5.1,2.5) node[right] {$\partial \Omega$} ; 
      \draw(4.5,4.5) node {$\Omega$}; 
    \end{tikzpicture}
}
    \hfill
\subfloat[\label{fig:1c}]{
\includegraphics[scale=0.21]{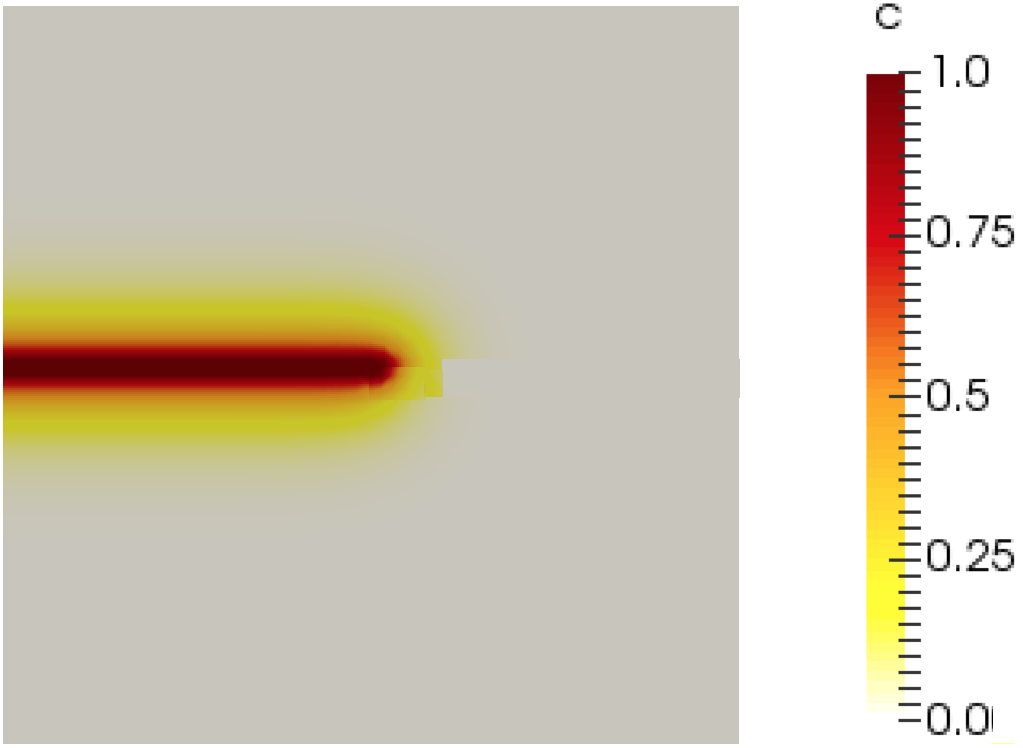}  
}
\caption{a) Computational domain $\Omega$ with internal discontinuity/fracture $\Gamma$. 
b) Phase field regularization of fracture surface, denoted by $\Gamma_{l_s}$.  Size of the regularized zone is controlled by the length-scale parameter $l_s$. 
c) Simulation realization of phase-field fracture problem. Damage parameter $c$ indicates fracture zone.  }
\label{fig:notation}
\end{figure}
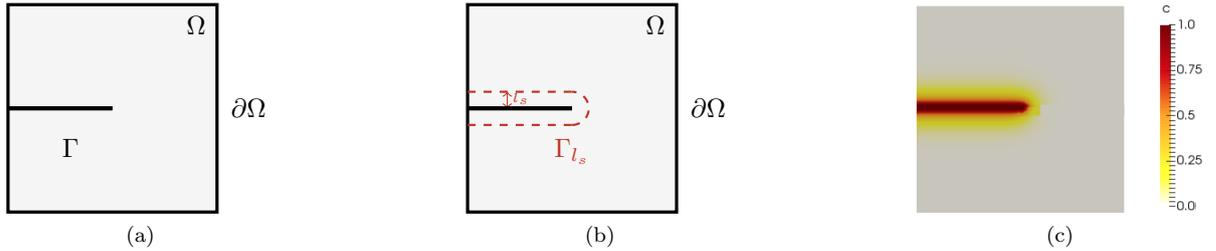
In this section, we briefly review the variational formulation of the second order phase-field model for brittle fracture. 
Our presentation focuses on quasi-static fracture problems and neglects inertia effects. 
A pseudo-time $t$ serves for indexing the deformation state of the loading process.

Let $\Omega \subset \mathbb{R}^d, d = 2, 3$ be the computational domain with internal fracture surface $\Gamma \subset \mathbb{R}^{d-1}$.
The computational domain $\Omega$ is assumed to have Lipschitz boundary $\partial \Omega$, which is decomposed into the Dirichlet $\partial \Omega_D$ and the Neumann $\partial \Omega_N = \partial \Omega \setminus \partial \Omega_D$ part. 
The variational approach to brittle fracture predicts the crack evolution by minimizing  the total potential energy of the body $\Omega$, given by
\begin{align}
\Psi(\u) :=  \underbrace{\int_{\Omega} \psi_e(\eps(\u))  \ d \Omega}_{\text{elastic energy}}  +\underbrace{ \int_{\Gamma} \pazocal{G}_c \ d \Gamma.}_{\text{fracture energy}}
\label{eq:potential_energy_basic}
\end{align}
 The energy functional \eqref{eq:potential_energy_basic} incorporates an elastic energy  and a fracture energy. 
In this work, the elastic energy describes small deformations and utilizes the linearized strain tensor $ \eps(\u) = 0.5 ( \nabla \u +  (\nabla \u)^T)$, where $\u$ is the displacement of a material point $\x\in\Omega$ at given pseudo-time $t$.

The term describing the fracture energy originates from the classical theory of brittle fracture \cite{griffith1921phenomena,  flugge2013elasticity}.
It states that the energy required to create a unit area of fracture surface equals the critical energy release $\pazocal{G}_c \in \R$. 
As a consequence, the fracture energy is defined as a contribution of $\pazocal{G}_c$ integrated over the fracture surface $\Gamma$. 
The presence of a fracture surface integral in \eqref{eq:potential_energy_basic} requires precise tracking of the crack surface paths and necessitates tedious remeshing, which leads into computationally exhausting algorithms.

\subsection{Energy formulation of phase-field approach to fracture}
\begin{figure}
\subfloat[\label{fig:2a}]{
\begin{tikzpicture}[domain=-2.0:2.0, samples=2000]
\draw[->, font=\small] (-2.,0) -- (2.,0) node[right] {$x$};
\draw[->, font=\small] (0,-0.1) -- (0,1.5) node[left] {$c(x)$};
\draw[-, font=\small] (-0.1, 1) -- (0.1, 1)node[right] {1};
\draw[-, font=\small] (-0.29, 0.05) -- (-0.29, -0.05) node [below] {-$l_s$};
\draw[-, font=\small] (0.29, 0.05) -- (0.29, -0.05) node [below] {$l_s$};
\draw[color=myred, style=thick] plot (\x, {exp(- (abs(\x)/0.3))} ) node[above] {};
\end{tikzpicture}
}
~\hfill
\subfloat[\label{fig:2b}]{
\begin{tikzpicture}[domain=-2.:2., samples=2000]
\draw[->, font=\small] (-2.,0) -- (2.,0) node[right] {$x$};
\draw[->, font=\small] (0,-0.1) -- (0,1.5) node[left] {$c(x)$};
\draw[-, font=\small] (-0.1, 1) -- (0.1, 1)node[right] {1};
\draw[-, font=\small] (-0.1, 0.1) -- (-0.1, -0.1);
\draw[-, font=\small] (0.1, 0.1) -- (0.1, -0.1);
\node at (-0.1, -0.3) {-$l_s$};
\node at (0.25, -0.3) {$l_s$};
\draw[color=myred, style=thick] plot (\x, {exp(- (abs(\x)/0.1))} ) node[above] {};
\end{tikzpicture}
}
~\hfill
\subfloat[\label{fig:2c}]{
\begin{tikzpicture}[domain=-2.:2., samples=2000]
\draw[->, font=\small] (-2.,0) -- (2.,0) node[right] {$x$};
\draw[->, font=\small] (0,-0.1) -- (0,1.5) node[left] {$c(x)$};
\draw[-, font=\small] (-0.1, 1) -- (0.1, 1)node[right] {1};
\draw[-, font=\small] (-0.01, 0.01) -- (-0.01, -0.01);
\draw[-, font=\small] (0.01, 0.01) -- (0.05, -0.05);
\node at (-0.1, -0.3) {-$l_s$};
\node at (0.25, -0.3) {$l_s$};
\draw[color=myred, style=thick] plot (\x, {exp(- (abs(\x)/0.01))} ) node[above] {};
\end{tikzpicture}
}
\caption{One dimensional phase-field approximation to the crack surface at $x=0$, see~\cite{miehe2010thermodynamically}. 
 As $l_s \rightarrow 0$, the crack approximation approaches the sharp crack surface.
 The crack approximations with parameter $l_s$ set to a) $0.3$  b) $0.1$  c) $0.05$. }
\label{figure:pf_approx_fun}
\end{figure}
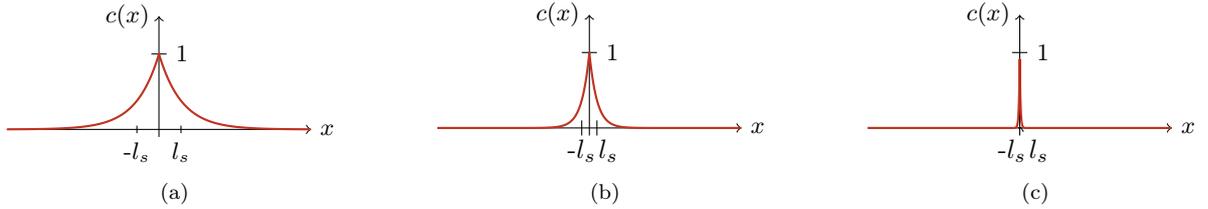

The phase field approach to fracture \cite{bourdin2008variational} overcomes difficulties related to the remeshing of the crack surface $\Gamma$
 by introducing the phase-field variable $c \in [0,1]$, which characterizes the material state from unbroken $c = 0$ to fully broken $c = 1$ and "damaged'' in between.
The regularization of the fracture is then performed by replacing fracture surface integral in \eqref{eq:potential_energy_basic} by its volumetric approximation
\begin{align}
\int_{\Gamma} \pazocal{G}_c \ d \Gamma  \approx \underbrace{ \int_{\Omega}  \frac{\pazocal{G}_c }{2 l_s } \ c^2  \ d \Omega}_{\text{reaction term}} + \underbrace{\int_{\Omega} \frac{\pazocal{G}_c l_s}{2} \ | \nabla c |^2 \ d \Omega.}_{\text{diffusion term}}
\label{eq:volumetric_approximation}
\end{align}
 Approximation \eqref{eq:volumetric_approximation} takes into account the phase-field $c$, its spatial gradient $\nabla c$ and regularization/length-scale parameter $l_s \in \R^+$.  
The reaction term in \eqref{eq:volumetric_approximation} corresponds to the area of the region, where material is damaged, e.g. $c \gg 0$. 
The diffusion term describes the area of the transition between the broken and intact part of the domain, where the gradient is high.
The length-scale parameter $l_s$ controls the thickness of the damaged region by weighting the contributions of the reaction and the diffusion terms.
Increasing $l_s$ makes the crack zone more diffused, while $l_s \rightarrow 0$ gives rise to a  sharp crack surface $\Gamma$. 
The effect of the length scale parameter $l_s$ is demonstrated in  Figure \ref{figure:pf_approx_fun}.

Taking into account \eqref{eq:volumetric_approximation}, the total potential energy from \eqref{eq:potential_energy_basic} can be reformulated as follows
\begin{align}
\Psi(\u, c) :=\underbrace{ \int_{\Omega} [g(c)(1-k) + k] \psi_e^+(\eps(\u)) + \psi_e^-(\eps(\u)) \  d \Omega}_{\text{elastic energy}} + \underbrace{\int_{\Omega} \pazocal{G}_c \Big[   \frac{1}{2 l_s }c^2 + \frac{l_s}{2} | \nabla c |^2 \Big] \ d \Omega,}_{\text{fracture energy}}
\label{eq:potential_energy_pf}
\end{align}
where $g(c) := (1-c)^2$ describes the degradation of the stored energy due to evolving damage. 
It should be noted, that the degradation function $g(c)$ acts only on the tension part of the elastic energy. 
This is due to the fact, that in the reality, fracture develops only in the direction of tension \cite{bourdin2000numerical}.  
The linear isotropic material models do not capture this feature and allow for fracture to occur in both tension and compression.
In order to avoid unrealistic crack patterns, the anisotropic material models have to be employed, see for example \cite{amor}, \cite{ambati2015review}.
We employ an anisotropic model introduced by Miehe et. al in \cite{miehe2010thermodynamically}, where $\psi_e(\eps(\u))$ is additively decomposed into $\psi_e^+(\eps(\u))$ and $\psi_e^-(\eps(\u))$
due to tension and compression, respectively.  For more details about employed energy split, we refer interested reader to \ref{section:split}. 

\begin{remark}
For simplicity, formulation \eqref{eq:potential_energy_pf} neglects body forces and surface tractions. 
Artificial parameter $k \approx \mathcal{O}(\epsilon)$ in  \eqref{eq:potential_energy_pf}  circumvents a full degradation of energy at the completely broken state and helps to improve the condition number of the numerical system \cite{ambrosio1990approximation}. 
\end{remark}

\subsection{Minimization problem}
\label{section:minimization_problem}
At each loading step $t$, the displacement $\u$ and phase field $c$ are characterized as a minimizer of the following minimization problem: \\
Find $(\u, c) \in \ \bm{V}^t \times H^1(\Omega) $, such that
\begin{equation}
\begin{aligned}
(\u, c) \in \ \text{argmin}  \ \  \Psi(\u, c) \\
\text{subject to} \ \  \partial_t c \geq 0,
\label{eq:min_problem_old}
\end{aligned}
\end{equation}
where $\bm{V}^{t} := \{ \u \in \bm{H}^1(\Omega) \ | \  \u = \bm{g}^t    \  \text{on} \ \  \partial \Omega_D  \}$  is an admissible space for the displacement field.
The symbol $H^1$ denotes the standard Sobolev space of weakly differentiable functions in $L^2$ with weak derivatives also in $L^2$. 
Since the minimization problem \eqref{eq:min_problem_old} is quasi-static, the time $t$ enters the above formulation only via the prescribed Dirichlet boundary condition $\bm{g}^t$ and the time-dependent irreversibility constraint $\partial_t c \geq 0$. 
This irreversibility constraint ensures positive evolution of the phase-field and prevents the crack from healing. Discretization of 
$\partial_t c \approx \frac{c^t - c^{t-1}}{\Delta t}, $
with $\Delta t$ denoting time-step size allows for the reformulation of the minimization problem $\eqref{eq:min_problem_old}$ as follows:\\
Find $\bm{U}= (\u, c) \in \bm{V}^t \times H^1(\Omega)$ such that
\begin{equation}
\begin{aligned}
\bm{U} \in  \ \text{argmin}  \ \  \Psi(\bm{U}) \\
\text{subject to} \ \  \bm{U} \geq \bm{U^{t-1}} \ \ \ \text{on} \  \Phi,
\label{eq:min_problem}
\end{aligned}
\end{equation}
where $\bm{U^{t-1}}$ represents solution on the previous time-step and $\Phi = \bm{0} \times H^1(\Omega) $. 
In this way, the inequality constraint acts only the phase-field variable $c$. 

\subsubsection{Properties of the minimization problem}
In order to solve the above minimization problem \eqref{eq:min_problem} numerically, we need to discretize underlying mathematical model. 
We perform spatial discretization of the problem by using finite-element method inside of the MOOSE framework \cite{gaston2009moose}. 
We refer interested reader to \ref{section:weak_form} for the detailed formulation of the Euler-Lagrange equations.

Solving energy minimization problem \eqref{eq:min_problem} is numerically challenging for the following reasons: 
\begin{itemize}
\itemsep-0.25em
\item The functional \eqref{eq:min_problem} is non-convex and therefore its minimization admits many local minimizers. 
    The most critical term is $ [g(c)(1-k) + k] \psi_e^+(\eps(\u))$ which is well-know to be non-convex in both variables $\u$ and $c$ \cite{wick2017modified}.
    The most popular solution strategy for solving nonlinear equations is Newton's method, as it is simple, robust and quadratically convergent \cite{nocedal2006numerical}.  
    However, a convergence properties of the standard Newton's method rely on a good initial guess. 
    In the context of phase-field fracture, solver failures are typically reported in the case of brutal crack evolution and in the post-peak regime \cite{GERASIMOV2016276}.    
    The convergence can be achieved by using small continuation steps \cite{alexander1978homotopy} or by applying some globalization strategy, such as trust region  \cite{conn2000trust}, line-search methods \cite{nocedal2006numerical}, cubic regularizations \cite{nesterov2006cubic} or affine invariant variants of the Newton method \cite{deuflhard2011newton}.
\item Resolving the regularized crack surface $\Gamma_{l_s}$ demands a fine mesh with a high spatial resolution. 
In fact, the mesh size $h$ has to fulfill $h \leq \frac{l_s}{2}$. Here, the $l_{s}$ is the length-scale parameter from \eqref{eq:volumetric_approximation} controlling the width of the diffusive crack zone.
As a consequence,  we end up with large-scale coupled nonlinear systems of equations, which are computationally demanding to solve. 
\item Arising systems of equations might also be ill-conditioned as the presence of the crack leads to the strongly varying elastic stiffness.
A suitable preconditioner, e.g. mutilevel \cite{mccormick1987multigrid}, or  domain decomposition methods \cite{saad2003iterative} can to be employed to address this issue. 
\item The irreversibility condition defines a lower bound on the phase-field variable $c$.  
Therefore, any applied solution strategy, for example line-search method, has to take into account point-wise constraint. 
The popular choice include for instance active set methods \cite{facchinei1998active}, \cite{ulbrich2011semismooth},  interior point methods \cite{byrd2000trust}, or sequential quadratic programming methods \cite{conn2000trust}. 
\end{itemize}
This work addresses all aforementioned points by employing the recursive multilevel trust region algorithm (RMTR)  \cite{gratton2010numerical, Gross2009}.
On one side,  the multilevel framework addresses issues of the ill-conditioning and the large number of dofs.  
On other side, the trust region framework tackles the non-convexity of the energy functional. 
In addition, the trust region method with the trust region radius defined by the $\ell_{ \infty}$ norm allows for a natural treatment of the point-wise constraints \cite{conn2000trust}, 
see the Section \ref{section:rmtr}.

\section{Recursive Multilevel Trust Region}
\label{section:rmtr}
In this work, we propose to solve the energy minimization problem \eqref{eq:min_problem} with the recursive multilevel trust region method (RMTR) originally introduced in  \cite{Gratton2008recursive}, \cite{Gross2009}. 
Extended details about the method and its convergence properties can be also found in \cite{gratton2010numerical} and \cite{gratton2008_inf}. 
The variant of the RMTR method considered here is associated with the following optimization problem
\begin{equation}
\begin{array}{ll}
\underset{\x \in \R^n}{\text{minimize}} &  f(\x) \\*[2ex]
\text{subject to} & \lb \leq \x. 
\end{array}
\label{eq:problem_min_tr}
\end{equation}
The objective function $f: \R^n \rightarrow \R$ in \eqref{eq:problem_min_tr} denotes the energy functional \eqref{eq:potential_energy_pf}. 
The solution vector $\x \in \R^n$ contains coefficients of the displacement and phase-field.
The lower bound $\lb \in \R^n$ originates from the crack irreversibility condition.
One should note, that the constraint in \eqref{eq:problem_min_tr} is defined component-wise, thus 
each component $(\x)_k$ of $\x$, where $k \in \{ 1, \dots, n \}$,  has to satisfy  $(\lb)_k \leq (\x)_k$.

\subsection{Trust region method}
\label{section:tr}
A trust region method (TR) is a step-size control strategy for minimizing nonlinear,  convex/non-convex objective functions \cite{conn2000trust}.
At each iteration $i$, the TR method approximates the objective function $f$ by a model $m_i$.  
The model $m_i$ is often constructed as a second order Taylor approximation of $f_i$ around the current iterate $\x_i \in \R^n$
\begin{align}
\label{eq:model}
m_i(\s_i ) = f_i + \langle \g_i, \s_i \rangle + \frac{1}{2} \langle \s_i,  \Hi \s_i \rangle, 
\end{align}
where we use the shortened notation $f_i = f(\x_i)$,  $\g_i = \nabla f(\x_i)$ and $\H_i = \nabla^2 f(\x_i)$. 
In this case, the update, i.e. the correction $\s_i$,  is then obtained by minimizing the model $m_i$ within a certain region. 
This region is called the trust region and it is defined around current iterate $\xxi$. 
The size of the region is defined by the means of radius $\Delta_i > 0$, which provides a constraint on the step size, e.g. $\| \s_i \|_p \leq \Delta_i$. 

The choice of norm $\| \cdot \|_p$ defines shape of the trust region. 
In this work, we employ $\ell_{\infty}$ norm, e.g. the trust region is defined as
\begin{align}
\pazocal{B}_i := \{ \x_i + \s \in \R^n  \ | \  \| \s \|_{\infty} \leq \Delta_i \}. 
\label{eq:b_box_tr}
\end{align}
As we will see later, defining trust region $\pazocal{B}_i$ as in \eqref{eq:b_box_tr} simplifies the treatment of point-wise constraint, e.g. irreversibility constraint. 
We can view \eqref{eq:b_box_tr} as point-wise constraints on step $\si$, thus $- \Delta_i \leq (\si)_k  \leq \Delta_i,      \ \ \forall  k \in \{ 1, \dots, n \}$.

In order to keep our iterates feasible throughout the computation, the minimization of model $m_i$ has to take place in the intersection of the feasible set 
$\pazocal{F}:= \{ \x   \in \R^n  \ | \  \lb \leq \x  \}$ with the trust region 
$\pazocal{B}_i$. For this reason, we define the active working set $\pazocal{W}_i$ as follows
\begin{align}
\pazocal{W}_i = \pazocal{B}_i \cap \pazocal{F}. 
\label{eq:intersection_BF}
\end{align}
Please note, that the  working set $\pazocal{W}_i$ needs to be recomputed in each iteration, as the trust region $\pazocal{B}_i$  is iteration dependent.
The evaluation of $\pazocal{W}_i$  is straightforward and can be performed component-wise.

At each iteration $i$, we obtain the step direction $\s_i$ by approximately solving the following constrained quadratic minimization problem
\begin{equation}
\begin{array}{ll}
\underset{\s \in \R^n}{\text{minimize}} &  m_i(\s) \\*[2ex]
\text{such that } & \x_i + \s \in  \pazocal{W}_i.   \\*[2ex]
\end{array}
\label{eq:tr_min}
\end{equation}
The iterate obtained in this way (also called the trial point), $\x_i + \s_i$, is not taken immediately.
Instead, the TR algorithm decides whether the trial point should be accepted and how to adjust the trust-region radius $\Delta_i$.  
To this aim, a trust-region ratio $\rho_i$, is defined as the quotient of actual reduction in $f_i$ and the reduction predicted by the model $m_i$. 
Given a step $\s_i$, $\rho_i$ is defined as
\begin{align}
\label{eq:rho}
   \rho_i = \frac{f(\x_i) - f(\x_i + \s_i)}{m_i(\0) - m_i(\s_i)} = \frac{\text{actual reduction}}{\text{predicted reduction}}.
\end{align}
The value of the trust region ratio $\rho_i$ close to 1 indicates good agreement between the
model $m_i$ and the function $f$. In this case, it is safe to accept the trial step and expand the trust region radius.
On the other hand, if value of $\rho_i$ is negative or close to 0, algorithm must reject trial step and shrink the trust region radius. 
Algorithm  \ref{alg:TR} summarizes the described process \cite{nocedal2006numerical}. 
\begin{algorithm}
\caption{Local\_TR}
\label{alg:TR}
\begin{algorithmic}[1]
\Require{$f:\R^{n} \rightarrow \R, \   {\x}_0 \in \R^{n},  \  \pazocal{F},  \ \Delta_0, \ \epsilon_g \in \R, \ i_{\max} \in \mathbb{N}$}
\Constants {  $\eta_1 \in \R$, where $0 < \eta_1 < 1$}
\For{$i=0, ..., i_{\max}$}
\State  $\pazocal{W}_i := \pazocal{B}_i \cap \pazocal{F}$  \Comment{Generate active working set}
\State
\State $\underset{\s_i \in \R^n}{\text{minimize}}  \   m_i(\si) = f_i + \langle \gi, \si \rangle + \frac{1}{2} \langle \si, \Hi \si \rangle$ \Comment{Solve constrained QP problem}
\State such that $\x_i + \s_i \in \pazocal{W}_i$
\State
\State  $ \rho_i = \frac{f(\xxi) - f(\xxi + \si)}{m_i(0) - m_i(\si)}$  \Comment{Trust region ratio computation}
\State
\If{$\rho_i > \eta_1$}
\State $\xxip := \xxi + \si$              \Comment{Trial point acceptance}
\Else
\State $\xxip := \xxi $			 \Comment{Trial point rejection}
\EndIf
\State
\State $\Delta_{i+1}$ = Radius\_update($\rho_i, \Delta_i$)   \Comment{Trust region radius update}
\State
\If{$ \pazocal{E}(\x_{i+1}, f) \leq \epsilon_{g}$} 
 \State \Return $\x_{i+1}$, $\Delta_{i+1}$  \Comment{Termination condition}
\EndIf
\EndFor
\State
\Return $\xxip$, $\Delta_{i+1}$
\end{algorithmic}
\end{algorithm}

\begin{minipage}{0.5\linewidth}
\begin{algorithm}[H]
\caption{Radius\_update}
\label{alg:multilevel_update}
\begin{algorithmic}
\Require{$ \rho_i, \ \Delta^{l}_{\nu, i} \in \R$}
\Constants {  $\eta_1, \eta_2, \gamma_1, \gamma_2 \in \R$, where \ $0 < \eta_1\leq \eta_2  < 1$ \textcolor{white}{.} \hspace{1.3cm}   and $0 < \gamma_1 < 1 < \gamma_2$}
\State
\State  $\Delta^{l}_{ i+1} = 
\begin{cases}
\gamma_1  \Delta^{l}_{\nu, i} & \rho_i < \eta_1 \\
\Delta_{l,i}  & \rho_i \in  [ \eta_1, \eta_2]   \\
 \gamma_2 \Delta^{l}_{\nu, i} & \rho_i >\eta_2
\end{cases}
$
\State
\State
\Return $\Delta^{l}_{ i+1}$
\end{algorithmic}
\end{algorithm}
\end{minipage}
~{}
\begin{minipage}{0.4\linewidth}
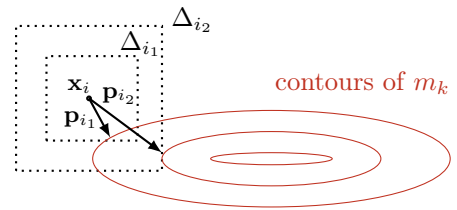
\begin{figure}[H]
\centering
	\begin{tikzpicture}[scale=0.8]
		\draw[color=black, dotted, thick] (-0.2,-0.2) rectangle (2.2,2.2);
		\draw[color=black, dotted, thick] (0.3,0.3) rectangle (1.8,1.7);

		\draw (2.7, 2.3) node {{\textcolor{black}{$\Delta_{i_2}$}}};
		\draw (1.85, 1.9) node {{\textcolor{black}{$\Delta_{i_1}$}}};

		\draw (5.5, 1.25) node {{\textcolor{myred}{contours of $m_k$}}};

		 \draw(4,0)[color=myred,style = thin] ellipse (2.94cm and 0.8 cm);
		\draw(4,0)[color=myred,style = thin] ellipse (1.8cm and 0.45 cm);
		\draw(4,0)[color=myred, style = thin] ellipse (1cm and 0.1 cm);

		\draw (0.84, 1.2) node {\small $\x_{i}$};
		 \draw(1,1)[thick] circle (0.03cm);

		\draw[color=black, thick, -latex] (1,1)--(2.2,0.1);
		\draw[color=black, thick, -latex] (1,1)--(1.35, 0.35);

		\draw (0.87, 0.65) node {\small {$\s_{i_1}$}};
		\draw (1.5, 1.05) node {\small {$\s_{i_2}$}};
	\end{tikzpicture}
\caption{
Step direction $\s_i$ is found as a minimizer of a quadratic model within the trust region radius $\Delta_i$.  }
\label{fig:delta_depict}
\end{figure}
\end{minipage}

\subsection{Nonlinear multilevel framework}
\label{section:multilevel}
Multilevel solution strategies take into account a hierarchy $\{ \pazocal{X}^l \}_{l=0}^L$
 of usually nested finite element, or approximation spaces. 
 Here, we use the subscript $l$ to denote a given level and assume $L \geq 1$ to be the finest level.
 In the case of geometric multilevel methods, the hierarchy $\{ \pazocal{X}^l \}_{l=0}^L$, with $\text{dim}( \pazocal{X}^l ) > \text{dim}( \pazocal{X}^{l-1} )$ , is constructed from a hierarchy of meshes $\{ \pazocal{T}^l \}_{l=0}^L$. 
For purpose of this work, the hierarchy of meshes $\{ \pazocal{T}^l \}_{l=0}^L$ was obtained by uniformly refining an initial coarse mesh $\pazocal{T}^{0}$.

Transfer of the data between subsequent levels of multilevel hierarchy is ensured by transfer operators. 
In the presented work, we use three types of the transfer operators:
\begin{itemize}
\itemsep-0.25em
\item The prolongation operator $\I_{l}^{l+1}: {\pazocal{X}}^{{l}} \rightarrow  {\pazocal{X}}^{{l+1}}$ 
is designed to transfer primal variables, such as  corrections from a coarse level to a finer level. 
\item The restriction operator   $\Rev^l_{l+1} : \pazocal{X}^{{l+1}} \rightarrow \pazocal{X}^{{l}}$ transfers dual variables, e.g. gradients, from a fine level to a coarser level.
In this work, the operators $\Rev^l_{l+1}$ and $\I_{l}^{l+1} $ are related by relation  $\Rev^l_{l+1} = (\I_{l}^{l+1})^T$. 
\item The projection operator $\P^l_{l+1} : \pazocal{X}^{{l+1}} \rightarrow \pazocal{X}^{{l}}$ transfers primal variables, for example the current iterate, from a fine level to a coarse level. 
\end{itemize}
In our approach, all three transfer operators are computed based on pseudo-L$^2$-projections, see Section \ref{section:transfer}.

The multilevel solution process is built as an interplay between smoothing and coarse grid correction. 
The smoothing step reduces high-frequency error related to each level.
In the context of the trust region framework considered here, smoothing is associated with performing a small number of trust region iterations, where the arising quadratic models \eqref{eq:model} are solved approximately.  
On the coarsest level, the minimization of coarse level objective function is again achieved by the trust region method, but constrained quadratic subproblems are solved accurately.  

Minimization on each level $l$ is associated with level-dependent objective function $\h^l$.
The definition of the level-dependent objective functions $\{ \h^l \}_{l=0}^{L-1}$ is often based on the coarse level representation/discretization of the original model, denoted by $\{ f^l \}_{l=0}^{L}$. 
A detailed discussion on how to construct  the level-dependent objective functions $\{ \h^l \}_{l=0}^{L-1}$ will be provided in Section \ref{section:coarse_level_models}. 

\subsection{RMTR algorithm}
\label{section:rmtr_algo}
Our description of the RMTR, Algorithm~\ref{alg:rmtr}, has the form of V-cycle and closely follows \cite{Gross2009}. 
In the following, we use a triplet (cycle, level, iterate) to denote all computational quantities. 
For example the solution computed as the $i$-th iterate during the $\nu$-th cycle on level $l$ is denoted by $\x^{l}_{\nu, i}$. 

Each V-cycle consists of nonlinear smoothing steps and a coarse level solve. 
Starting from the finest level, $l=L$, the RMTR algorithm performs $\mu_1 \in \N$  pre-smoothing steps. 
This is done by a trust region solver, see Algorithm~\ref{alg:TR}. 
The result of pre-smoothing is an intermediate iterate $\x^{l}_{\nu, \mu_1}$. 
The projection of this iterate to the next coarser level yields the initial guess on level $l-1$, i.e. 
\begin{align}
\x^{l-1}_{\nu, 0} =  \P^{l-1}_{l} \x^{l}_{\nu, \mu_1}.
\label{eq:init_iterate_coarse}
\end{align}
Once we have reached the coarser level, we can perform the minimization of the level dependent objective function $h^{l-1}_{\nu}$. 
The result of this minimization, a coarse level correction, $\s_{\nu}^{l-1} = \x^{l-1}_{\nu, \mu_1} - \x^{l-1}_{\nu, 0}$, is afterwards prolongated back to the finer level.
As a consequence, we have to guarantee that the the prolongated coarse level correction satisfies the following requirements:
\begin{enumerate}
\setlength\itemsep{0em}
\item   It does not exceed the current fine level trust region radius, e.g.
$\| \I_{l-1}^{l} \s^{l-1}_{\nu}  \|_{\infty} \leq \Delta^{l}_{\nu, \mu_1}.$
\item  It does not violate the fine level irreversibility condition, e.g.
$\I_{l-1}^{l} \s^{l-1}_{\nu}  \in \pazocal{F}^l_{\nu}.$
\end{enumerate}
We can meet both requirements  by minimizing the coarse level objective function within a feasible set $\pazocal{L}_{\nu}^l$
 as follows
\begin{equation}
\begin{array}{ll}
\underset{\s^{l-1}_{\nu} \in \R^{n^{l-1}}}{\text{minimize}} &  \h_{\nu}^{l-1}(\x_{\nu, 0}^{l-1} + \s^{l-1}_{\nu}) \\*[2ex]
\text{such that} &  \x_{\nu, 0}^{l-1} + \s^{l-1}_{\nu} \in \pazocal{L}_{\nu}^l.
\end{array}
\label{eq:coarse_level_min}
\end{equation}
We will discuss in Section \ref{section:constrains}  how to construct set $\pazocal{L}_{\nu}^l$, such that $\| \I_{l-1}^{l} \s^{l-1}_{\nu}  \|_{\infty} \leq \Delta^{l}_{\nu, \mu_1}$ and $\I_{l-1}^{l} \s^{l-1}_{\nu}  \in \pazocal{F}^l_{\nu}$. 

The minimization of $\h^l_{\nu}$ on a given level $l$ is performed by means of the trust region, Algorithm \ref{alg:TR}. The correction $\s_{\nu,i}^l$ is obtained by minimizing the quadratic model $m$ as
\begin{align}
\underset{\s_{\nu,i}^l \in \R^{n^l}}{\text{minimize}}  & \ \ \h_{\nu,i}^l + \langle \g_{\nu,i}^l, \s_{\nu,i}^l \rangle + \frac{1}{2} \langle \s_{\nu,i}^l,  \H_{\nu,i}^l  \ \s_{\nu,i}^l \rangle \\
\text{such that} &  \ \  \x_{\nu,i}^l + \s_{\nu,i}^l \in \pazocal{W}_{\nu,i}^l, 
\end{align}
where $\pazocal{W}_{\nu,i}^l$ is the working active set. 
As for the single level case, the set $\pazocal{W}_{\nu,i}^l$ is result of the intersection of the feasibility set $\pazocal{L}_{\nu}^l$ with the current trust region  $\pazocal{B}_{c, i}^l:= \{ \x_{\nu,i}^l + \s_{\nu,i}^l \in \R^{n^l}  \ | \  \| \s_{\nu,i}^l \|_{\infty} \leq \Delta_{\nu,i}^l \}$, e.g. $\pazocal{W}_{\nu,i}^l:= \pazocal{L}_{\nu}^l \cap \pazocal{B}_{c, i}^l$.

\subsubsection{Acceptance of the coarse level correction}
\label{section:trial_point_acceptance}
After the approximate minimization on level $(l-1)$ is performed, we prolongate the resulting coarse level correction back  to the finer level, i.e.
\begin{align}
\s^{l}_{\nu,\mu_1 +1} = \I_{l-1}^{l}(\x^{l-1}_{c, \mu^{l-1}} - \x^{l-1}_{\nu, 0}). 
\label{eq:coarse_level_cor}
\end{align}
The symbol $\mu^{l-1}$ in \eqref{eq:coarse_level_cor} denotes the number of all computed corrections on level $l-1$ during cycle $\nu$.
The coarse level correction is not accepted immediately, but only if it provides a decrease in the fine level objective function; thus only if
\begin{align}
 \h_{\nu, \mu_1 +1}^{l}  <  \h_{\nu, \mu_1}^{l}.  
\end{align}
The quality of the prolongated coarse level correction $\s^{l}_{\nu,\mu_1+1}$ is determined by means of the multilevel trust region ratio $\rho^{l}_{\nu,\mu_1}$, 
which measures the agreement between the fine level and the coarse level reduction
\begin{align*}
\rho^{l}_{\nu,\mu_1+1}= \frac{\h^{l}_{\nu}(\x^{l}_{\nu,\mu_1}) - \h^{l}_{\nu}(\x^{l}_{\nu,\mu_1} + \s_{\nu,\mu_1}^{l})}{\h^{l-1}_{\nu}( \P^{l-1}_{l} \x^{l}_{\nu,\mu_1}) - \h^{l-1}_{\nu+1}(\x^{l-1}_{\nu, \mu^{l-1}})}  = \frac{\text{fine level reduction}}{\text{coarse level reduction}}. 
\end{align*}
The correction $\s_{\nu, \mu_1+1}^{l}$ will be accepted only if $\rho_{\nu, \mu_1+1}^{l} > \eta_1$.
This is in contrast to the classical nonlinear multilevel schemes, e.g. FAS \cite{brandt1977multi},  or NMG \cite{hackbusch2013multi}, where the coarse level correction
is used independently of its quality.  A similar strategy has been used in \cite{RKornhuber_RKrause_2001} for controlling global convergence. 

Once the correction $\s_{\nu, \mu_1+1}^{l}$ has been accepted or rejected, the RMTR algorithm adjusts the trust region radius $\Delta_{\nu}^l$ accordingly, see Algorithm \ref{alg:multilevel_update}. 
The V-cycle completes with $\mu_2 \in \N$ post-smoothing steps.
The outlined process is performed recursively on each level, except on the coarsest one, where no recursion is invoked. 
On the finest level, the last iterate of cycle $\nu$ is taken over to the next cycle, $(\nu+1)$, as the initial iterate. 

Termination of the RMTR algorithm is specified in terms of the criticality measure  $\pazocal{E}$. 
Here, we follow \cite{conn2000trust} and define 
$\pazocal{E}(\x, f) := \| \mathcal{P}(\x - \nabla f(\x) ) - \x \|,$
where $\mathcal{P}$ is an orthogonal projection onto the feasible set. The details on how to compute $\mathcal{P}$ can be found in \ref{section:projection_set}. 

\begin{algorithm}[ht]
\caption{RMTR}
\label{alg:rmtr}
\begin{algorithmic}[1]
\Require{$ l \in \N, \x^l_{\nu, 0},  \g^{l+1}_{\nu, \mu_1} \in \R^n, \H^{l+1}_{\nu, \mu_1} \in \R^{n \times n},  \Delta^l_{\nu,0}, \epsilon^l_g \in \R $}
\Constants { $\mu_1, \mu_2 \in \N$,  $\eta_1 \in  \R$, where $0 < \eta_1  < 1$}
\State
\State Generate $\h_{\nu}^l$ by means of \eqref{eq:coarse_objective_first_order}, \eqref{eq:coarse_objective_galerkin}, \eqref{eq:coarse_objective_second_order}, or \eqref{eq:coarse_objective_modified2nd_order} 
\Comment{Initialization of the objective function \label{start}}
\State Generate $\pazocal{L}_{\nu}^l$ by means of \eqref{eq:define_l}
\Comment{Initialization of the feasible set}
\State
\State $ [\x_{\nu, \mu_1}^l, \  \Delta_{\nu, \mu_1}^l] = \text{ Local}\_\text{TR}(\h_{\nu}^l,  \ \x_{\nu, 0}^l, \  \pazocal{L}_{\nu}^l,  \  \Delta_{\nu, 0}^l, \ \epsilon^l_g ,  \  \mu_1 )$
\Comment{Pre-smoothing}
\State
\If{$l > 0$}
\State $[\x_{\nu, \mu_{l-1}}^{l-1}]$=RMTR($l-1,  \P^{l-1}_{l} \x_{\nu,\mu_1}^l, \Rev^{l-1}_{l} \g_{\nu,\mu_1}^l, \Rev^{l-1}_{l} \H^{l}_{\nu,\mu_1} \I^{l-1}_{l},   \Delta_{\nu,\mu_1}^l, \epsilon^{l-1}_{g}) $ 
\State $\s_{\nu,\mu_1+1}^l = \I_{l-1}^{l}(\x_{\nu, \mu_{l-1}}^{l-1} - \x_{\nu, 0}^{l-1})$ \Comment{Prolongation of coarse level correction}
\State 
\State $\rho^l_{\nu, \mu_1 +1} =  \frac{\h^l_{\nu}(\x^l_{\nu,\mu_1}) - \h^l_{\nu}(\x^l_{\nu,\mu_1} + \s_{\nu,\mu_1+1}^l)}{\h_{\nu}^{l-1}( \P^{l-1}_{l} \x^l_{\nu,\mu_1}) - \h^{l-1}_{\nu}(\x^{l-1}_{c, \mu^{l-1}})}\mu^l$  \Comment{Computation of multilevel TR ratio}
\State
\If{$\rho^l_{\nu, \mu_1+1} > \eta_1$}
\State $\x_{\nu,\mu_1 +1}^l := \x_{\nu,\mu_1}^l + \s_{\nu, \mu_1 +1}^l$ \Comment{Acceptance of coarse level correction}
\Else 
\State $\x_{\nu,\mu_1 +1}^l := \x_{\nu,\mu_1}^l $ \Comment{Rejection of coarse level correction}
\EndIf
\EndIf 
\State
\State $[\Delta_{\nu, \mu_1 +1}^l]$ = Radius\_update($l, \  \rho_{\nu,\mu_1 +1}^l, \ \Delta_{\nu, \mu_1}^l, \ \Delta^{l+1}$)  
\State 
\If{$ \pazocal{E}(\x_{\nu,\mu_1 +1}^l, h_{\nu}^l) \leq \epsilon^l_{g}$} 
\State \Return $\x_{\nu, \mu_1+1}^l$  \Comment{Termination condition}
\EndIf
\State
\State $ [\x_{\nu, \mu^l}^l,  \ \Delta_{\nu, \mu^l}^l] =\text{ Local}\_\text{TR}(\h_{\nu}^l,  \x_{\nu, \mu_{1} +1}^l, \  \pazocal{L}_{\nu}^l, \ \Delta_{\nu, \mu_{1} +1}^l, \ \epsilon^l_g,  \  \mu_{2} )$ \Comment{Post-smoothing}
\State
\If{$l=L$}
\State Set $ \Delta_{\nu+1, 0}^{l} = \Delta_{\nu, \mu^l}^{l}, \ \x_{\nu+1,0}^{l}= \x_{\nu, \mu^l}^{l}, \ \nu = \nu + 1 $ \Comment{Initialization of next V-cycle}
\State \Goto{start}  \Comment{Go to next V-cycle}
\Else
\State
\Return $\x_{\nu, \mu^l}^l$ \Comment{Continue recursion and go to next finer level}
\EndIf
\end{algorithmic}
\end{algorithm}

\section{Level dependent objective functions}
\label{section:coarse_level_models}
On each level $l$, the RMTR Algorithm \ref{alg:rmtr} minimizes some nonlinear objective function $\h_{\nu}^l$. 
The result of this minimization, the coarse level correction, is expected to provide decrease in the objective function on the subsequent finer level.
The overall efficiency of the multilevel method relies on the quality of this prolongated coarse level correction. 
Our aim is to find coarse level model (level dependent objective function), which yields adequate coarse level correction. 
In this section, we review several models, which are often used in literature for constructing $\{ \h_{\nu}^l \}_{l=0}^L$. 
Our discussion mostly focuses on the convergence properties and on the computational complexity of the given models.
Later in this section, we propose an alternative model, which is specially designed for the phase-field fracture models considered in this work. 
All models discussed in this section assume $\h^L := f^L$.

\subsection{First order consistency}
\label{section:fisrt_order}
Classical nonlinear multilevel schemes, such as FAS~\cite{brandt1977multi}, NMG~\cite{hackbusch2013multi},  NMM~\cite{krause2009nonsmooth} or MG/OPT~\cite{Nash2000multigrid}, often employ 
a first order consistency model, where the level dependent function  $\h^l: \R^{n^l} \rightarrow \R$  is defined as
\begin{align}
\h^{l}_{\nu}(\x_{\nu}^{l}) := \underbrace{ f_{\nu}^{l}(\x_{\nu}^{l})}_{\text{coarse model}} + 
\underbrace{ \langle \delta \g^l_{\nu}, \x^l_{\nu} - \x^{l}_{\nu,0} \rangle,}_{\text{1st order coupling}}
\label{eq:coarse_objective_first_order}
\end{align}
where 
\begin{align}
\delta \g^l_{\nu} := 
\begin{cases}
\Rev^l_{l+1} \g_{\nu,\mu_1}^{l+1}  - \nabla f_{\nu}^{l}(\x^{l}_{\nu,0}) &  \text{if}  \ \ \ l < L\\
0  & \text{if}  \ \ \ l = L.
\end{cases}
\label{eq:delta_g}
\end{align}
Here, the first term, $f_{\nu}^{l}(\x_{\nu}^{l})$, is the original  model problem discretized on level $l$, thus $f_{\nu}^{l}$ is the energy functional $\Psi(\u, c)$ evaluated at the level $l$. 
The second term, $\langle \delta \g^l_{\nu}, \x^l_v - \x^{l}_{\nu,0} \rangle$, couples the coarse level with the fine level. 
This is achieved by using the quantity $\delta \g^l_{\nu}$, which expresses the difference between the restricted fine level gradient and the initial coarse level gradient (evaluated at $\x^{l}_{\nu, 0} =  \P^{l}_{l+1} \x^{l+1}_{\nu, \mu_1}$).
The presence of this term enforces the relation  $\g^{l-1}_{\nu,0} = \nabla \h_{\nu,0}^{l-1}(\x^{l-1}_{\nu, 0}) = \Rev^{l-1}_{l} \g_{\nu, \mu_1}^{l}$.
We refer to the model \eqref{eq:coarse_objective_first_order} as first order consistency since the first step of the coarse level minimization process will be done in the direction of the negative restricted fine level gradient. 
   
Model \eqref{eq:coarse_objective_first_order} is known to yield fast convergence, even with a poor initial guess, under the assumption that 
$h_{\nu}^{l}$ provides a good approximation of the error after pre-smoothing. 
 If the coarse level model $h_{\nu}^{l}$ has poor approximation properties, then the multilevel method might produce an inadequate coarse level correction. 
 This results in a slow convergence rate, even if the initial error is very low and nonlinearity is weak  \cite{Yavneh2006}.
 
\begin{figure}
\subfloat[ $l=1$ \label{fig:f1}]{
\includegraphics[scale=0.035]{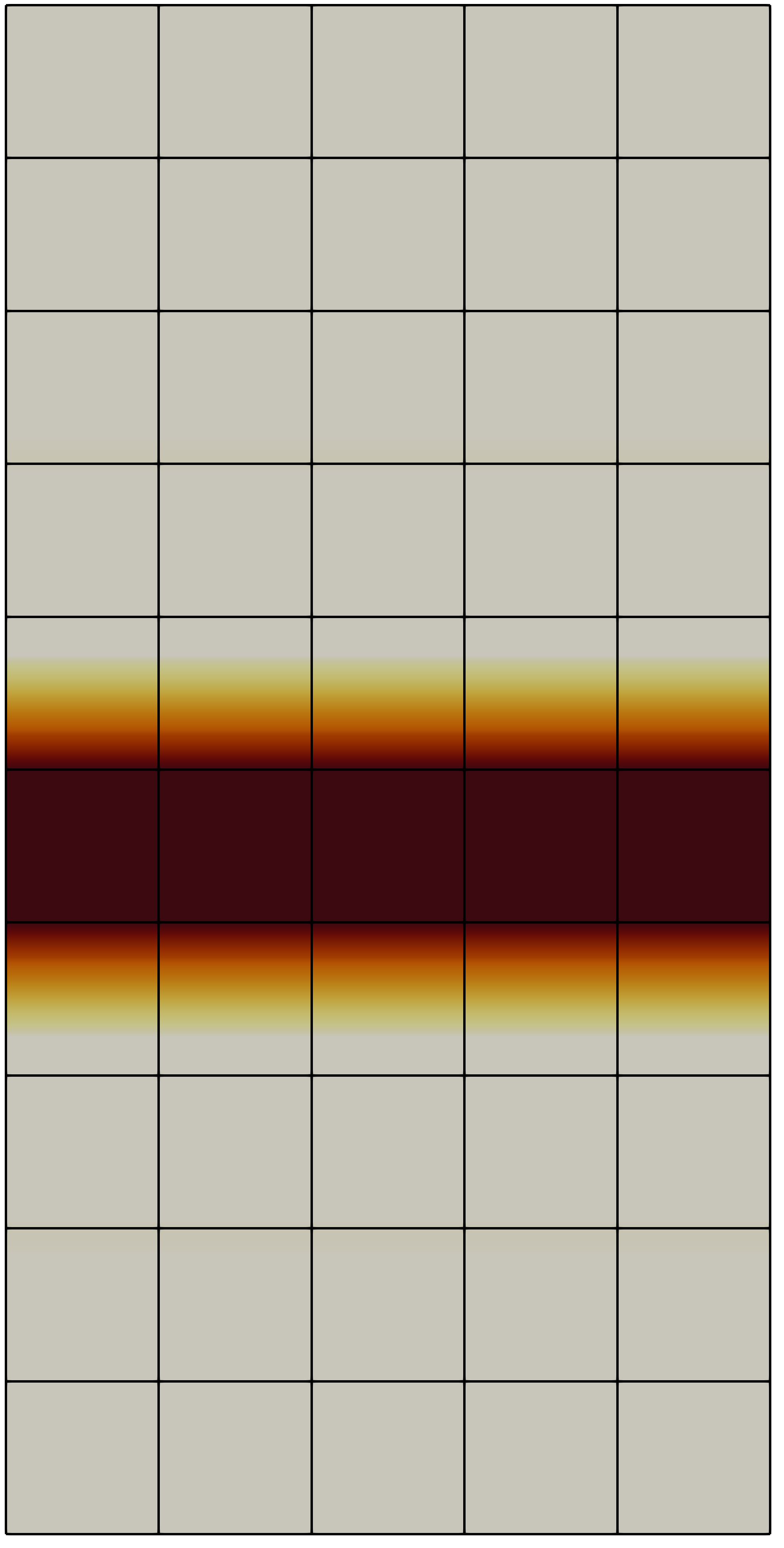}  
}
\hfill
\subfloat[$l=2$ \label{fig:f2}]{
\includegraphics[scale=0.035]{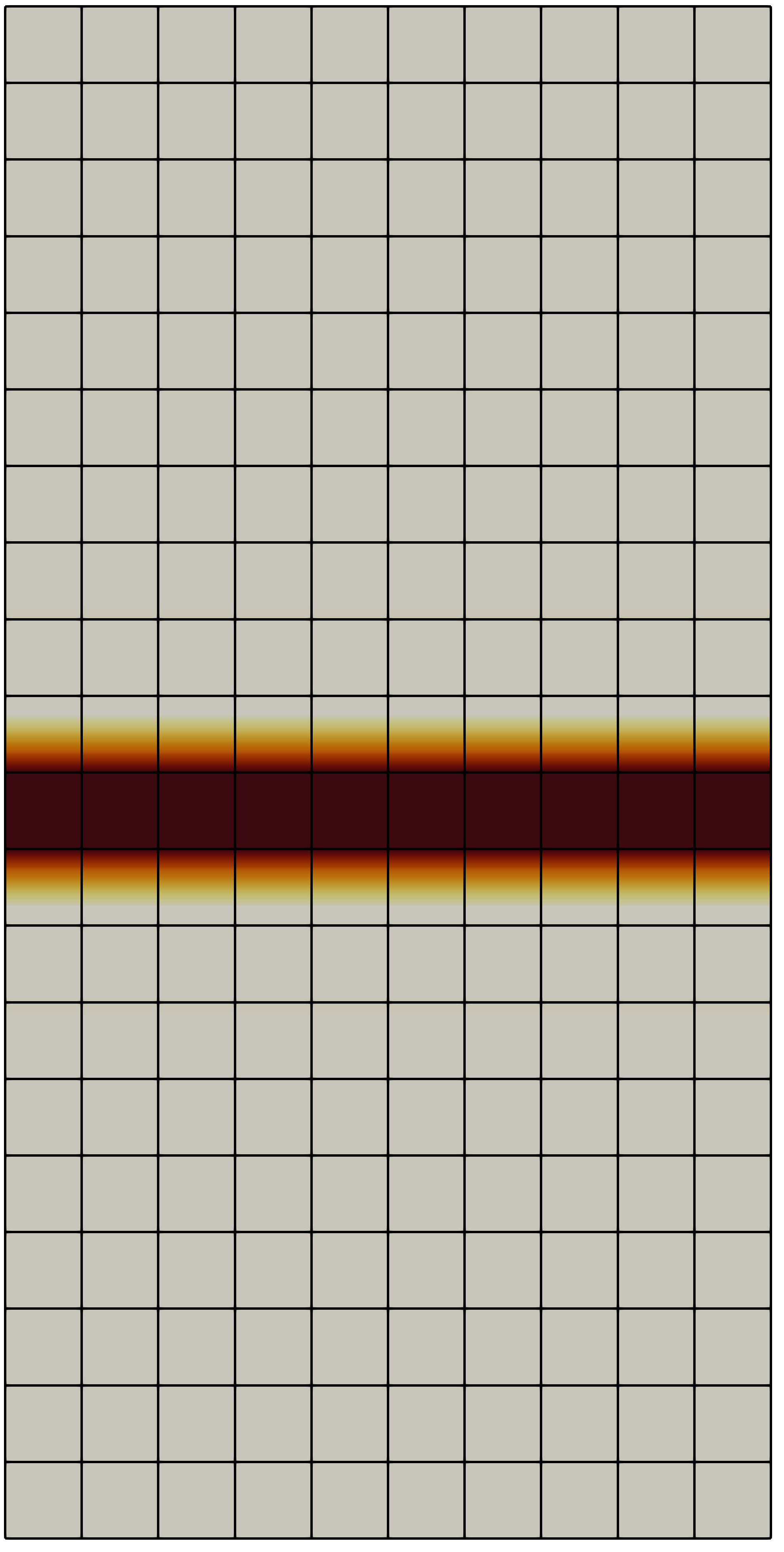}  
}
\hfill
\subfloat[ $l=3$ \label{fig:f3}]{
\includegraphics[scale=0.035]{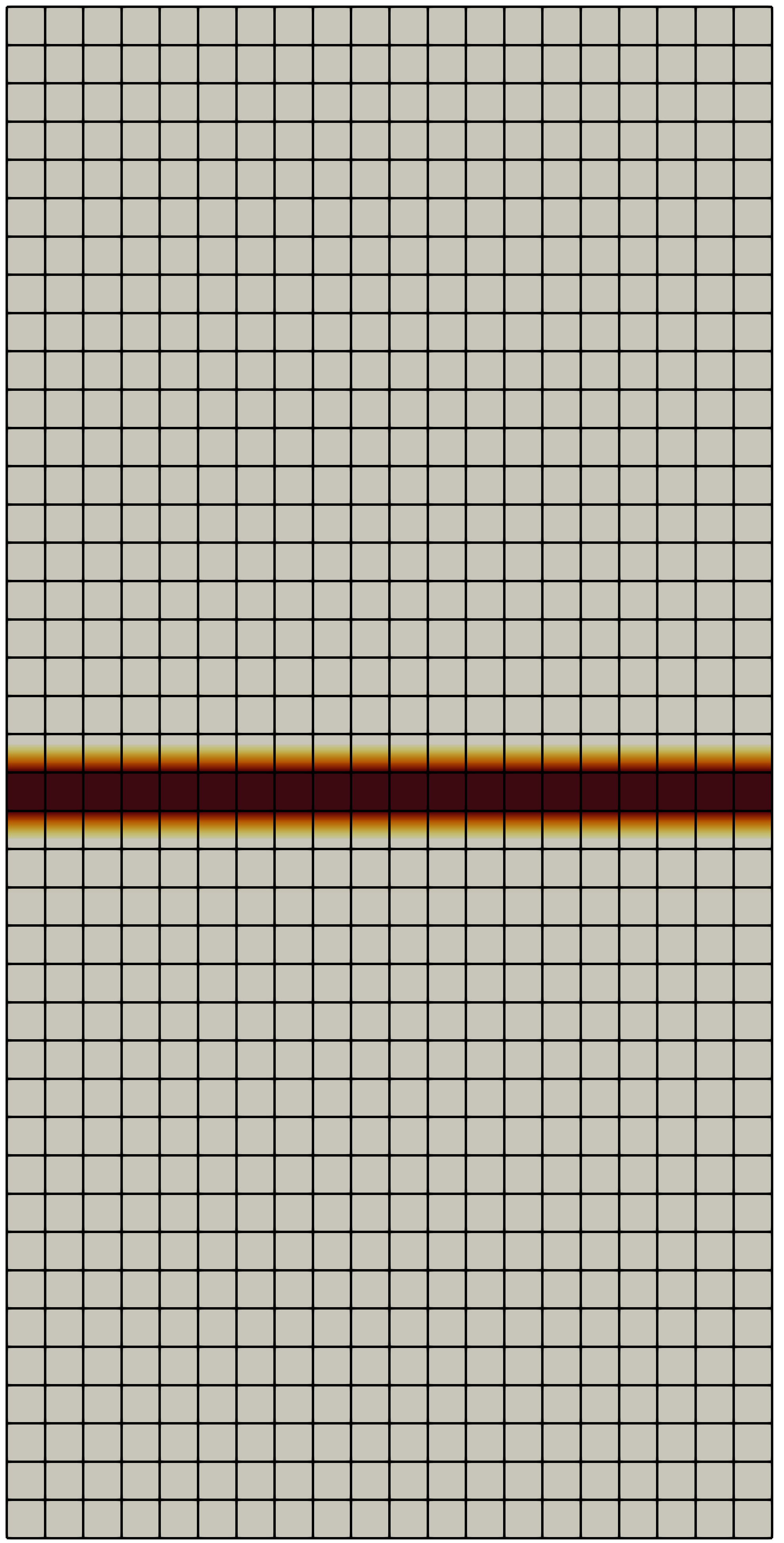}  
}
\hfill
\subfloat[$l=4$ \label{fig:f4}]{
\includegraphics[scale=0.035]{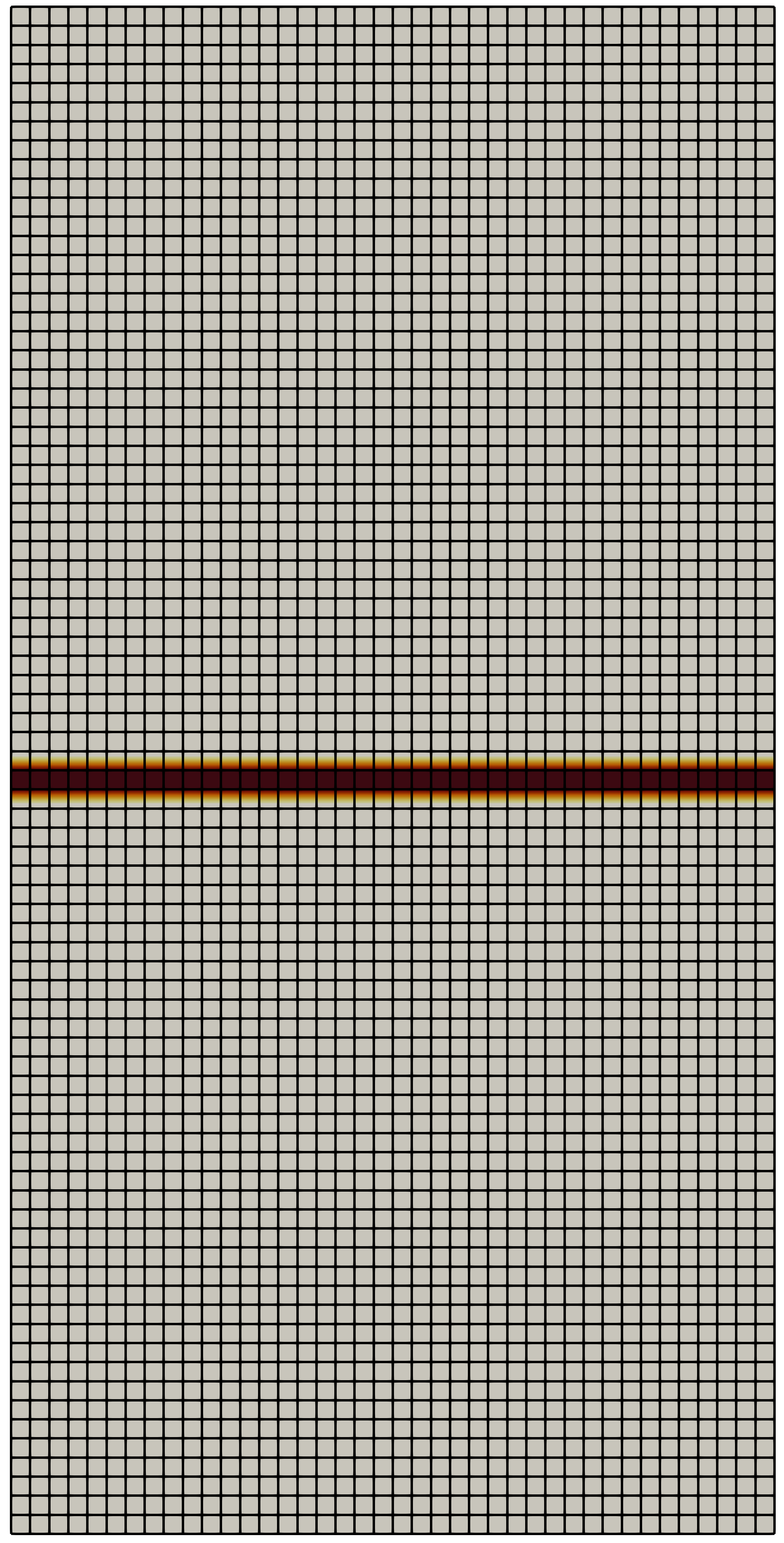}  
}
\hfill
\subfloat{
\includegraphics[scale=0.052]{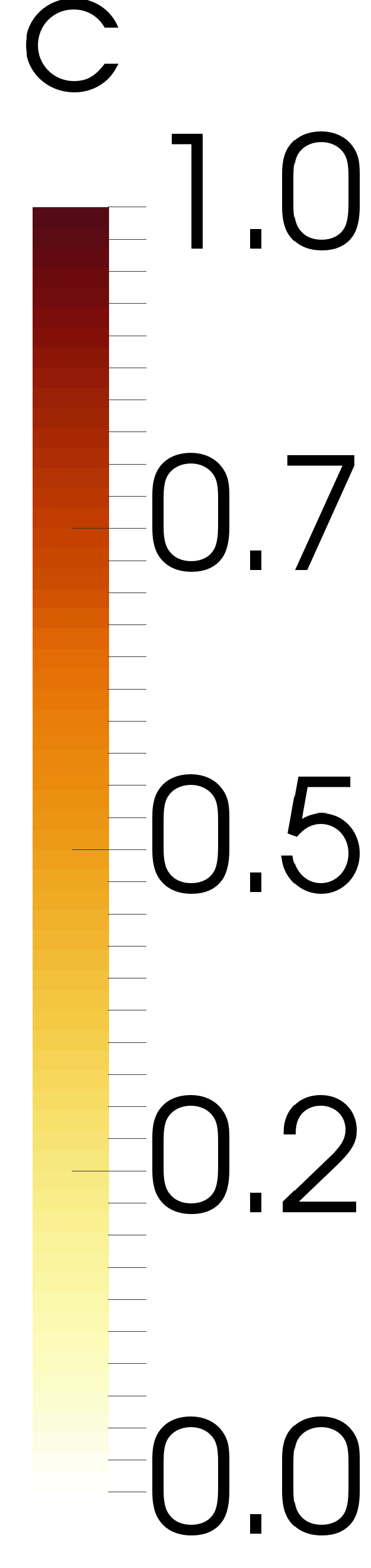}
}
\caption{ Two dimensional tension test. 
The solution, crack pattern, depends on the mesh size and reassembles sharp fracture surface as $h \rightarrow 0$. }
\label{fig:coarse_representations_of_fracture}
\end{figure}

In terms of the phase-field fracture problems considered here, the model $h^l_{\nu}$ based on the $f^l_{\nu}$, the coarse level discretization of the energy functional  $\Psi(\u, c)$, cannot be expected to provide good coarse level correction.
This is due to the fact that volumetric approximation to the fracture surface integral \eqref{eq:volumetric_approximation} is determined by the mesh dependent length scale parameter $l_s$.
Therefore, we are not able to represent the same fracture zones on the coarse level, as we can on the finest level. 
Figure \ref{fig:coarse_representations_of_fracture} demonstrates this effect on two dimensional tension test. 
In particular, we consider a rectangular computational domain discretized with quadrilateral finite elements. 
The body is subjected to tensile loading, which results in a straight crack path inside the specimen. 
We perform the same experiment repeatedly while increasing the resolution of the mesh. 
In the same time, we also modify the length-scale parameter $l_s$ as follows $l_s = 2h$, where $h$ represents mesh size.
As we can see, the obtained crack paths differ for different resolution levels. 
Moreover, even larger inconsistencies between different refinement levels can be expected for more complex loading scenarios. 

\subsection{Galerkin approximation}
Exploiting ideas from linear multilevel methods, the coarse level model $\h^l: \R^{n^l} \rightarrow \R$ can be created as a restricted version of the 
fine level model
\begin{align}
\h^{l}_{\nu}(\x_{\nu}^{l}) :=  \langle \Rev^{l}_{l+1} \g_{\nu, \mu_1}^{l+1}, \x^l_{\nu} - \x^{l}_{\nu,0} \rangle + \frac{1}{2} \langle  \x^l_{\nu} - \x^{l}_{\nu,0}, ( \Rev^l_{l+1} \H_{\nu, \mu_1}^{l+1} \I_l^{l+1})  (\x^l_{\nu} - \x^{l}_{\nu,0} )\rangle,
\label{eq:coarse_objective_galerkin}
\end{align}
where $\Rev^{l}_{l+1} \g_{\nu, \mu_1}^{l+1}$ represents the restricted fine level gradient and $\Rev^l_{l+1} \H_{\nu, \mu_1}^{l+1} \I_l^{l+1}$ is the restricted fine level Hessian.
In the unconstrained convex case, the formulation \eqref{eq:coarse_objective_galerkin} is equivalent to applying a Galerkin linear multigrid on the associated Newton's equation.

The overall performance of the method is determined by the accuracy of the fine level linearization, which depends on the properties of the Hessian and on the given initial guess. 
If the initial guess is in the local neighborhood of the solution, the method can be expected to achieve the quadratic convergence rate. 
However, for generic non-linear, non-convex functions, the 1st order consistency formulation \eqref{eq:coarse_objective_first_order}, which tackles nonlinearity directly on the coarse level, yields usually better convergence rates \cite{Yavneh2006}. 
From the computational point of view, the model \eqref{eq:coarse_objective_galerkin} is simple to implement, as the multilevel solver does not require knowledge about coarse level discretization.    
Concerning the phase-field fracture simulations, the coarse level model \eqref{eq:coarse_objective_galerkin} is able to represent fracture zones well, as it is built from the fine level data.

\subsection{Second order consistency}
Yavneh and Dardyk in \cite{Yavneh2006} explore an idea of combining the first order consistency, e.g. \eqref{eq:coarse_objective_first_order}, with the Galerkin approximation approach \eqref{eq:coarse_objective_galerkin} by
defining $\h^l: \R^{n^l} \rightarrow \R$ as follows
\begin{align}
\h^{l}_{\nu}(\x_{\nu}^{l}) := \underbrace{f_{\nu}^{l}(\x_{\nu}^{l})}_{\text{coarse model}} + 
\underbrace{ \langle \delta \g^l_{\nu}, \x^l_{\nu} - \x^{l}_{\nu,0} \rangle }_{\text{1st order coupling}} + 
\underbrace{\frac{1}{2} \langle  \x^l_{\nu} - \x^{l}_{\nu,0}, \delta \H  (\x^l_{\nu} - \x^{l}_{\nu,0} )\rangle.}_{\text{2nd order coupling}} 
\label{eq:coarse_objective_second_order}
\end{align}
Here, the term $\delta \g^l_{\nu}$  is as defined in \eqref{eq:delta_g}. 
The quantity $\delta \H^l_{\nu}$ expresses the difference between restricted fine level Hessian and the initial coarse level Hessian
\begin{align}
\delta \H^l_{\nu} := 
\begin{cases}
\Rev^l_{l+1} \H^{l+1}_{\nu,\mu_1} \I^l_{l+1}  - \nabla^2 f_{\nu}^{l}(\x^{l}_{\nu,0}) &  \text{if}  \ \ \ l < L\\
0  & \text{if}  \ \ \ l = L.
\end{cases}
\label{eq:delta_H}
\end{align}
The presence of this term enforces the first coarse level correction to be the restricted Newton's direction. 
Subsequent coarse level corrections involve knowledge from the coarse level discretization of the energy functional. 
The multilevel method with $\{ \h^l \}_{l=0}^{L}$ defined as in \eqref{eq:coarse_objective_second_order} is demonstrated to guarantee at least as good and usually better convergence rate than multilevel method with first order consistency or with the Galerkin approximation model, for the details see \cite{Yavneh2006}.  
Despite its efficiency, the method is computationally more costly than previously discussed methods. 
This is due to the fact, that the model \eqref{eq:coarse_objective_second_order} requires knowledge about the coarse level discretization as well as the information about the matrix $\delta \H^l_{\nu}$. 
In addition, the product consisting of $\delta \H^l_{\nu}$ is used to update the coarse level gradient and the Hessian at each successful coarse level iteration.

\subsection{Solution dependent second order consistency}
\label{section:modified_second_order}
In this section, we introduce a novel model in order to define of coarse level objective function $\h^l: \R^{n^l} \rightarrow \R$.
The proposed model is basically variation of the \eqref{eq:coarse_objective_second_order} defined as  following
\begin{align}
\h^{l}_{\nu}(\x_{\nu}^{l}) := \underbrace{\tilde{f}_{\nu}^{l}(\x_{\nu}^{l})}_{\text{modified coarse model}} + 
\underbrace{ \langle \delta \g^l_{\nu}, \x^l_{\nu} - \x^{l}_{\nu,0} \rangle }_{\text{1st order coupling}} + 
\chi_1(\mathbf{c}, \epsilon) \underbrace{ \frac{1}{2} \langle  \x^l_{\nu} - \x^{l}_{\nu,0}, \delta \H  (\x^l_{\nu} - \x^{l}_{\nu,0} )\rangle.}_{\text{2nd order coupling}} 
\label{eq:coarse_objective_modified2nd_order}
\end{align}
The two main ingredients of the formulation \eqref{eq:coarse_objective_modified2nd_order} are
\begin{itemize}
\itemsep-0.25em
\item the indicator function $\chi_1$ is used to turn on and off the second order coupling term. 
\item the objective function ${f}_{\nu}^{l}(\x_{\nu}^{l})$ is replaced by its modified variant $\tilde{f}_{\nu}^{l}(\x_{\nu}^{l})$.
\end{itemize}

The first modification is justified by the fact that during the pre-stage phase of the simulation, the crack has not yet propagated. 
Thus, the coarse level discretization is a valid representation of the fine level model, with respect to the discretization error. 
In this case, the use of the first order consistency approach is sufficient and computationally cheaper.
To this aim, we define indicator function
\begin{align}
\chi_1(\mathbf{c}, \epsilon) := 
\begin{cases}
0 &  \text{if}  \ \ \  \underset{i = 0, \cdots,  n^L}{\max(c_i)} > \epsilon   \\
1  &  \text{otherwise},
\end{cases}
\label{eq:indicator1}
\end{align}
which allows for turning on and off the second order consistency term. 
The symbol $c_i$ in \eqref{eq:indicator1} denotes the nodal coefficients of the fine level phase-field.
A threshold $\epsilon \in \R$,  $0 \leq \epsilon < 1$ defines a limit, after which we believe that the fine level model can not be represented well by the coarse level discretization. 
In our experience, $\epsilon \approx 0.85$ provides reasonable choice.


Our second modification consists of replacing the objective function ${f}_{\nu}^{l}(\x_{\nu}^{l})$ by its modified variant $\tilde{f}_{\nu}^{l}(\x_{\nu}^{l})$.
As the staggered solution scheme \cite{miehe2010phase}, we split the total potential energy \eqref{eq:potential_energy_pf} into the elastic and the fracture energy. 
For the elastic part of the energy, we employ a coarse level discretization, since it can represent the fine level elastic energy well.
On the other hand, for the fracture energy our coarse model uses only information restricted from the fine level, since the regularized fracture energy differs at the different discretization levels.
To incorporate both of the above requirements, the modified objective function $\tilde{f}_{\nu}^{l}$ is defined as the modified energy functional 
\begin{align}
\tilde{\psi}(\u, c) := \int_{\Omega} \Big(  [g(c)+k] \psi_e^+(\eps(\u)) + \psi_e^-(\eps(\u))  \Big) \  d \Omega +
\chi_1  \int_{\Omega} \Big(  \pazocal{G}_c \Big[   \frac{1}{2 l_s }c^2 + \frac{l_s}{2} | \nabla c |^2 \Big]  \Big) \ d \Omega, 
\label{eq:modified_energy_ML}
\end{align}
where the indicator function  $\chi_1$ is defined as in \eqref{eq:indicator1}. 

As it will be demonstrated in Section \ref{section:results_different_models}, the coarse level model \eqref{eq:coarse_objective_modified2nd_order}, provides better coarse level error approximation when compared to the coarse level models \eqref{eq:coarse_objective_galerkin}, \eqref{eq:coarse_objective_first_order}, \eqref{eq:coarse_objective_second_order}. 
Conceptually similar approach was applied in the context of multilevel method for the frictional contact in \cite{RKornhuber_RKrause_2001}, see also \cite{krause2009nonsmooth}.

\section{Transfer operators}
\label{section:transfer}
The pseudo-$L^2$-projection belongs to the more generic group of mortar/variational transfer operators, which are commonly used in the context of non-conforming domain decomposition methods and contact problems \cite{wohlmuth2000mortar}.   
The mortar projections allow for transferring the information of functional quantities from one approximation space to the another. 
Following the naming convention in the literature on the variational transfer, we transfer the data from the master space to the slave space.
Therefore, in the following discussion, we use subscripts $m$ and $s$  to denote quantities related to the master and slave side, respectively. 

\subsection{Variational transfer}
Let us consider meshes $\pazocal{T}^m$ and $\pazocal{T}^s$  with corresponding finite element spaces $\pazocal{X}^m$ and $\pazocal{X}^s$.
The mortar projection $\Pi: \pazocal{X}^m \rightarrow \pazocal{X}^s$ maps a function from  $\pazocal{X}^m$ to $\pazocal{X}^s$. 
For a function $m_h \in \pazocal{X}^m$, we are looking for $s_h = \Pi(m_h) \in \pazocal{X}^s$, such that the following equality condition holds
\begin{align}
(\Pi(m_h), q_h)_{L_2(\pazocal{I})} = (s_h, q_h)_{L_2(\pazocal{I})} \ \ \ \ \ \  \forall q_h \in \pazocal{X}^q,
\label{eq:eq_cond}
\end{align}
where $\pazocal{I} = \pazocal{T}^m \cap \pazocal{T}^s$.
The space $\pazocal{X}^q$ in \eqref{eq:eq_cond} denotes the multiplier space.
Reformulating \eqref{eq:eq_cond} results in the following weak equality condition
\begin{align}
\int_{\pazocal{I}} (m_h - \Pi(m_h)) q_h \ d \x  = \int_{\pazocal{I}}  (m_h - s_h) q_h \ d \x  = 0  \ \ \ \ \ \  \forall q_h \in \pazocal{X}^q.
\label{eq:weak_eq}
\end{align}
Expressing $m^h \in \pazocal{X}^m$, $s^h \in \pazocal{X}^s$ and $q^h \in \pazocal{X}^q$ in terms of their respective basis, $\{ N^m_i \}_{i=1}^{n^m}$, $\{ N^s_j \}_{j=1}^{n^s}$, $\{ N^q_k \}_{k=1}^{n^q}$, 
we obtain 
$m^h = \sum_{i=1}^{n^m} N^m_i m_i$,  $s^h = \sum_{j=1}^{n^s} N^s_j s_j$,  $q^h = \sum_{k=1}^{n^q} N^q_k q_k$,
where $\{ m_i \}_{i=1}^{n^m}$, $\{ s_j \}_{j=1}^{n^s}$, $\{ m_k \}_{k=1}^{n^q}$ are the coefficients.
This allows us to reformulate equation \eqref{eq:weak_eq} as
\begin{align}
\sum_{i=1}^{n^m} m_i \int_{\pazocal{I}}  N^m_i   N^q_k \ d \x          = \sum_{j=1}^{n^s}  s_j  \int_{\pazocal{I}}  N^s_j   N^q_k \ d \x   \ \ \ \ \ \text{for} \ \ k \in \{ 1, \dots, n^q \},
\label{eq:transfer}
\end{align}
which can be expressed by using vector-matrix notation as 
\begin{align}
\B \m = \D \mathbf{s}, 
\label{eq:transfer_mat}
\end{align}
where  the vectors $\m$ and $\mathbf{s}$ contain the coefficients $m_i$ and $s_j$.
The components of the coupling matrix $\B$ are defined as $B_{k,i}= \int_{\pazocal{I}}  N^m_i   N^q_k \ d \x$, while the matrix $\D$ contains the entries $D_{k,j} = \int_{\pazocal{I}}  N^s_j   N^q_k \ d \x$. 
As a result, the  discrete projection operator has the following algebraic representation
\begin{align}
\PI = \D^{-1} \B.
\label{eq:transfer_discret}
\end{align}

\subsection{Transfer operators via pseudo-$L^2$-projection}
Different choices of the multiplier space $\pazocal{X}^q$ result in different transfer operators.
For example, if we set $\pazocal{X}^q = \pazocal{X}^s$, the operator $\PI$ is nothing other than the well-known $L^2$ projection. 
In this particular case, the matrix $\D$ is the mass matrix on the slave space. 
For the purpose of this work, we consider the computationally cheaper alternative, e.g. pseudo-$L^2$ projection, where the multiplier space $\pazocal{X}^q$ is associated with the dual basis \cite{dickopf2010multilevel, oswald2001polynomial}. 
The dual basis functions satisfy the following bi-orthogonality condition
\begin{align}
(N_j^s, N_k^q)_{L_2(\pazocal{I})} = \delta_{jk} \int_{\pazocal{I}} N_j^s \ d \x, 
\end{align}
which gives rise to a diagonal matrix $\D$ in our setting.
This is very convenient, as it makes an evaluation of the inverse in \eqref{eq:transfer_discret} computationally cheap. 
In addition, it produces a transfer operator with a sparse structure, which is not the case for the $L^2$-projection.
The parallel implementation of the mortar projections is not a trivial task, however, a library MOONoLith \cite{krause2016parallel} provides an efficient and scalable implementation.

As described in the Section \ref{section:multilevel}, the RMTR Algorithm \ref{alg:rmtr} makes use of three different transfer operators: prolongation, restriction, and projection. 
The prolongation operator $\I_{l}^{l+1}: {\pazocal{X}}^{{l}} \rightarrow  {\pazocal{X}}^{{l+1}}$  transfers  primal variables from a coarse level to the subsequent fine level. 
The assembly of the prolongation operator $\I_{l}^{l+1}$ is performed by the formula \eqref{eq:transfer_discret}, with $\pazocal{X}^m = \pazocal{X}^{{l}}$ and $\pazocal{X}^s =  {\pazocal{X}}^{{l+1}}$.
On the other hand, in order to assemble the projection operator $\P^l_{l+1} : \pazocal{X}^{{l+1}} \rightarrow \pazocal{X}^{{l}}$, we associate master and slave spaces as follows:
$\pazocal{X}^m = \pazocal{X}^{{l+1}}$ and $\pazocal{X}^s =  {\pazocal{X}}^{{l}}$.
The resulting operator $\P^l_{l+1} $ is then used to transfer iterates from the fine level to the coarse level.
The restriction operator $\Rev^l_{l+1} : \pazocal{X}^{{l+1}} \rightarrow \pazocal{X}^{{l}}$, which is designed for  transferring the dual variables from the fine level to the coarse is not assembled by using \eqref{eq:transfer_discret}. Instead, it is obtained as the adjoint of the prolongation operator, thus $\Rev^l_{l+1} = (\I_{l}^{l+1})^T$. 

It is important to realize, that $\P^l_{l+1} \neq (\I_{l}^{l+1})^T$. 
As noted by Gro{\ss} and Krause in  \cite{Gross2009}, the restriction operator $\Rev^l_{l+1} = (\I_{l}^{l+1})^T$ is not suitable for transferring the primal variables from the fine level to the coarse level. 
This is due to the fact, that the restriction operator $\Rev^l_{l+1}$  is designed as a dual operator. 
The numerical results reported in \cite{Gross2009} also demonstrate that the usage of the projection operator $\P^l_{l+1} $ based on $L^2$-projection accelerates the convergence behavior of the RMTR method.  The detailed discussion about how to construct suitable operator $\P^l_{l+1} $  can be also found in \cite{Kopanicakova_2019a}.

\begin{remark}
We note, that we refer to the operator $\P^l_{l+1}$ as to the projection.  
However, both transfer operators $\I_{l}^{l+1}$ and $\P^l_{l+1}$ are technically projections. 
\end{remark}
\begin{remark}
As described in Section  \ref{section:multilevel}, the mesh $\pazocal{T}^{l+1}$  is obtained by uniformly refining the mesh $\pazocal{T}^{l}$. 
In this case, i.e.  $\pazocal{T}^{l} \subset \pazocal{T}^{l+1}$, the prolongation operator $\I_{l}^{l+1}: {\pazocal{X}}^{{l}} \rightarrow  {\pazocal{X}}^{{l+1}}$ obtained as described above,  coincides with standard finite element interpolation. 
\end{remark}
\begin{remark}
In order to obtain a computationally efficient simulation, we perform the assembly of the transfer operators just once, e.g. at the beginning of the simulation. 
This is possible, as meshes $\{ \pazocal{T}^l \}_{l=0}^{l=L}$ are not altered during the simulation.
\end{remark}

\section{Multilevel treatment of point-wise constraints}
\label{section:constrains}
As discussed in Section \ref{section:multilevel}, the minimization on level $l$ is performed with respect to some feasible set $\pazocal{L}^l_{c}$. 
The feasible set $\pazocal{L}^l_{c}$ is constructed as an intersection of two sets
\begin{align}
\pazocal{L}^l_{\nu} := \pazocal{F}^l_{\nu} \cap \pazocal{S}^l_{\nu}, 
\label{eq:define_l}
\end{align}
where $\pazocal{F}^l_{\nu}$ defines the restricted lower bound and $\pazocal{S}^l_{\nu}$ carries over the restricted trust region bounds. 
We note, that the feasible set $\pazocal{L}^l_{c}$ is V-cycle dependent.
Therefore, it has to be re-evaluated every time the algorithm enters the given level. 

In order to define the set $\pazocal{F}^l_{\nu}$, we follow an approach of Gelman et. al \cite{Gelman1990}. 
The alternatives can be found for example in \cite{krause2009nonsmooth, RKornhuber_RKrause_2001, Kornhuber1994}. 
The level dependent feasible set $\pazocal{F}^l_{\nu}$ is then defined as follows
\begin{align}
\pazocal{F}^l_{\nu} := \{ \x^l \ | \ \lb^l_{\nu} \leq \x^l  \},  
\end{align}
where the lower bound $\lb^l_{\nu}$ is obtained component-wise as 
\begin{align}
(\lb^l_{\nu})_k := (\x_{\nu,0}^l)_k + \underset{j = 1, \cdots, n^l}{\max} [ (\lb^{l+1}_{\nu} - \x^{l+1}_{\nu, \mu_1} )_j ]. 
\label{eq:hard_constraints}
\end{align}
Formula \eqref{eq:hard_constraints} ensures that the prolongated coarse level corrections do not violate the fine level irreversibility constraint, 
thus the following relation  holds
\begin{align}
\x^{l+1}_{\nu, \mu_1} + \I_{l}^{l+1} ( \x^l_{\nu}  - \P_{l+1}^{l} \x^{l+1}_{\nu, \mu_1}  )   \in \pazocal{F}^{l+1}    \ \ \ \ \ \text{for all}  \ \ \x^l_{\nu} \in \pazocal{F}^l. 
\end{align}

\begin{remark}
The restriction of the point-wise constraints to the coarse level by using formula \eqref{eq:hard_constraints} is known to provide narrow approximation on the feasible set \cite{Kornhuber1994}. 
However, test runs without the coarse level constraints \eqref{eq:hard_constraints} showed similar convergence behavior, in terms of number of iterations.
This can be taken as an indicator that, on the one hand, the coarse level constraints are not too restrictive and that, on the other hand, the constructed model hierarchy provides a good multi-scale decomposition of the fine level model.
\end{remark}

The set  $\pazocal{S}^l_{\nu}$ containing trust-region constraints inherited from the $(\nu, \mu_1, l+1)$ iterate
 is also represented component-wise as
\begin{align}
\pazocal{S}^l_{\nu} := \{ \x^l_{\nu} \ | \ \tl^l_{\nu} \leq \x^l_{\nu}   \leq   \tu^l_{\nu}  \},
\end{align}
where 
\begin{align}
(\tl^l_{\nu})_k := \sum_{j=1}^{n^{l+1}} (P_{l+1}^{l})_{kj} ( \max [ \tl^{l+1}_{\nu}, \x^{l+1}_{\nu, \mu_1} - \Delta^{l+1}_{\nu, \mu_1} \e ] )_j,
\label{eq:lower_tr_a}
\end{align}
and
\begin{align}
(\tu^l_{\nu})_k := \sum_{j=1}^{n^{l+1}} (P_{l+1}^{l})_{kj} ( \min [  \tu^{l+1}_{\nu}, \x^{l+1}_{\nu, \mu_1}   +  \Delta^{l+1}_{\nu, \mu_1} \e ] )_j. 
\label{eq:upper_tr_b}
\end{align}
The $\max$ and $\min$ operations in \eqref{eq:lower_tr_a} - \eqref{eq:upper_tr_b} are taken component-wise and $\e \in \R^{n^{l+1}}$ denotes the vector of ones. 
The restrictions \eqref{eq:lower_tr_a} - \eqref{eq:upper_tr_b} of the fine level trust region bound are weaker than the restriction of the irreversibility constraint since it is
not necessary that they are satisfied by all iterates. 
This is due to the fact that a bounded violation of the trust region constraints is acceptable by the RMTR Algorithm \ref{alg:rmtr}, see \cite{gratton2008_inf} for the rigorous details.

\begin{remark}
We assume, that components of the lower bound vectors $\lb^{L+1}_{\nu}$, $\tl^{L+1}_{\nu}$ equal to $- \infty$, while components of the upper bound vector $\tu^{L+1}_{\nu}$ are $\infty$. 
\end{remark}
 
\section{Numerical experiments}
\label{section:numerical_examples}
In this section, we present several numerical examples, which we use in order to examine the performance of RMTR method.
All of the presented examples are three dimensional and employ the length scale parameter $l_s = 2 h$, where $h$ denotes a size of the biggest element. 
The set of material parameters used during our simulations can be found in Table \ref{table:params}.
\begin{table}
\centering
\caption{Setup of parameters used in numerical experiments.}
\begin{tabular}{ ccccc} 
\hline
\multirow{2}{*}{ \parbox{1.7cm}{\textcolor{myred}{Material parameter}}}        &  \multicolumn{4}{c}{\textcolor{myred}{Test case}}               			             						 			\\ 
       						& \textcolor{myred}{Fracture modes}    		& \textcolor{myred}{Mixed-mode}    				& \textcolor{myred}{Conchoidal fracture}  	& \textcolor{myred}{Pressurized fracture	}\\ \hline
$\lambda(\text{GPa}) $ 	& 	   $  12.1 $   				&	8.8					& 	100				&		$12$			\\ 
$\mu(\text{GPa}) $ 		&  	   $ 	7.9  $				&	1.3					& 	100				&		$8$			\\
$\pazocal{G}_c(\text{GPa})$  	& 	   $5 \times 10^{-4}$			&	$ 1.1 \times 10^{-4}$	& 	$1 \times 10^{-3}$	&		$1 \times 10^{-3}$			\\ \hline
\end{tabular}
\label{table:params}
\end{table}

\subsection{Fracture modes}
\label{section:fracture_modes}
\begin{figure}
\hspace{1cm}
\subfloat[\label{fig:3modes_a}]{
\begin{tikzpicture}
\node[anchor=south west,inner sep=0] at (0,0) {\includegraphics[scale=0.1]{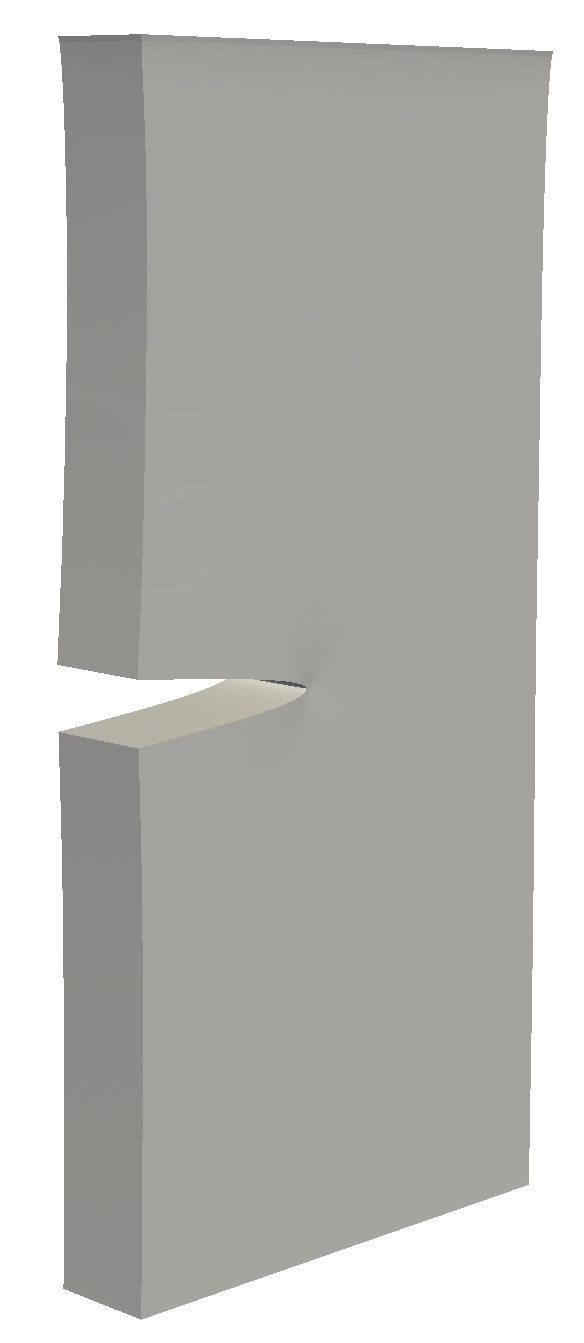}};
\node[color=black] at (1.1, 5.4) {$\bm{u}$};
\draw[style=thick, ->, yshift=115pt, xshift=20pt, xslant=0.0,color=black] (0, 0.6) -- (0, 1.1);
\draw[style=thick, ->, yshift=115pt, xshift=20pt, xslant=0.0, color=black] (0.4, 0.6) -- (0.4, 1.1);
\draw[style=thick, ->, yshift=115pt, xshift=20pt, xslant=0.0, color=black] (0.8, 0.6) -- (0.8, 1.1);
\end{tikzpicture}
}
\hspace{2.5cm}
\subfloat[\label{fig:3modes_b}]{
\begin{tikzpicture}
\node[anchor=south west,inner sep=0] at (0,0) {\includegraphics[scale=0.1]{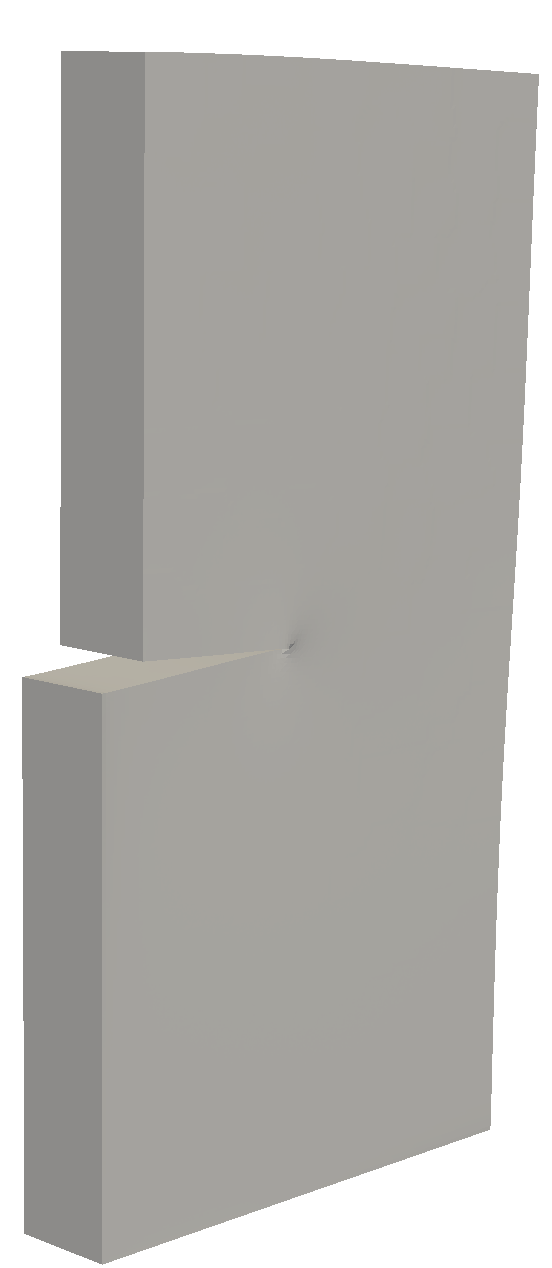}};
\node[color=black] at (1.15, 4.85) {$\bm{u}$};
\draw[style=thick, ->, yshift=115pt, xshift=20pt, xslant=0.0, color=black] (-0.3, 0.48) -- (0.1, 0.44);
\draw[style=thick, ->, yshift=115pt, xshift=20pt, xslant=0.0, color=black] (0.2, 0.45) -- (0.6, 0.40);
\draw[style=thick, ->, yshift=115pt, xshift=20pt, xslant=0.0, color=black] (0.7, 0.4) -- (1.1, 0.35);
\end{tikzpicture}
}
\hspace{2.5cm}
\subfloat[\label{fig:3modes_c}]{
\begin{tikzpicture}
\node[anchor=south west,inner sep=0] at (0,0) {\includegraphics[scale=0.1]{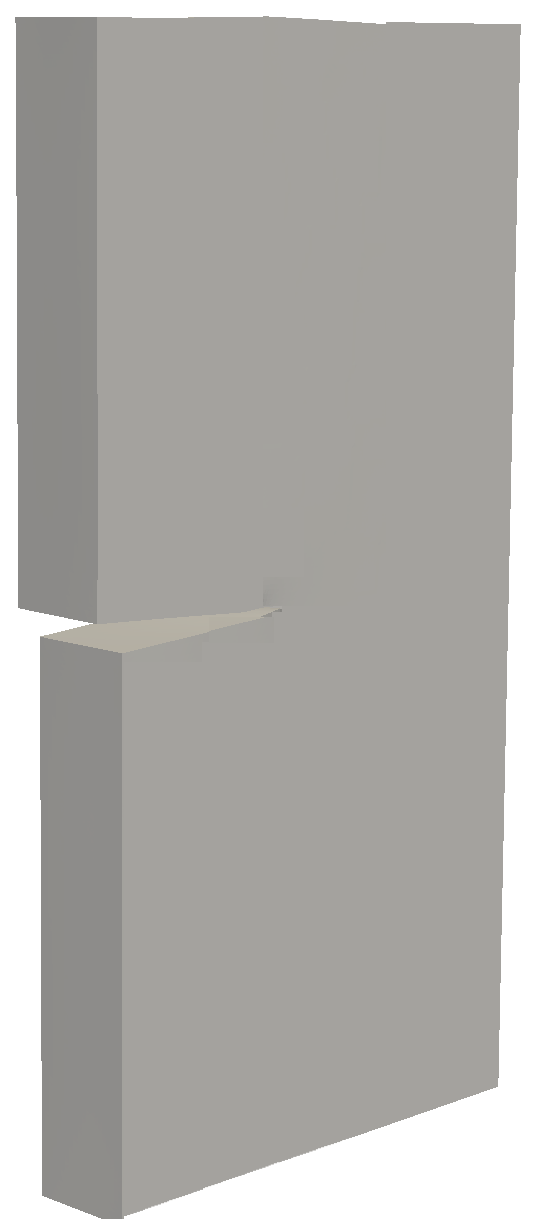}};
\node[yshift=115pt, xshift=20pt, color=black] at (0.4,  -2-0.47-0.47) {$\bm{u}$};
\node[yshift=115pt, xshift=20pt, color=black] at (-1.4, -2+0.47+0.47+0.47) {$\bm{u}$};
\draw[style=thick, ->, yshift=115pt, xshift=20pt, xslant=0.0, color=black] (-0.25, -2-0.47-0.47-0.47) -- (0.15, -2-0.47-0.47-0.47);
\draw[style=thick, ->, yshift=115pt, xshift=20pt, xslant=0.0, color=black] (-0.25, -2-0.47-0.47) -- (0.15, -2-0.47-0.47);
\draw[style=thick, ->, yshift=115pt, xshift=20pt, xslant=0.0, color=black] (-0.25, -2-0.47) -- (0.15, -2-0.47);
\draw[style=thick, <-, yshift=115pt, color=black] (-0.35, -2+0.47+0.47+0.47+0.47) -- (0.05, -2+0.47+0.47+0.47+0.47);
\draw[style=thick, <-, yshift=115pt, color=black] (-0.35, -2+0.47+0.47+0.47) -- (0.05, -2+0.47+0.47+0.47);
\draw[style=thick, <-, yshift=115pt, color=black] (-0.35, -2+0.47+0.47) -- (0.05, -2+0.47+0.47);
\end{tikzpicture}
}
\caption{ Boundary value problem set-up for three fracture modes. 
a) Mode I, tension test. We impose an incremental displacement $\Delta \u = 1 \times 10^{-4}$mm. 
b) Mode II., shear  test. We prescribe an incremental displacement  $\Delta \u = 1 \times 10^{-4}$mm.
c) Mode III.,   tearing test. An incremental displacement  $\Delta \u = 1 \times 10^{-4}$mm is enforced.}
\label{fig:3modes_init_setup}
\end{figure}
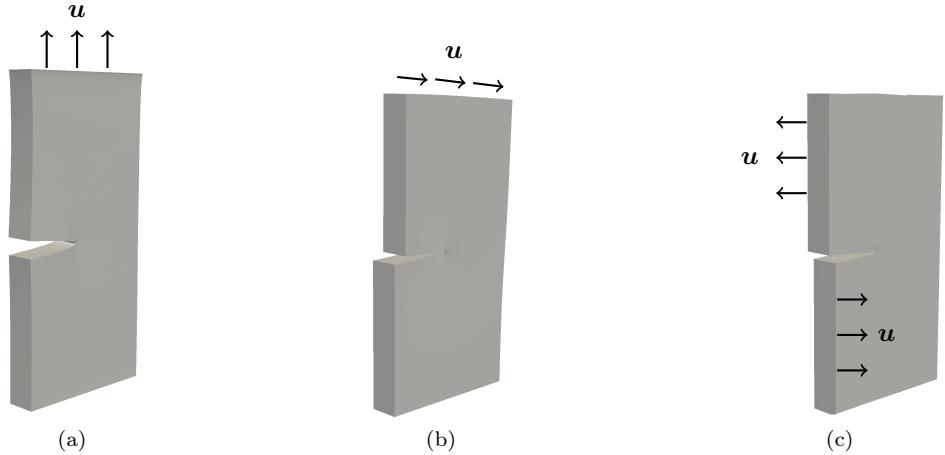

\begin{figure}
\subfloat[\label{fig:3modes_r1}]{
\includegraphics[scale=0.085]{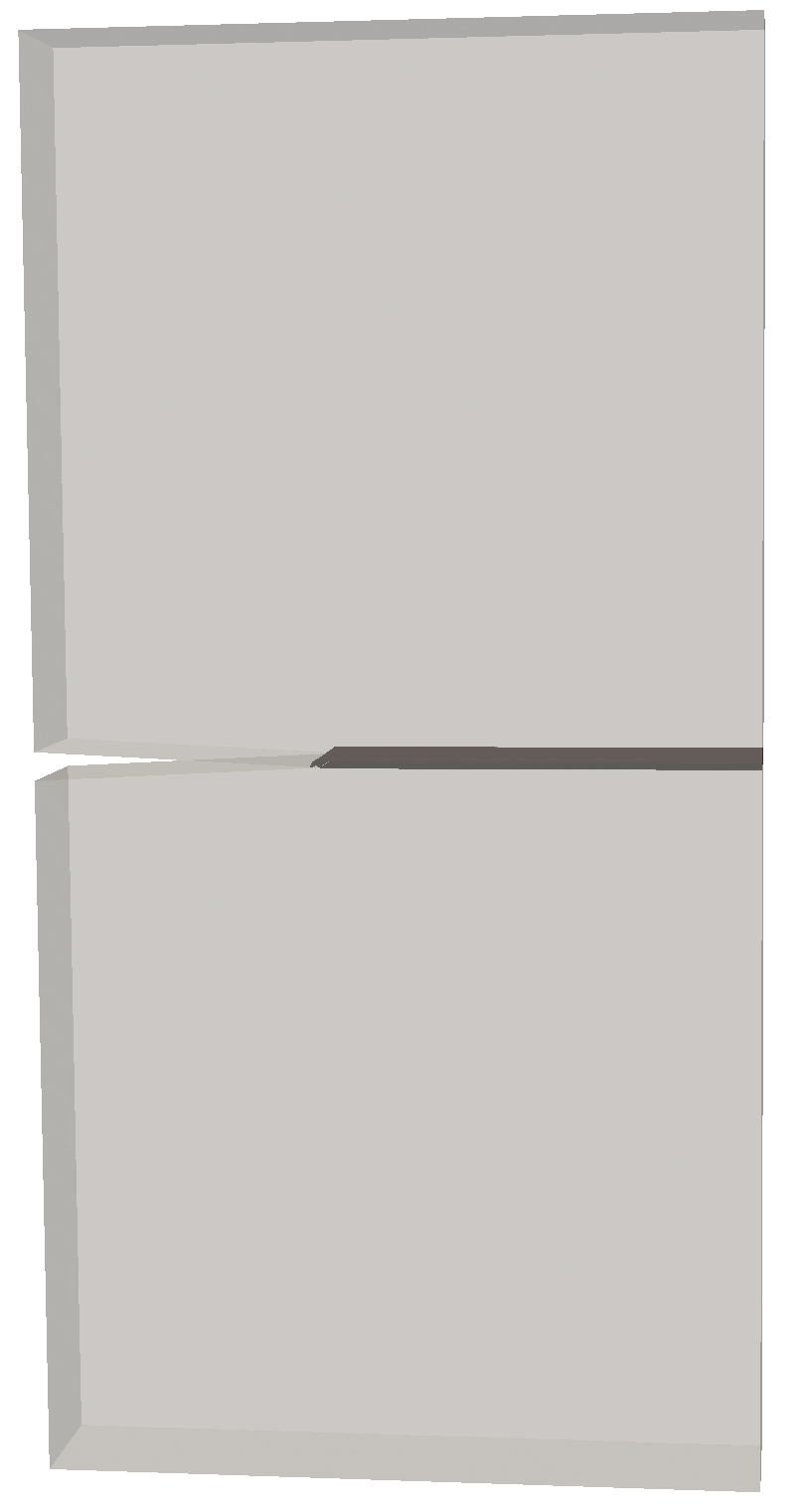}
}
\hfill
\subfloat[\label{fig:3modes_r2}]{
\includegraphics[scale=0.085]{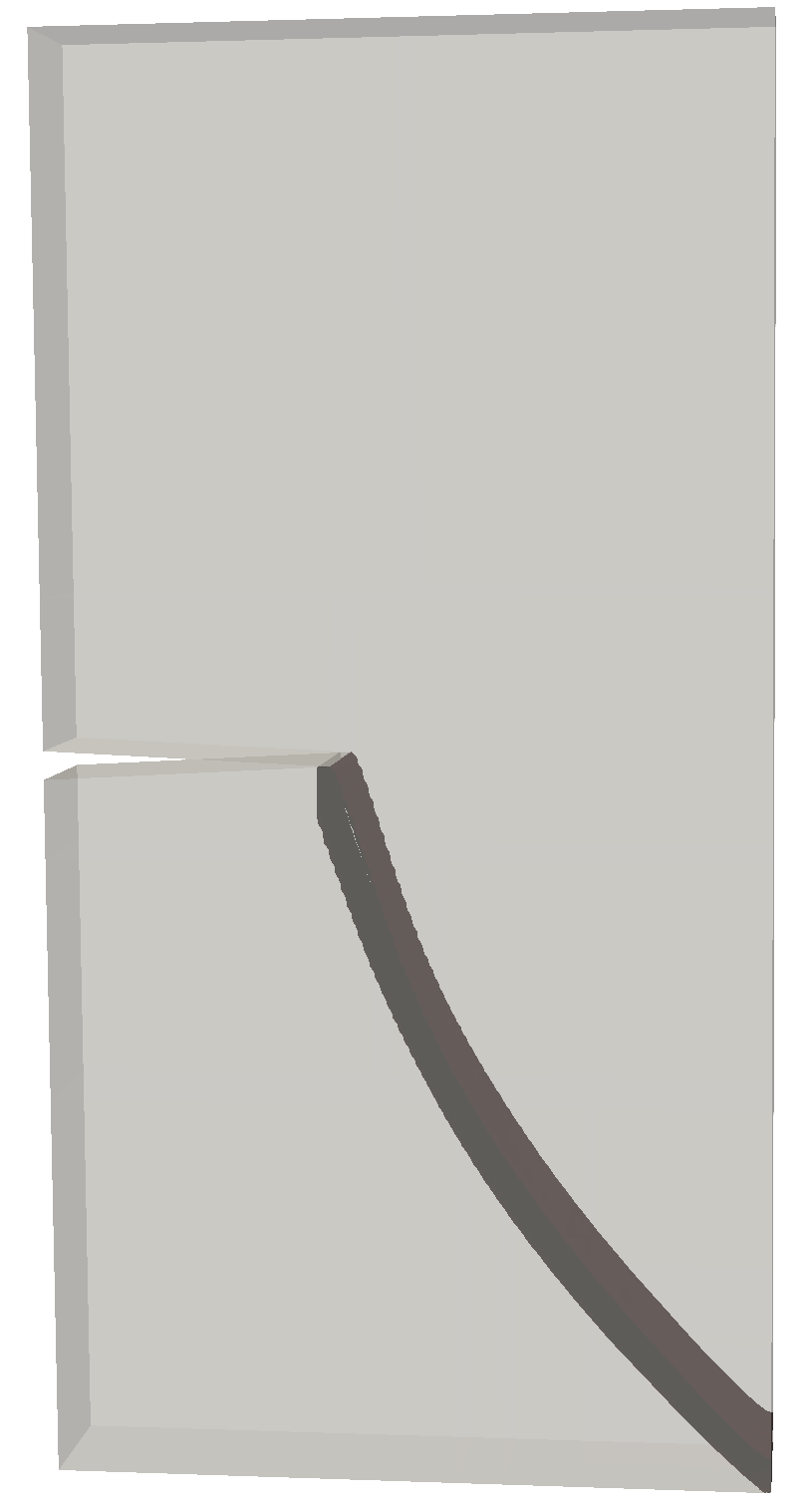}
}
\hfill
\subfloat[\label{fig:3modes_r3}]{
\includegraphics[scale=0.085]{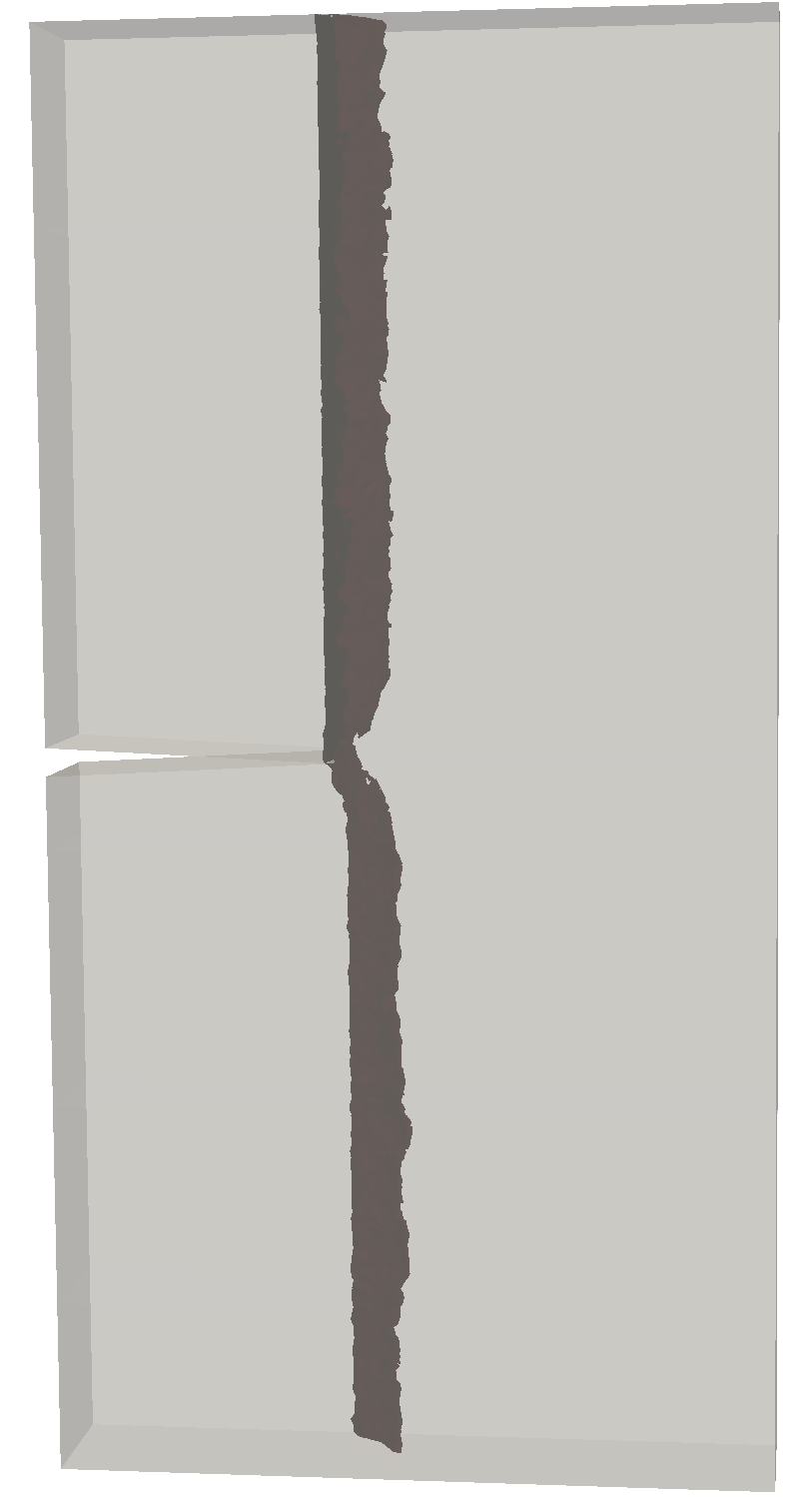}
}
\hfill
\subfloat{
\includegraphics[scale=0.05]{c_scale.png}
}
\caption{Iso-surfaces of the damage in the reference configuration for $c = 0.99$ for a) Tension b) Shear c) Tearing tests. 
The presented results originate from an numerical experiment with approx. $106$ thousands dofs.
}
\label{fig:3modes_result}
\end{figure}

At the beginning, we investigate the performance of the method on three fracture modes. 
We consider a three dimensional brick of size $10 \times 5 \times 2$ mm with $2$ mm initial notch. 
The specimen is then subjected to three different loading scenarios, which can force  crack propagation, namely 
\begin{enumerate}
\itemsep-0.25em
  \item Mode I.  \ \ - Tension (tensile stress).
  \item Mode II. \ - Shear (in-plane shear stress).
  \item Mode III. - Tearing (out-of-plane shear stress).
\end{enumerate}
Figure \ref{fig:3modes_init_setup} depicts an initial set-up and details on prescribed Dirichlet boundary conditions. 
The final crack patterns are shown on Figure \ref{fig:3modes_result} by depicting the iso-surface for  $c = 0.95$. 

\subsection{Mixed mode}
\begin{figure}
\subfloat[\label{fig:mixed_1}]{
	\begin{tikzpicture}[scale=1.8]

		\draw[color=gray, fill=gray] (0,-0.05) rectangle (2.05,0);
		\draw[color=gray, fill=gray] (2.0,0) rectangle (2.05, 0.975);		
		
		\draw[color=gray, fill=gray] (-0.05, 2) rectangle (2,2.05);
		\draw[color=gray, fill=gray] (-0.05, 2) rectangle (0.0, 1.025);		
	
		\draw[style=very thick, color=black] (0, 0) -- (2, 0);
		\draw[style=very thick, color=black] (2, 0) -- (2, 0.975);
		\draw[style=very thick, color=black] (2, 0.975) -- (1.75, 0.975);
		\draw[style=very thick, color=black] (1.75, 0.975) -- (1.75, 1.025);
		\draw[style=very thick, color=black] (1.75, 1.025) -- (2, 1.025);				
		\draw[style=very thick, color=black]  (2, 1.025)--(2,2) ;				
		\draw[style=very thick, color=black]  (2, 2)--(0,2) ;			
		\draw[style=very thick, color=black]  (0, 2)--(0,1.025) ;			
		\draw[style=very thick, color=black]  (0,1.025)  -- (0.25, 1.025);			
		\draw[style=very thick, color=black]  (0.25, 1.025) -- (0.25, 0.975);					
		\draw[style=very thick, color=black]  (0.25, 0.975) -- (0.0, 0.975);	
		\draw[style= very thick, color=black]  (0.0, 0.975) -- (0.0, 0.0);

		\draw[color=black, <->, latex-latex,  opacity=0.7] (0, 1.75) -- (2, 1.75);	
		\node[color=black] at (1, 1.65) {$200$};		
		
		\draw[color=black, <->, latex-latex, opacity=0.7] (0.6, 2) -- (0.6, 0);	
		\node[color=black, rotate=90] at (0.7, 1) {$200$};

		\draw[color=black, <->, latex-latex, opacity=0.7] (2, 0.925) -- (1.75, 0.925);
		\draw[color=black, -, opacity=0.7] (1.75, 0.975) -- (1.75, 0.875);				
		\node[color=black] at  (1.87, 0.81) {$25$};

		\draw[color=black, -, opacity=0.7]  (0.3, 1.025) -- (0.3, 0.975);
		\draw[color=black, -, opacity=0.7]  (0.25, 1.025) -- (0.35, 1.025);
		\draw[color=black, -, opacity=0.7]  (0.25, 0.975) -- (0.35, 0.975);		

		\draw[color=black, <-, latex-,  opacity=0.7]  (0.3, 1.025) -- (0.3, 1.175);
		\draw[color=black, <-, latex-,  opacity=0.7]  (0.3, 0.975) -- (0.3, 0.855);						
		\node[color=black] at  (0.41, 1) {$5$};

		\draw[color=black, ->,very  thick]  (1, 2.05) -- (1, 2.25);			
		\node[color=black] at (1.2, 2.15) {${F}_t$};

		\draw[color=black, ->,very  thick]  (1, -0.05) -- (1, -0.25);			
		\node[color=black] at (1.2, -0.15) {${F}_t$};

		\draw[color=black, <-,very  thick]  (2.05, 0.5) -- (2.25, 0.5);			
		\node[color=black] at (2.15, 0.3) {${F}_s$};

		\draw[color=black, <-,very  thick]  (-0.05, 1.5) -- (-0.25, 1.5);			
		\node[color=black] at (-0.15, 1.3) {${F}_s$};				
	\end{tikzpicture}
}
\hfill
\subfloat[\label{fig:mixed_2}]{
\includegraphics[scale=0.032]{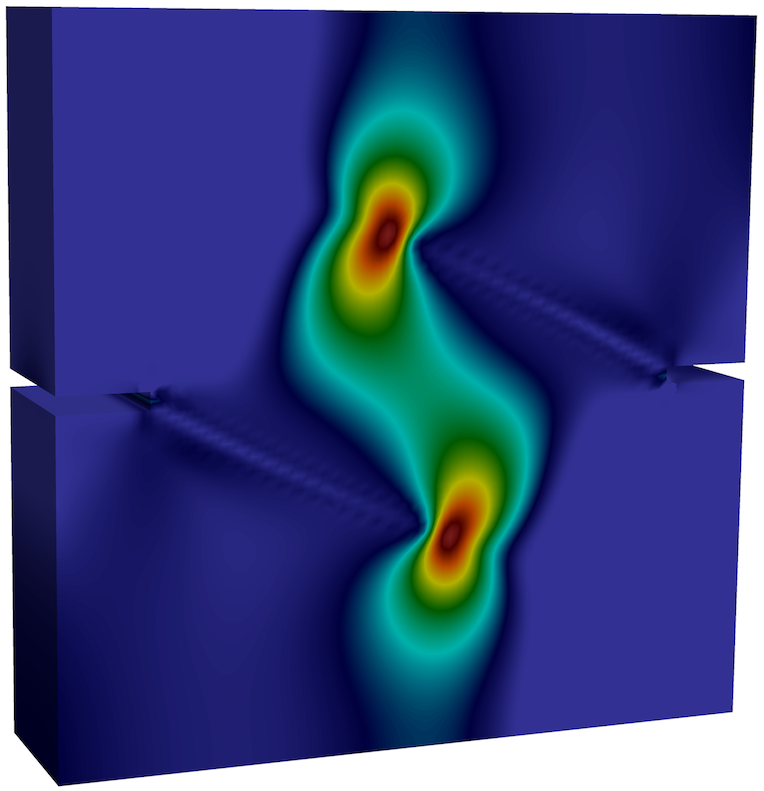}
}
\hfill
\subfloat{
\includegraphics[scale=0.045]{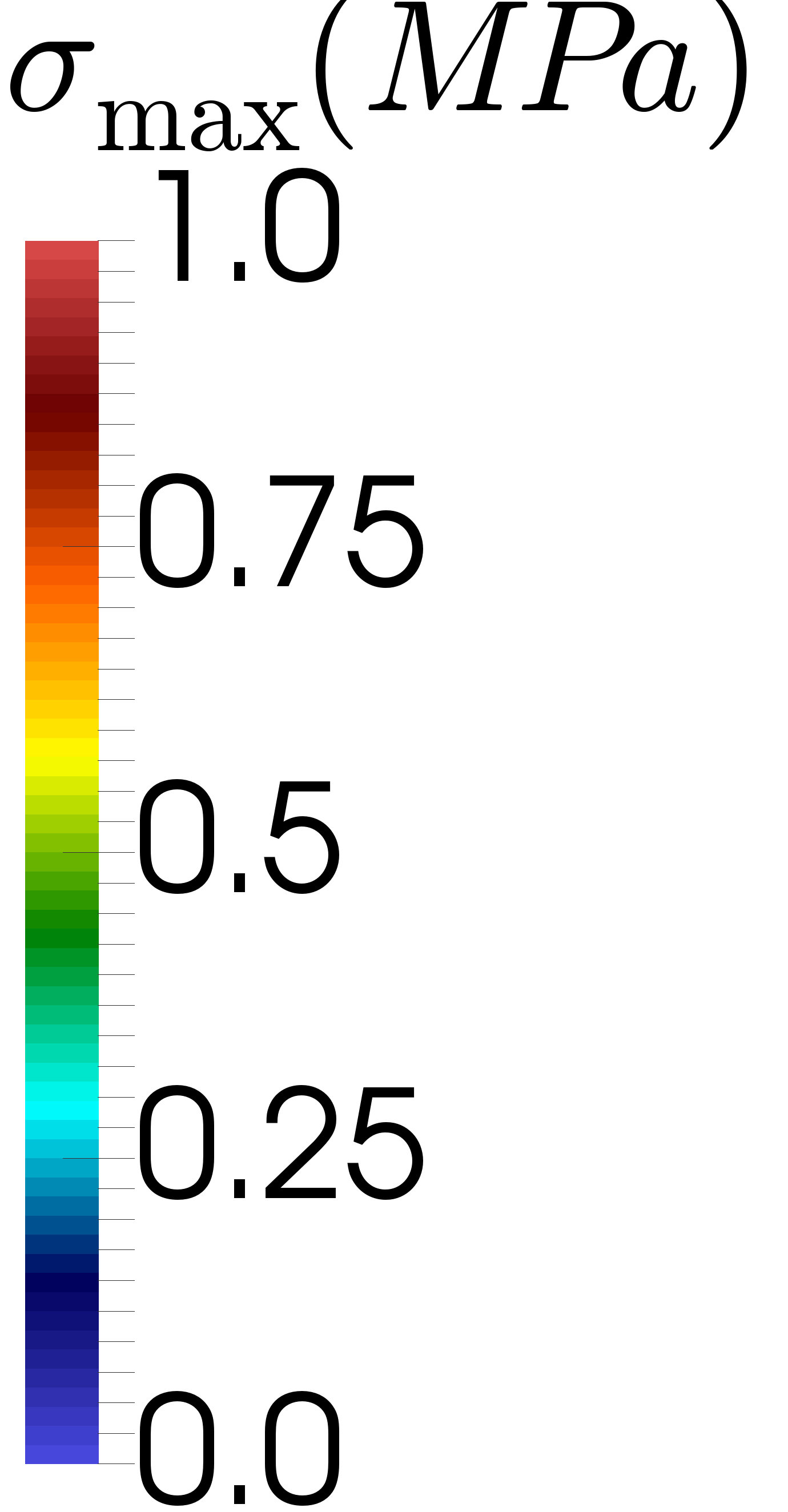}
}
\caption{Mixed-mode test.  a) Geometry and boundary value problem setup.  All dimensions are given in mm. The thickness of the cube is $50$mm. 
c) Principal stresses for at $t=7.3 \times 10^{-3}$. }
\label{fig:mixed_mode}
\end{figure}
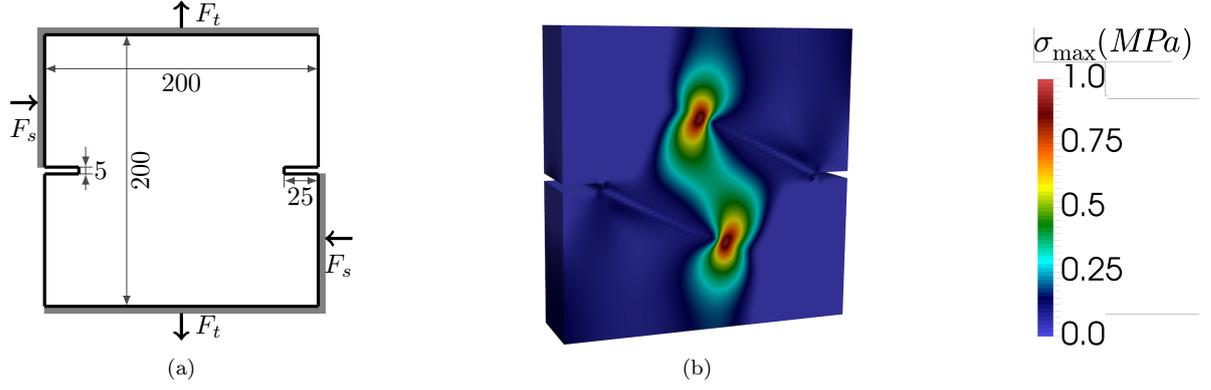

In our next test consider a mixed mode scenario, also know as the Nooru-Mohamed test.
The numerical experiment applies both tensile and shear loads to the double-notched specimen simultaneously. 
Here, we utilize scheme 6a) as proposed and empirically investigated in \cite{nooru1992mixed}. 
The dimensions of the specimen, as well as loading and boundary conditions, are illustrated in the Figure \ref{fig:mixed_1}. 
In particular, we apply shear $F_s$ and tensile $F_t$ forces, such that the ratio of axial and lateral deformation remains constant throughout the test, precisely $F_s/F_t =0.5$. 
The Nooru-Mohamed experiment was investigated in terms of phase-field fracture models in \cite{Zhang2017}. 

\begin{figure}
\subfloat[\label{fig:mixed_v1}]{
\includegraphics[scale=0.03]{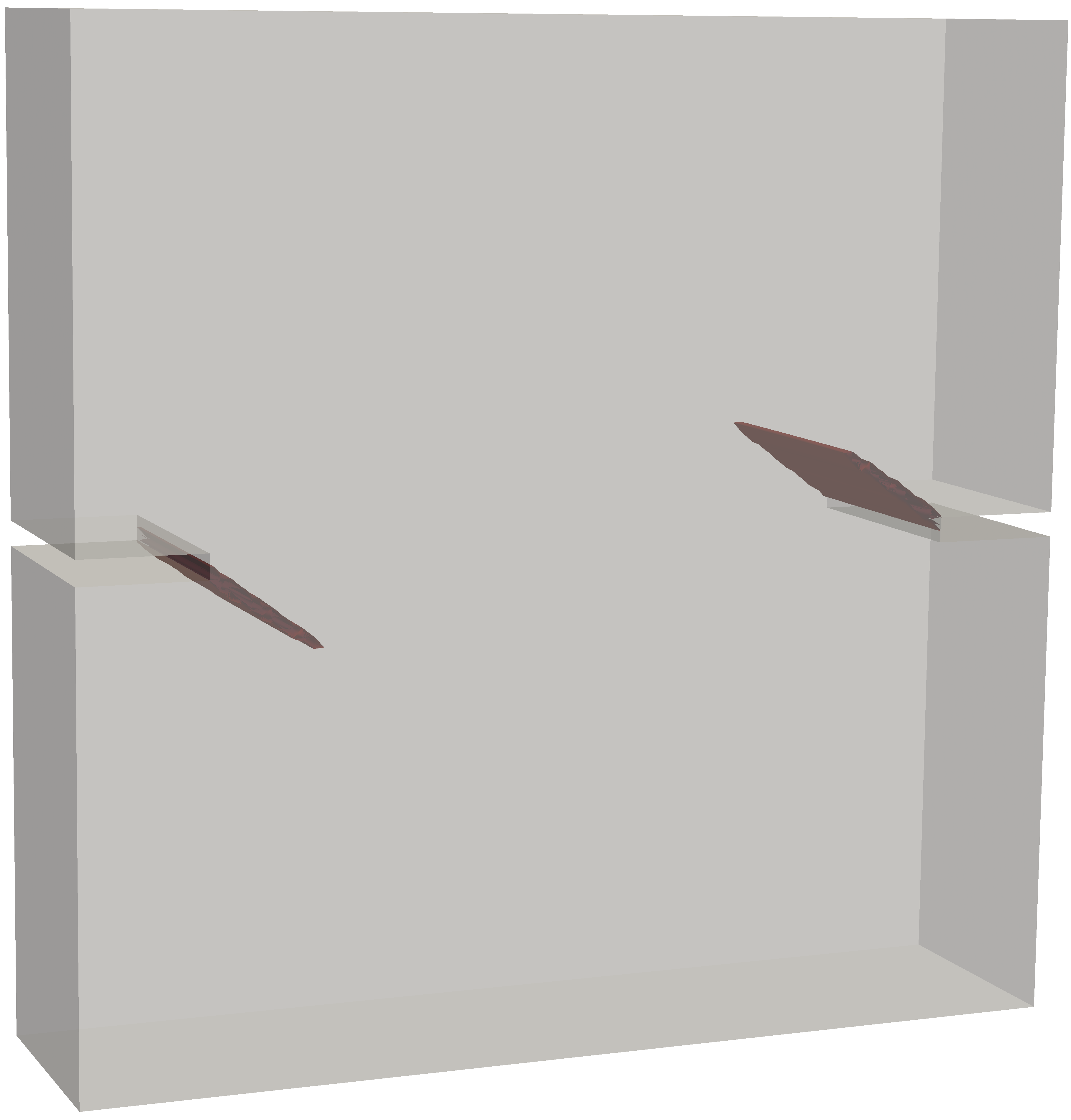}
}
\hfill
\subfloat[\label{fig:mixed_v2}]{
\includegraphics[scale=0.03]{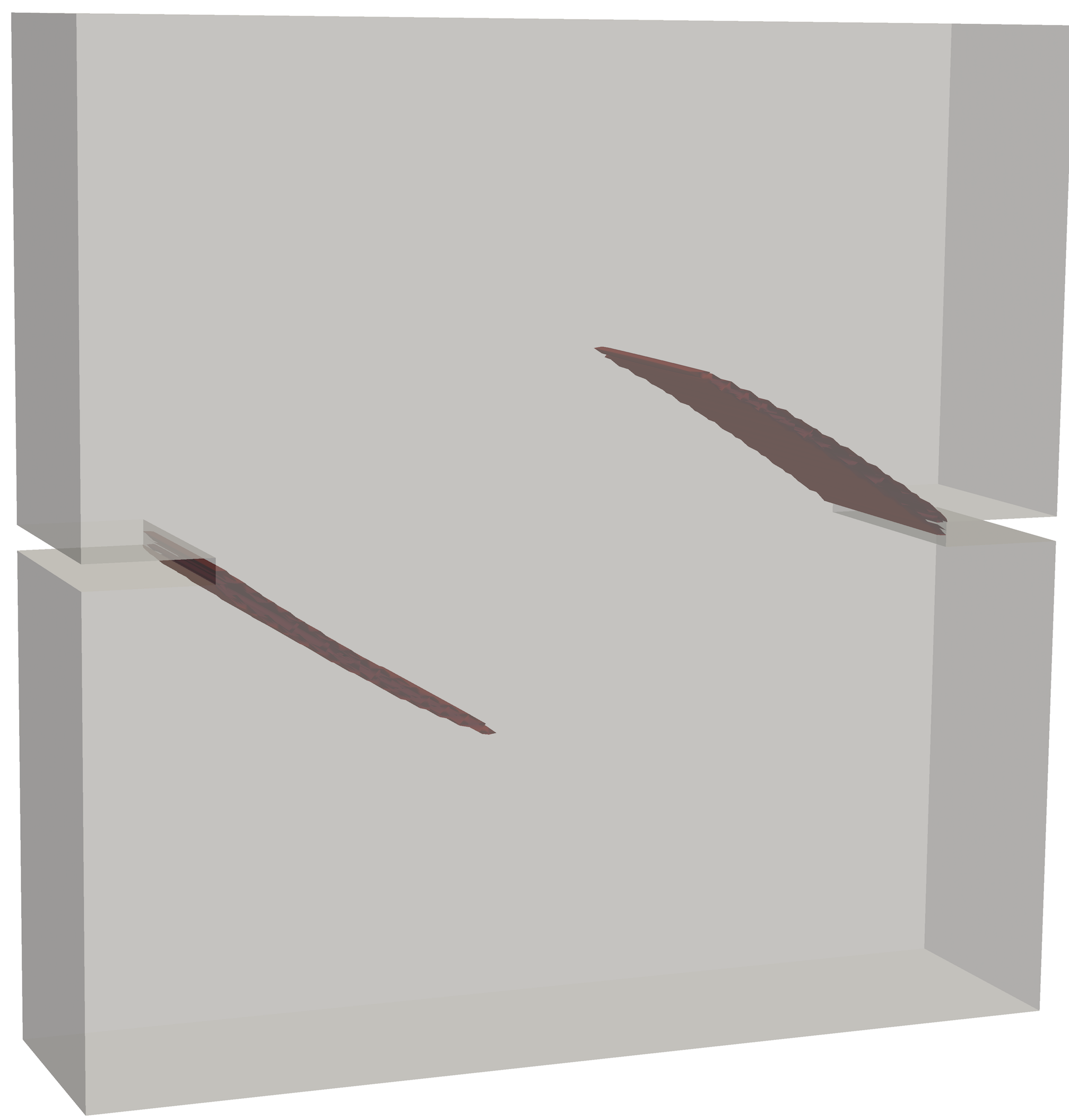}
}
\hfill
\subfloat[\label{fig:mixed_v3}]{
\includegraphics[scale=0.03]{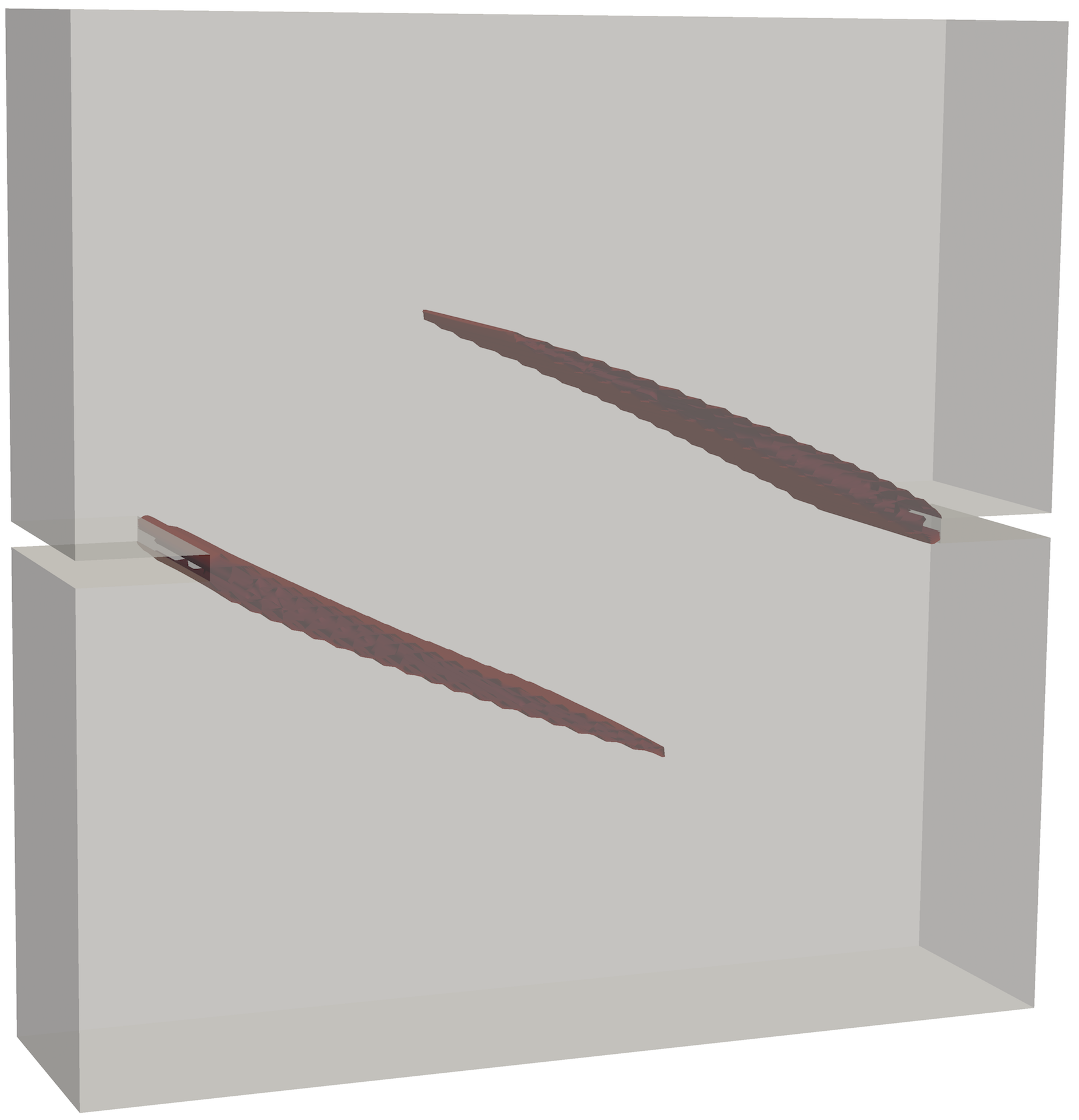}
}
\hfill
\subfloat{
\includegraphics[scale=0.04]{c_scale.png}
}
\caption{Mixed-mode test, fracture evolution in terms of crack iso-surfaces for $c=0.9$ at time:
 a) $t=3.5 \times 10^{-3}$   b) $t=4.5 \times 10^{-3}$  c)  $t=9 \times 10^{-3}$. The numerical experiment was performed with approx. $1,2$ million of dofs.}
\label{fig:mixed_mode_sim_result}
\end{figure}

Thanks to the specific loading protocol, the principal stresses rotate during the test, resulting in two  non-interacting,  non-planar,  curvilinear cracks, see Figure \ref{fig:mixed_mode_sim_result}.
The Figure \ref{fig:mixed_2}  illustrates the maximum principal Cauchy stress in the initial configuration. 

\subsection{Conchoidal fracture}
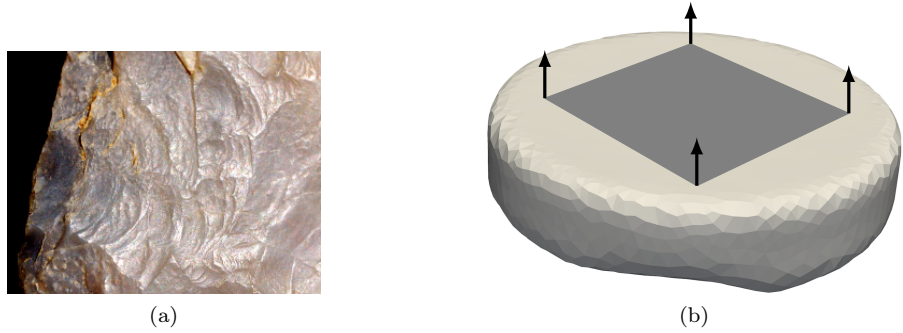
\begin{figure}[H]
\centering
\hfill
\subfloat[\label{fig:cf_1}]{
			\includegraphics[scale=0.4]{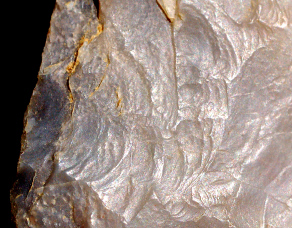}
}
\hfill
\subfloat[\label{fig:cf_2}]{
\begin{tikzpicture}[scale=0.8]
	\node[anchor=south west,inner sep=0] at (0,0) {\includegraphics[scale=0.044]{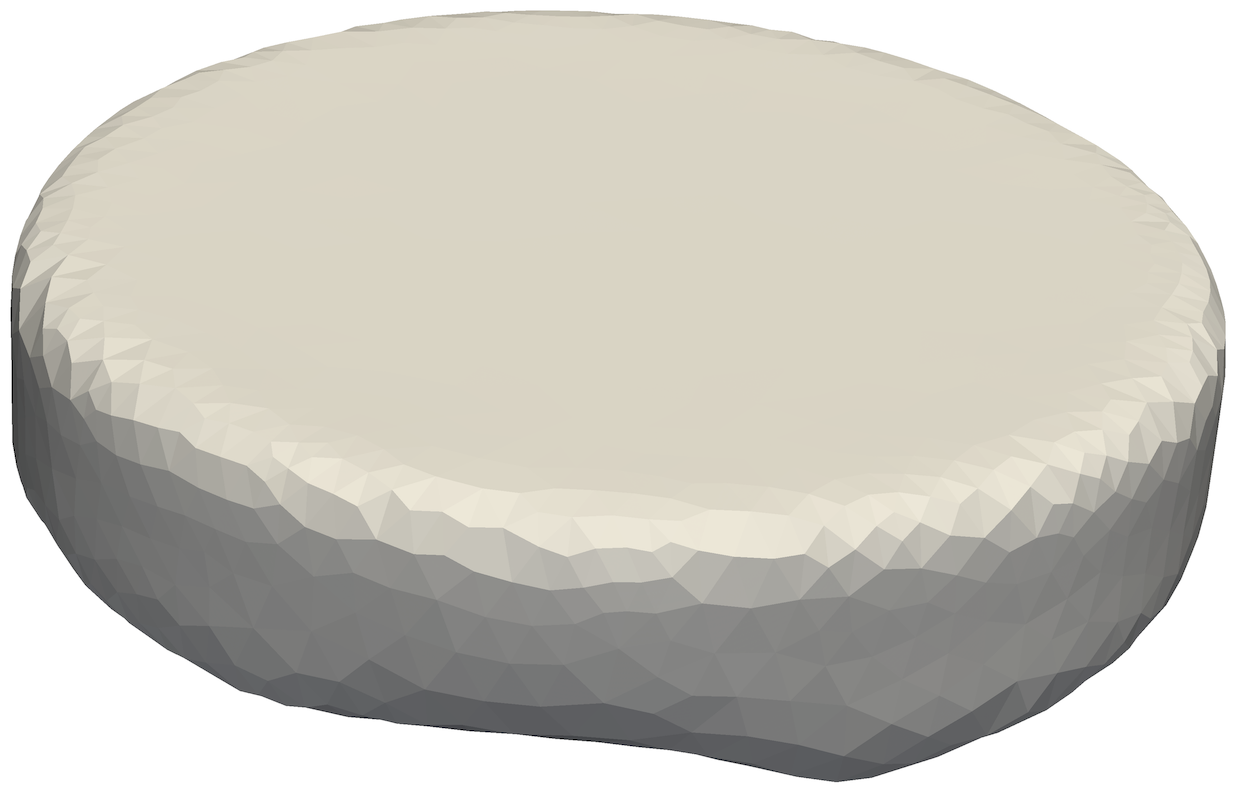}};
	\draw[fill=gray, gray] (3.5,1.8) -- (1.0, 3.25) -- (3.4, 4.15) -- (6, 3.0) -- (3.5, 1.8);

	\draw[color=black, ->,very  thick, -latex]   (3.5,1.8)--  (3.5, 2.6) ;	
	\draw[color=black, ->,very  thick, -latex]   (1.0, 3.25)--  (1.0, 4.05);
	\draw[color=black, ->,very  thick, -latex]   (3.4, 4.15)-- (3.4, 4.85);
	\draw[color=black, ->,very  thick, -latex]   (6, 3.0) --  ((6, 3.8);				
\end{tikzpicture}
}
\hfill
\subfloat{ }
\caption{a) Conchoidal fracture in quartz. Source of the original photo is \cite{quartz2018}.
 b) Rock geometry and the boundary value problem set-up.  The shaded area illustrates a part of the surface, where time-dependent Dirichlet boundary conditions are enforced.  We prescribe incremental displacement $\Delta \u = 5 \times 10^{-3}$mm while, we keep the phase field variable fixed, e.g. $c=0$.
 The bottom surface of the rock was held fixed during the entire simulation. }
\label{fig:CF_init}
\end{figure}

The  conchoidal fractures are commonly observed in fine-grained or amorphous materials such as rocks, minerals, and glasses.
They are fascinating because they do not follow any natural way of the separation and result in a curved breakage, which resembles the concentric undulations of a seashell \cite{hodgson1961classification}. 
Figure \ref{fig:cf_1} demonstrates an example of conchoidal fracture in quartz. 
From the modeling point of view, the conchoidal fracture is interesting because the crack initiates usually inside the body and not on the surface of the specimen. 
Therefore, the phase-field model needs to be capable of predicting the crack initiation and propagation without any initial crack or notch.

\begin{figure}[H]
\centering
\hfill
\subfloat[\label{fig:cf_v1}]{
			\includegraphics[scale=0.055]{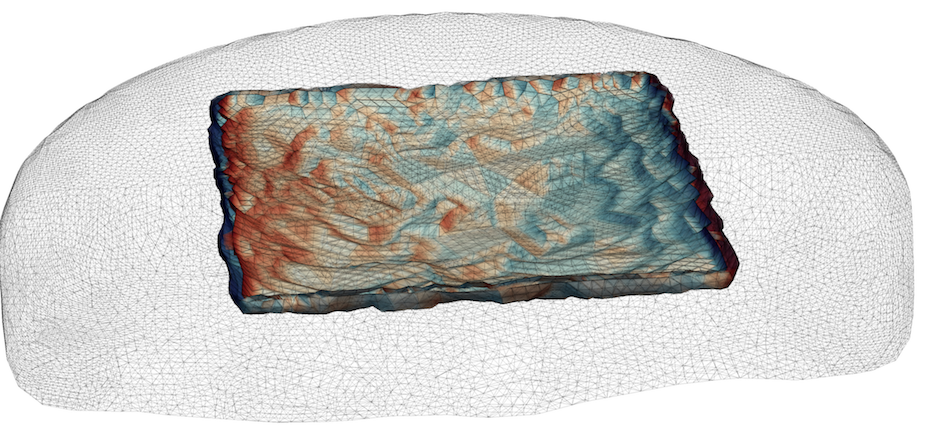}
}
\hfill
\subfloat[\label{fig:cf_v2}]{
			\includegraphics[scale=0.06]{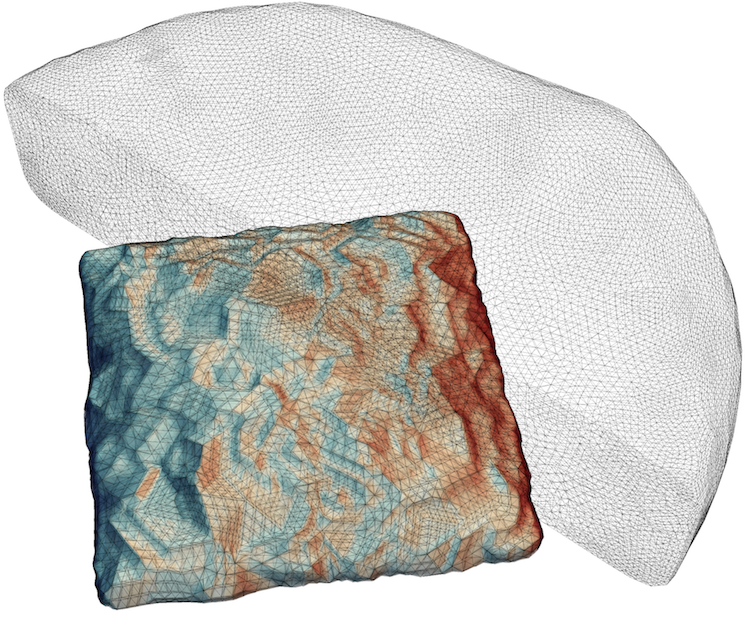}
}
\hfill
\subfloat{
\includegraphics[scale=0.04]{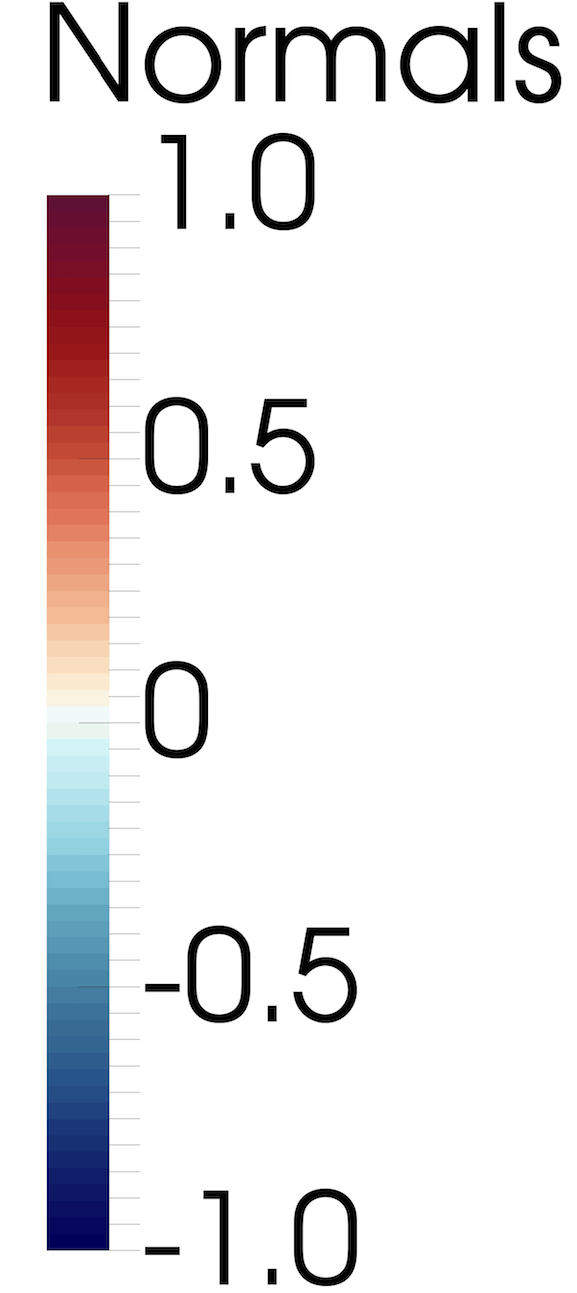}
}
\hfill
\subfloat{

}
\caption{Result of conchoidal fracture simulation, a numerical experiment with approx. $7,9$ millions number of dofs. 
A crack domain is shown inside of the rock geometry form a) top view b) bottom view. 
One should note, that the fracture reassembles curved patterns as expected for conchoidal fracture.}
\label{fig:CF_result}
\end{figure}

Our numerical experiment considers real rock geometry as depicted on Figure \ref{fig:cf_2}. 
The rock mesh is unstructured and contains rough surfaces, which can only be represented well with a relatively fine mesh. 
Following \cite{muller2016}, we simulate a conchoidal fracture by prescribing the displacement on the part of the upper boundary, see Figure \ref{fig:cf_2}.
On the same part of the boundary, the phase-field variable is enforced by prescribing Dirichlet boundary condition $c=0$.
On the lower boundary, we prescribe zero Dirichlet boundary condition for both displacement and the crack phase field.
The same approach, but on a parallelepipedal domain was studied by Bilgen et. al in \cite{bilgen2017phase}.  
The result of the simulation is demonstrated in Figure \ref{fig:CF_result}.

\subsection{Pneumatic fracture}
Our last numerical experiment considers pneumatic fracture. 
Pneumatic fracture involves the injection of air or fluid at a sufficient pressure which induces crack propagation. 
In order to simulate pneumatic fracturing, we follow the approach derived in  \cite{Wick15}, where a pressure term, $p$,  is added into energy functional ${\Psi}(\u, c) $ as follows
\begin{align}
{\Psi}(\u, c) := \int_{\Omega}   [g(c)(1-k) + k] \psi_e^+(\eps(\u)) + \psi_e^-(\eps(\u))  \  d \Omega +
\int_{\Omega}  \pazocal{G}_c \Big[   \frac{1}{2 l_s }c^2 + \frac{l_s}{2} | \nabla c |^2 \Big]  \ d \Omega + \int_{\Omega} (1-c)^2 p \nabla \cdot \u \ d \Omega. 
\label{eq:pressure_energy}
\end{align}

\begin{figure}[H]
\subfloat[\label{fig:pressure_1}]{
			\includegraphics[scale=0.027]{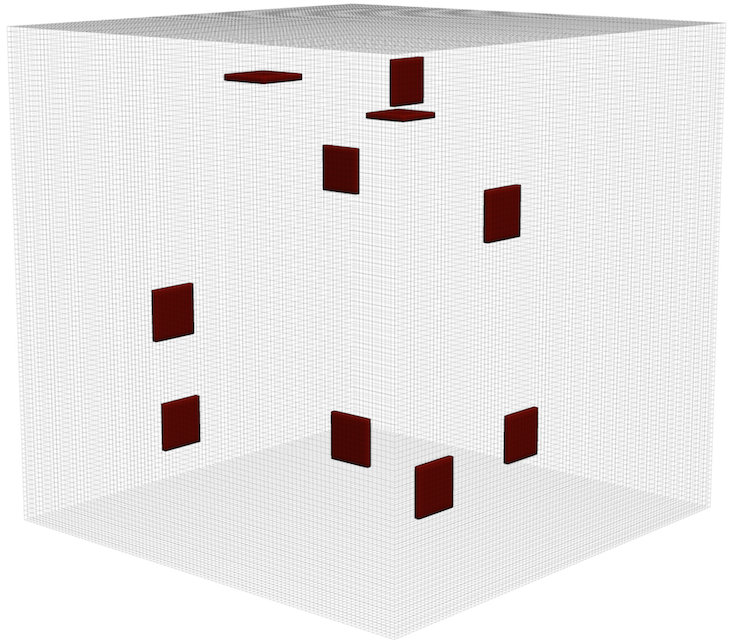}
}
\hfill
\subfloat[\label{fig:pressure_2}]{
			\includegraphics[scale=0.027]{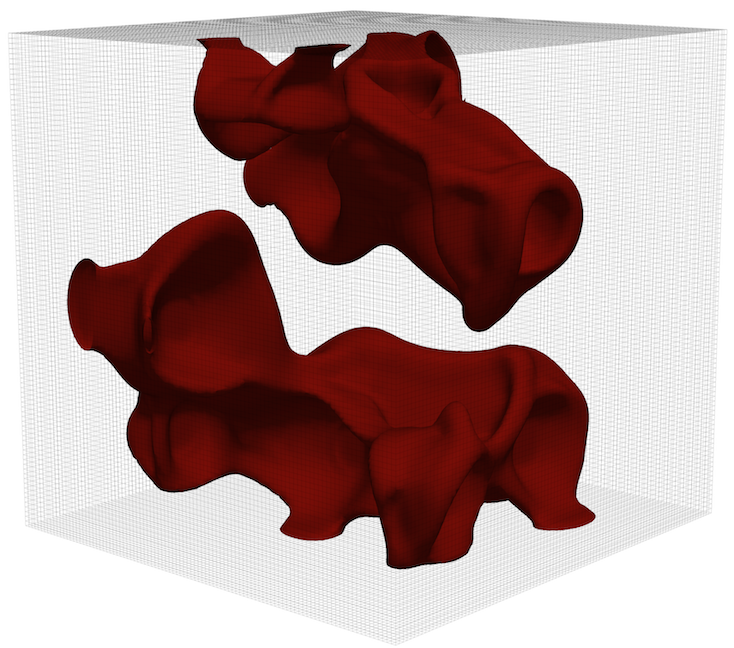}
}
\hfill
\subfloat[\label{fig:pressure_3}]{
			\includegraphics[scale=0.027]{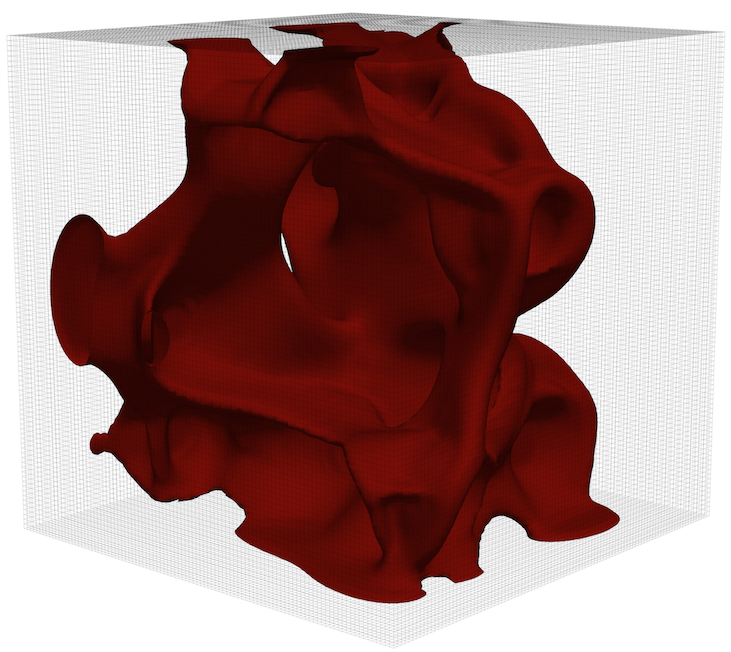}
}
\hfill
\subfloat[\label{fig:pressure_4}]{
			\includegraphics[scale=0.03]{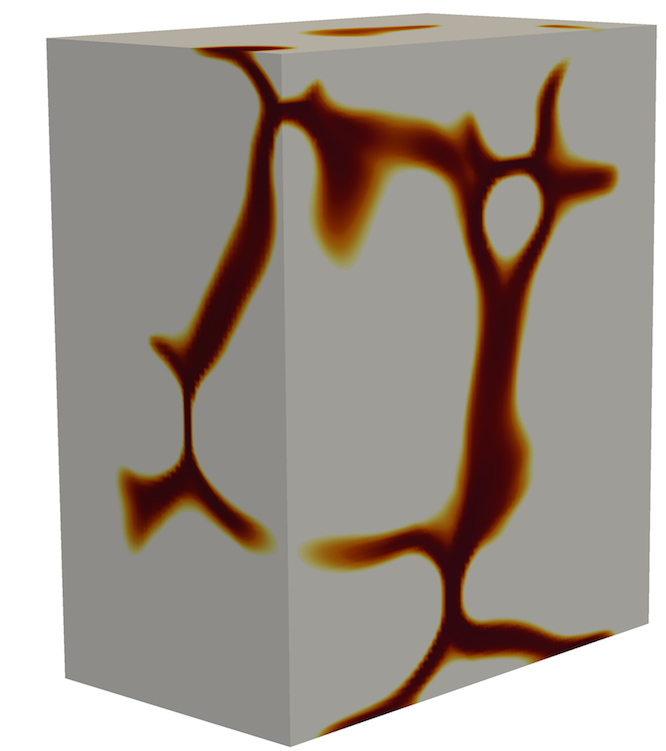}
}
\hfill
\subfloat{
\includegraphics[scale=0.032]{c_scale.png}
}
\caption{
Pressure induced fracture,  an example with approx. $4,1$ millions of dofs.
a) Initial set-up with ten randomly distributed fracture zones. 
b) Crack patterns at time $t=11$s. 
 c) Simulation result at time $t = 12$s.  One can observe, robust crack propagation between two subsequent time-steps.
 d) Cut across the computational domain at the end of the simulation, thus at $t = 12$s. The interesting crack patterns as well as crack merging and crack branching can be observed on several places. }
\label{fig:pressure_result}
\end{figure}

The numerical experiment was performed on a cube of size $1 \times 1 \times 1$ mm.
Initial set-up of the simulation takes into account ten randomly distributed fractured zones, see Figure \ref{fig:pressure_1}.
These initial cracks were specified by setting the value of the phase field to $c=1$ in those ten zones.
Through the simulation, we apply zero Dirichlet boundary conditions for the displacement field on all six sides of the domain. 
At each pseudo-time step $t$, the pressure load is linearly increased as $p(t) = p_0 + t p_0$, with an initial pressure $p_0=1$.
The simulation time step was set-up to $\Delta t=1$s. 
A similar approach was considered also in \cite{bilgen2017_pneumatic}. 

Figures \ref{fig:pressure_2} - \ref{fig:pressure_3} demonstrate evolution of the fracture by illustrating crack iso-surfaces for $c=0.9$. 
Figure \ref{fig:pressure_4} depicts the obtained results by  introducing a cut across the computational domain.
This allows us to observe interesting crack patterns and interacting fractures  inside of the computational domain.

\section{Convergence and performance study}
\label{section:numerical_results}
In this section, we investigate the performance of the RMTR method, Algorithm \ref{alg:rmtr}. 
Our investigation starts by analyzing the effect of different coarse level models on the performance of the RMTR. 
Further, we compare the RMTR method with the single level trust region solver, Algorithm \ref{alg:TR}. 
At the end, we examine the level-independence, complexity and scalability properties of the RMTR. 
While benchmarking, trust region based algorithms were set-up with the solver parameters as described in Table \ref{table:params_solver}. 
The solution strategies terminate, when the following stoping criterion is satisfied
\begin{align}
\| \pazocal{E}(\x_i, f) \| < 10^{-8} \  \text{or}  \  \| \x_{i-1} -  \x_{i} \|  < 10^{-14}  
\label{eq:stopping_criterium}
\end{align}
where the current iterate is denoted by $\x_{i}$. 
In the context of recursive multilevel trust region algorithm, $\x_{i}$ is considered to be the fine level iterate. 
The  criticality measure $\| \pazocal{E}(\x_i, f) \|$ is as defined in \eqref{eq:projection_onto_feas_set} and  $\| \x_{i-1} -  \x_{i} \|$ denotes the correction size.
The RMTR solver was configured as a V-cycle with one pre/post-smoothing step and two coarse solves. 
We solved the trust region subproblems on levels $l = 1, \dots, L$ approximately with 10 steps of the projected Conjugate-Gradient method preconditioned with Jacobi preconditioner \cite{nocedal2006numerical}. 
On the coarsest level, we employed the active set strategy \cite{benson2010toolkit} with the sparse linear solver MUMPS \cite{amestoy2000mumps}.
 \begin{table}
    \centering
    \caption{Choice of parameters used inside of TR/RMTR algorithms.}
    \label{table:params_solver}
    \begin{tabular}{ cccccccccc} 
    \toprule
    \textcolor{myred}{Parameter}     &        \textcolor{myred}{$\eta_1$}     &      \textcolor{myred}{$\eta_2$}  &     \textcolor{myred}{$\gamma_1$} &      \textcolor{myred}{$\gamma_2$} &     \textcolor{myred}{$\Delta_{0,0}^L$} &     \textcolor{myred}{$\mu_1$} &     \textcolor{myred}{$\mu_2$}  \\  \hline
       \textcolor{myred}{Value}    &     $0.1$&    $0.75$    &     $0.5$    & $2.0$    &     $1$  &     1 &    1 \\ \bottomrule

    \end{tabular}
 \end{table}

Our implementation of the parallel multilevel phase-field fracture simulations is based on the finite element framework MOOSE \cite{gaston2009moose}. 
Presented solution strategies were implemented as part of our open-source,  C++ library \utopia{} \cite{utopia}. 
\utopia{} is an embedded domain specific language (EDSL), designed to make implementations for parallel computing as transparent as possible. 
It enables us to write complex scientific code by using expressions similar to Matlab, whereas the complexity of parallelization is hidden in different back-ends.  
The presented results were computed with a \petsc{} \cite{balay2014petsc} back-end. 

All experiments were performed on the local cluster at the Institute of
Computational Science (ICS), Universit\`a della Svizzera Italiana. 
The cluster consists of  $42$ compute nodes, each equipped with 2 Intel R E5-2650 v3 processor with a clock frequency of 2.60GHz. 
Having 10 cores per processor leads to 20 cores per node. Memory per node is 64 GB.

\subsection{Effect of different coarse level models on the performance of the RMTR}
\label{section:results_different_models}
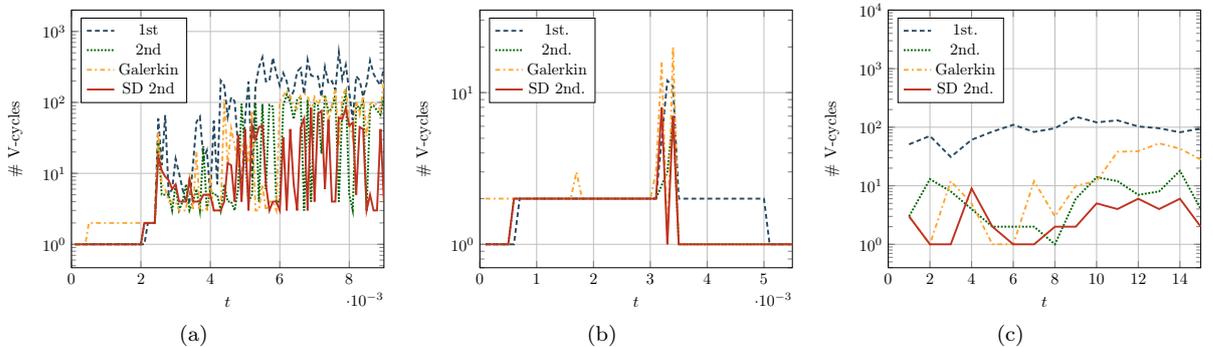
\begin{figure}[H]
\subfloat[\label{fig:tr_cm1}]{
\begin{tikzpicture}[scale=0.6]
\begin{axis}[legend entries={1st , 2nd , Galerkin, SD 2nd },
  legend pos=north west,
  xlabel={$t$}, ylabel={$\#$ V-cycles}, ymajorgrids=true, xmajorgrids=true, xmode=normal, ymode=log,     xmin=0, xmax=0.0090,  ymax=2000]

\addplot[color = myblue, very thick, densely dashed] table [y=its, x=time, col sep=comma]  {mixed_first.csv};
\addplot[color = mygreen, very thick, densely dotted] table [y=its, x=time, col sep=comma]  {mixed_second.csv};
\addplot[color = myyellow, very thick, dashdotted] table [y=its, x=time, col sep=comma]  {mixed_galerkin.csv};
\addplot[color = myred, very thick, solid] table [y=its, x=time, col sep=comma]  {mixed_mine.csv};

\end{axis}
\end{tikzpicture}
}
\hfill
\subfloat[\label{fig:tr_cm2}]{
\begin{tikzpicture}[scale=0.6]
\begin{axis}[legend entries={1st., 2nd., Galerkin, SD 2nd.},
  legend pos=north west,
  xlabel={$t$}, ylabel={$\#$ V-cycles}, ymajorgrids=true, xmajorgrids=true, xmode=normal, ymode=log, xmin=0, xmax=0.0055,  ymax=35]

\addplot[color = myblue, very thick, densely dashed] table [y=its, x=time, col sep=comma]  {conchoidal_first.csv};
\addplot[color = mygreen, very thick, densely dotted] table [y=its, x=time, col sep=comma]  {conchoidal_second.csv};
\addplot[color = myyellow, very thick, dashdotted] table [y=its, x=time, col sep=comma]  {conchoidal_galerkin.csv};
\addplot[color = myred, very thick, solid] table [y=its, x=time, col sep=comma]  {conchoidal_mine.csv};

\end{axis}
\end{tikzpicture}
}
\hfill
\subfloat[\label{fig:tr_cm3}]{
\begin{tikzpicture}[scale=0.6]
\begin{axis}[legend entries={1st. , 2nd. , Galerkin, SD 2nd.},
  legend pos=north west,
  xlabel={$t$}, ylabel={$\#$ V-cycles}, ymajorgrids=true, xmajorgrids=true, xmode=normal, ymode=log, ymax=10000, xmin=0, xmax=15]

\addplot[color = myblue, very thick, densely dashed] table [y=its, x=time, col sep=comma]  {pressure_first.csv};
\addplot[color = mygreen, very thick, densely dotted] table [y=its, x=time, col sep=comma]  {pressure_second.csv};
\addplot[color = myyellow, very thick, dashdotted] table [y=its, x=time, col sep=comma]  {pressure_galerkin.csv};
\addplot[color = myred, very thick, solid] table [y=its, x=time, col sep=comma]  {pressure_mine.csv};

\end{axis}
\end{tikzpicture}
}
\caption{Number of nonlinear V-cycles over time-steps for a) Mixed mode, specimen with $1,164,996$ dofs.  
 b) Conchoidal fracture, example with $7,866,166$ dofs.  c) Pressurized fracture, test case with $4,121,204$ dofs.
RMTR was set-up with 6 levels.
For each example, we depict  the performance of RMTR method with respect to different coarse level models. 
Please note, that x and y axes differ for all graphs.   }
\label{fig:num_rmtr_models}
\end{figure}

In this section, we investigate the performance of the RMTR method configured with different coarse level models, as described in Section \ref{section:coarse_level_models}. 
We monitor number of nonlinear V-cycles required by RMTR method to converge to the desired tolerance during the time step, see Figure \ref{fig:num_rmtr_models}.
As we can see from Figure \ref{fig:num_rmtr_models},  all variants of RMTR method converge to the desired tolerance.  However, convergence speed differs. 
Based on our numerical experiments, we conclude that the 1st order consistency approach \eqref{eq:coarse_objective_first_order}  yields the worst convergence rates. 
This is especially true in the post-stage, where fracture occurs and propagates. 
As mentioned in Section \eqref{section:fisrt_order}, the slow convergence is connected to the fact that the coarse level models are not suitable representations of the fine level fracture zones.  

Figure \ref{fig:num_rmtr_models} demonstrates that all three approaches: Galerkin, 1st order consistency and 2nd order consistency achieve worse convergence rates as our novel solution dependent 2nd order model \eqref{eq:coarse_objective_modified2nd_order}.  
We can compute achieved speed-up with respect to any other coarse level model as follows
\begin{align}
 \frac{\# \ \text{accumulated V-cycles}}{ \# \ \text{accumulated V-cycles with the solution dependent 2nd order model}}.
\end{align}
From the obtained results, Table \ref{table:speedup_of_models},  we can clearly see that usage of our solution dependent 2nd order coarse level model is beneficial and leads to significant speed-up. 
For example, for the mixed mode test, the speed-up equals $7.42$, when compared to 1st order consistency approach. 

All results reported in the reminder of this section employ our RMTR method with the coarse level models based on the solution dependent second order consistency approach, Section \ref{section:modified_second_order}. 

\begin{table*}[t]
\centering
\caption{Speedup of RMTR method configured with the solution dependent 2nd order model with respect to alternative variants, namely 1st order consistency, 2nd order consistency and Galerkin model.
Study performed on three examples: Mixed mode,  Conchoidal fracture, and Pressurized fracture.}
\label{table:speedup_of_models}
\begin{tabular}{ cccccccccc} 
\toprule
   								&    \textcolor{myred}{Mixed mode} 	&   \textcolor{myred}{Conchoidal fracture} &  \textcolor{myred}{Pressurized fracture}          \\ \hline
\textcolor{myred}{1st order}			&  $7.42$			& $1.36$  			& $27.86$	 	\\ 
\textcolor{myred}{2nd order}			& $2.13$ 		& $1.12$  			& $2.12$	 		\\ 
\textcolor{myred}{Galerkin}			& $2.65$ 		& $1.30$  			& $5.33$	 		\\ \bottomrule
\end{tabular}
\end{table*}

\subsection{Single level TR  vs. RMTR}
\begin{figure}[H]
\subfloat[\label{fig:tr_rmtr4}]{
\begin{tikzpicture}[scale=0.6]
\begin{axis}[legend entries={	TR, 
						RMTR},
  legend pos=north east,
  xlabel={$t$}, ylabel={$\#$ nonlin. it/ V-cycles}, ymajorgrids=true, xmajorgrids=true, xmode=normal, ymode=log,     xmin=0, xmax=0.0090,  ymax=10000]

\addplot[color = green4, very thick, densely dashdotted] table [y=its, x=time, col sep=comma]  {mixed_tr.csv};
\addplot[color = myred, very thick, solid] table [y=its, x=time, col sep=comma]  {mixed_mine.csv};

\end{axis}
\end{tikzpicture}
}
\hfill
\subfloat[\label{fig:tr_rmtr5}]{
\begin{tikzpicture}[scale=0.6]
\begin{axis}[legend entries={	TR, 
						RMTR},
  legend pos=north east,
  xlabel={$t$}, ylabel={$\#$ nonlin. it/ V-cycles}, ymajorgrids=true, xmajorgrids=true, xmode=normal, ymode=log,  ymax=10000, xmin=0, xmax=0.0055,  ymax=200]

\addplot[color = green4, very thick, densely dashdotted] table [y=its, x=time, col sep=comma]  {conchoidal_tr.csv};
\addplot[color = myred, very thick, solid] table [y=its, x=time, col sep=comma]  {conchoidal_mine.csv};

\end{axis}
\end{tikzpicture}
}
\hfill
\subfloat[\label{fig:tr_rmtr6}]{
\begin{tikzpicture}[scale=0.6]
\begin{axis}[legend entries={	TR, 
						RMTR},
  legend pos=north east,
  xlabel={$t$}, ylabel={$\#$ nonlin. it/ V-cycles}, ymajorgrids=true, xmajorgrids=true, xmode=normal, ymode=log, ymax=10000, xmin=0, xmax=15]

\addplot[color = green4, very thick, densely dotted] table [y=its, x=time, col sep=comma]  {pressure_tr.csv};
\addplot[color = myred, very thick, solid] table [y=its, x=time, col sep=comma]  {pressure_mine.csv};

\end{axis}
\end{tikzpicture}
}
\caption{Number of nonlinear iterations/V-cycles over time-steps for a) Mixed mode, specimen with $1,164,996$ dofs.  
 b) Conchoidal fracture, example with $7,866,166$ dofs.  c) Pressurized fracture, test case with $4,121,204$ dofs.
RMTR was set-up with 3 levels. Please note, that x and y axes differ for all graphs.   }
\label{fig:num_it_tr_rmtr_2}
\end{figure}
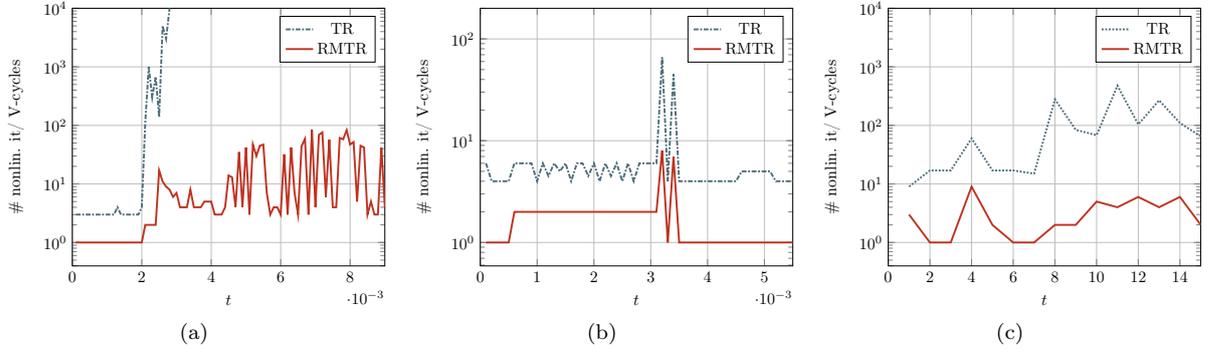

In this section, we compare the  performance of RMTR Algorithm \ref{alg:rmtr} with its single level variant, TR Algorithm  \ref{alg:TR}.
To this aim, we measure the number of nonlinear iterations required by TR method vs. the number of nonlinear V-cycles required by RMTR method. 
Figures \ref{fig:num_it_tr_rmtr} - \ref{fig:num_it_tr_rmtr_2} depict obtained results for different time-steps. 
As we can observe from Figures \ref{fig:num_it_tr_rmtr} - \ref{fig:num_it_tr_rmtr_2}, the number of nonlinear iterations grows once the simulation reaches the post-stage phase, where the crack finally propagates. 
This behavior is not surprising, as in the post-stage phase the non-convexity and the non-linearity become stronger.
For example, from the tension test in Figure \ref{fig:tr_rmtr1}, we can see that crack propagation starts around time $t=4.5 \times 10^{-3}$. 

Our numerical experiences indicate that the RMTR method always outperforms the single level TR method.
On the average, we observe a reduction in the number of required iterations approximately by factor of ten. 
The biggest decrease can be observed for the pressure example, where the number of iterations required by RMTR method is 32 times lower than the number of iterations required by the single level TR method.
Another interesting results were obtained for the mixed mode example, see Figure \ref{fig:tr_rmtr4}.
The single level TR algorithm encounters a failure, as the solver did not converge within $10, 000$ iterations.
At the same time, the multilevel variant converges without any problems. 
However, we have to remark that the single level TR method converges, while performing the simulation with smaller time-step.
\begin{figure}[H]
\subfloat[\label{fig:tr_rmtr1}]{
\begin{tikzpicture}[scale=0.59]
\begin{axis}[legend entries={	TR, 
						RMTR},
  legend pos=north west,
  xlabel={$t$}, ylabel={$\#$ nonlin. it/ V-cycles}, ymajorgrids=true, xmajorgrids=true, xmode=normal, ymode=log,     xmin=0, xmax=0.005,  ymax=1000]

\addplot[color = green4, very thick, densely dashdotted]  table [y=its, x=time, col sep=comma]  {tension_tr6.csv};
\addplot[color = myred, very thick, solid] table [y=its, x=time, col sep=comma] {tension_6grid.csv};
\end{axis}
\end{tikzpicture}
}
\hfill
\subfloat[\label{fig:tr_rmtr2}]{
\begin{tikzpicture}[scale=0.59]
\begin{axis}[legend entries={	TR, 
						RMTR},
  legend pos=north west,
  xlabel={$t$}, ylabel={$\#$ nonlin. it/ V-cycles}, ymajorgrids=true, xmajorgrids=true, xmode=normal, ymode=log,  ymax=500, xmin=0, xmax=0.025]
\addplot[color = green4, very thick, densely dashdotted] table [y=its, x=time, col sep=comma]   {shear_tr6.csv};
\addplot[color = myred, very thick, solid] table [y=its, x=time, col sep=comma]   {shear_6grid.csv};
\end{axis}
\end{tikzpicture}
}
\hfill
\subfloat[\label{fig:tr_rmtr3}]{
\begin{tikzpicture}[scale=0.59]
\begin{axis}[legend entries={	TR, 
						RMTR},
  legend pos=north west,
  xlabel={$t$}, ylabel={$\#$ nonlin. it/ V-cycles}, ymajorgrids=true, xmajorgrids=true, xmode=normal, ymode=log, ymax=100, xmin=0, xmax=0.003]
\addplot[color = green4, very thick, densely dashdotted]  table [y=its, x=time, col sep=comma]   {tear_tr6.csv};
\addplot[color = myred, very thick, solid] table [y=its, x=time, col sep=comma]   {tear_6grid.csv};
\end{axis}
\end{tikzpicture}
}
\caption{Number of nonlinear iterations/V-cycles required by TR/RMTR method  to converge for a) Tension  b) Shear  c) Tearing test. 
The RMTR was configured with 6 levels. Experiments performed on a brick with $6,010,884$ dofs. 
Please note, that scale for x and y axes are different for all graphs. }
\label{fig:num_it_tr_rmtr}
\end{figure}
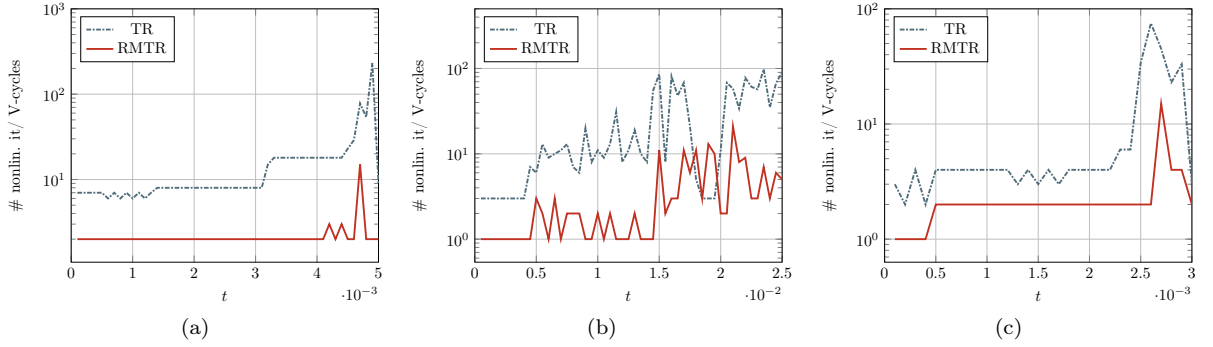

\subsection{Level-independence}
\label{section:robustness_levels}
The classical multilevel solvers, for instance multigrid are level-independent, in the sense, that  the number of iterations stays constant with an increased number of dofs. 
Here, we  investigate the level-independence of our RMTR algorithm by measuring number of required V-cycles for three fracture modes while increasing number of levels/dofs. 
The Table \ref{table:dofs} provides number of dofs used on each level of multilevel hierarchy as well as respective length-scale parameter $l_s$. 
It is of major importance to notice that the length-scale parameter $l_s$ decreases with increased resolution. 
As a consequence, volumetric approximation of the fracture reassembles sharp fracture surface more closely as refinement level increases. 

\begin{figure}[H]
\subfloat[\label{fig:rmtr_tension}]{
\begin{tikzpicture}[scale=0.59]
\begin{axis}[legend entries={4 levels,  5 levels,  6 levels},
  legend pos=north west,
  xlabel={$t$}, ylabel={$\#$ V-cycles}, ymajorgrids=true, xmajorgrids=true, xmode=normal, ymode=log,     xmin=0, xmax=0.005,  ymax=75]

\addplot[color = myblack, very thick, densely dotted] table [y=its, x=time, col sep=comma] {tension_4grid.csv};
\addplot[color = green3, very thick, solid] table [y=its, x=time, col sep=comma] {tension_5grid.csv};
\addplot[color = myred, very thick, densely dashed] table [y=its, x=time, col sep=comma] {tension_6grid.csv};

\end{axis}
\end{tikzpicture}
}
\hfill
\subfloat[\label{fig:rmtr_shear}]{
\begin{tikzpicture}[scale=0.59]
\begin{axis}[legend entries={4 levels,  5 levels,  6 levels},
  legend pos=north west,
  xlabel={$t$}, ylabel={$\#$ V-cycles}, ymajorgrids=true, xmajorgrids=true, xmode=normal, ymode=log,     xmin=0, xmax=0.025,  ymax=75]

\addplot[color = myblack, very thick, densely dotted] table [y=its, x=time, col sep=comma]  {shear_4grid.csv};
\addplot[color = green3, very thick, solid] table [y=its, x=time, col sep=comma] {shear_5grid.csv};
\addplot[color = myred, very thick, densely dashed] table [y=its, x=time, col sep=comma]  {shear_6grid.csv};

\end{axis}
\end{tikzpicture}
}
\hfill
\subfloat[\label{fig:rmtr_tear}]{
\begin{tikzpicture}[scale=0.59]
\begin{axis}[legend entries={ 4 levels,  5 levels,  6 levels},
  legend pos=north west,
  xlabel={$t$}, ylabel={$\#$ V-cycles}, ymajorgrids=true, xmajorgrids=true, xmode=normal, ymode=log, xmin=0, xmax=0.003,  ymax=75]

\addplot[color = myblack, very thick, densely dotted] table [y=its, x=time, col sep=comma]   {tear_4grid.csv};
\addplot[color = green3, very thick, solid] table [y=its, x=time, col sep=comma]{tear_5grid.csv};
\addplot[color = myred, very thick, densely dashed] table [y=its, x=time, col sep=comma]   {tear_6grid.csv};

\end{axis}
\end{tikzpicture}
}
\caption{ Number of nonlinear V-cycles required by RMTR method over time-steps for a) Tension  b) Shear  c) Tearing test. }
\label{fig:rmtr_h_ind}
\end{figure}
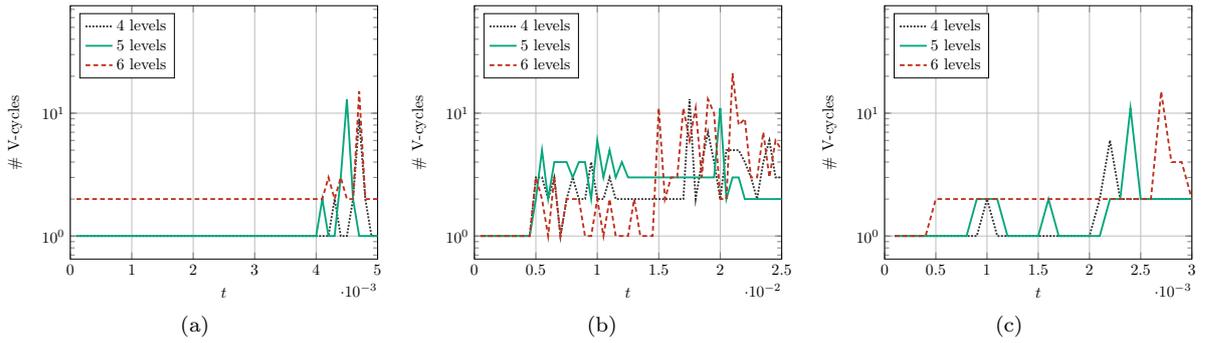

Figure \ref{fig:rmtr_h_ind} shows, that the number of required V-cycles increases slightly with the number of levels. 
This might be caused by two factors: 
Firstly, phase-field fracture problems discretized on different grid are not able to represent the same fractures. 
Therefore, simulating fracture propagation on different discretization level does not imply that we are solving the same problem. 
The second reason, why we are not able to observe level-independence might come directly from the use of the trust region globalization. 
During certain stages of the simulation,  the trust region radius might become too small. 
In this case, the sum of corrections coming from upper levels of the multilevel hierarchy already satisfy the constraints provided by the trust region radius on the finest level. 
As a consequence, recursion terminates before reaching the coarsest level. 
The low frequencies  of the error are not eliminated, which results in the higher number of nonlinear V-cycles. 
However, this is a very small price to pay in order to obtain a globally convergent multilevel solution strategy. 

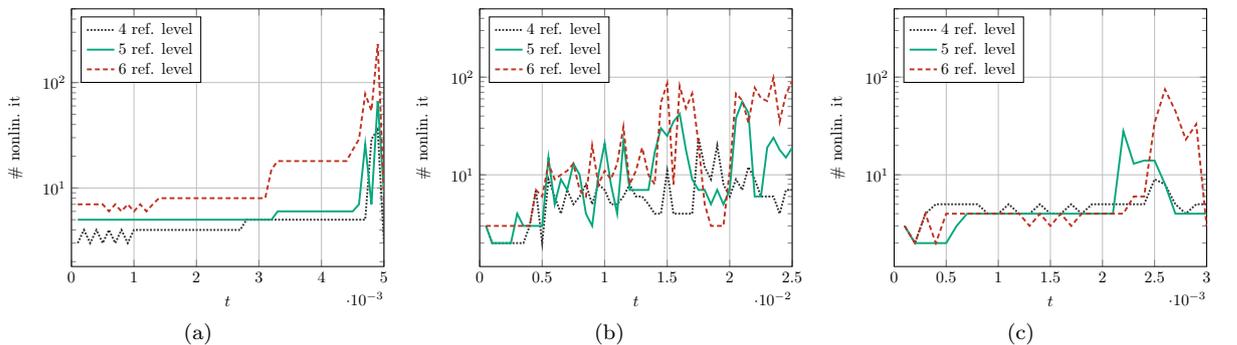
\begin{figure}[H]
\subfloat[\label{fig:tr_tension}]{
\begin{tikzpicture}[scale=0.6]
\begin{axis}[legend entries={4 ref.  level,  5 ref.  level,  6 ref.  level},
  legend pos=north west,
  xlabel={$t$}, ylabel={$\#$ nonlin. it}, ymajorgrids=true, xmajorgrids=true, xmode=normal, ymode=log,     xmin=0, xmax=0.005,  ymax=500]

\addplot[color = myblack, very thick, densely dotted] table [y=its, x=time, col sep=comma]  {tension_tr4.csv};
\addplot[color = green3, very thick, solid] table [y=its, x=time, col sep=comma]  {tension_tr5.csv};
\addplot[color = myred, very thick, densely dashed] table [y=its, x=time, col sep=comma] {tension_tr6.csv};

\end{axis}
\end{tikzpicture}
}
\subfloat[\label{fig:tr_shear}]{
\begin{tikzpicture}[scale=0.6]
\begin{axis}[legend entries={4 ref.  level,  5 ref.  level,  6 ref.  level},
  legend pos=north west,
  xlabel={$t$}, ylabel={$\#$ nonlin. it}, ymajorgrids=true, xmajorgrids=true, xmode=normal, ymode=log,     xmin=0, xmax=0.025,  ymax=500]

\addplot[color = myblack, very thick, densely dotted] table [y=its, x=time, col sep=comma]  {shear_tr4.csv};
\addplot[color = green3, very thick, solid] table [y=its, x=time, col sep=comma]  {shear_tr5.csv};
\addplot[color = myred, very thick, densely dashed] table [y=its, x=time, col sep=comma] {shear_tr6.csv};

\end{axis}
\end{tikzpicture}
}
\hfill
\subfloat[\label{fig:tr_tear}]{
\begin{tikzpicture}[scale=0.6]
\begin{axis}[legend entries={ 4 ref.  level,  5 ref.  level,  6 ref.  level},
  legend pos=north west,
  xlabel={$t$}, ylabel={$\#$ nonlin. it}, ymajorgrids=true, xmajorgrids=true, xmode=normal, ymode=log,  xmin=0, xmax=0.003,  ymax=500]

\addplot[color = myblack, very thick, densely dotted] table [y=its, x=time, col sep=comma]  {tear_tr4.csv};
\addplot[color = green3, very thick, solid] table [y=its, x=time, col sep=comma]  {tear_tr5.csv};
\addplot[color = myred, very thick, densely dashed] table [y=its, x=time, col sep=comma]  {tear_tr6.csv};

\end{axis}
\end{tikzpicture}
}
\caption{ Number of nonlinear iterations required by single level  TR method over time-steps for a) Tension  b) Shear  c) Tearing test. 
The performance of TR algorithm is monitored for different discretization levels. }
\label{fig:tr_h_ind}
\end{figure}

In comparison, Figure \ref{fig:tr_h_ind} demonstrates how number of required nonlinear iterations changes with increased refinement level for the single level TR method. 
We note, that growth in number of iterations is prevalent from one refinement level to the another. 
For example, for the tension example, Figure \ref{fig:tr_tension}, the accumulated number of iterations increases by factor $3.17$ from refinement level $5$ to refinement level $6$. 

\begin{table}[t]
\centering
\caption{Number of dofs used for testing three fracture modes. The number of dofs is displayed with respect to the different refinement level.
One should note, that the employed length-scale parameter $l_s$ also changes accordingly.  }
\label{table:dofs}
\begin{tabular}{ cccccccccc} 
\toprule
\textcolor{myred}{Level}    &    \textcolor{myred}{1} 	&   \textcolor{myred}{2} &  \textcolor{myred}{3} &  \textcolor{myred}{4} &  \textcolor{myred}{5}  &  \textcolor{myred}{6}           \\ \hline
\textcolor{myred}{\# dofs}			& $480$ 		& $2,484$  		& $15,300$	 	& $105,732$ 		& $782,340$ 	&  $6,010,884$ 		\\ 
\textcolor{myred}{$l_s$}			& $0.2$ 		& $0.1$  			& $0.05$	 		& $0.025$ 		& $0.0125$ 	&  $0.00625$ 			\\ \bottomrule
\end{tabular}
\end{table}

\subsection{Complexity and scalability}
\begin{figure}
\subfloat[\label{fig:perf1}]{
\begin{tikzpicture}[scale=0.8]
\begin{axis}[legend cell align={left}, legend entries={Simulation, Ideal}, xmode=log, ymode=log,   
 legend pos=north west, xlabel={Number of dofs}, ylabel={Time(s)}, ymajorgrids=true, xmajorgrids=true, 
    xmin=0, 
   ymin=0 
   ]
\addplot[color = myblue,  mark=oplus*, very thick] table [x=dofs, y=time, col sep=comma] {complexity.csv};
\addplot[color = myred,  mark=triangle*, very thick, densely dotted] table [x=dofs, y=optimal_time, col sep=comma] {complexity.csv};
\end{axis}
\end{tikzpicture}
}
\hfill
\subfloat[\label{fig:perf2}]{
\begin{tikzpicture}[scale=0.8]
\begin{axis}[legend cell align={left}, legend entries={Simulation, Ideal, New}, xmode=normal, ymode=normal,   
 legend pos=north west, xlabel={Number of cores}, ylabel={Relative speedup}, ymajorgrids=true, xmajorgrids=true, 
    xmin=0, 
   ymin=0 
   ]
\addplot[color = myblue,  mark=oplus*, very thick] table [x=cores, y=speedup, col sep=comma] {strong_scaling.csv};
\addplot[color = myred, mark=triangle*, very thick, densely dotted] table [x=cores, y=ideal, col sep=comma] {strong_scaling.csv};
\end{axis}
\end{tikzpicture}
}
\caption{Performance tests.  a) Computational complexity performed on the tension test.  b) Strong scaling test. Experiment was performed on tension test with $782, 340$ dofs.}
\label{fig:scaling}
\end{figure}
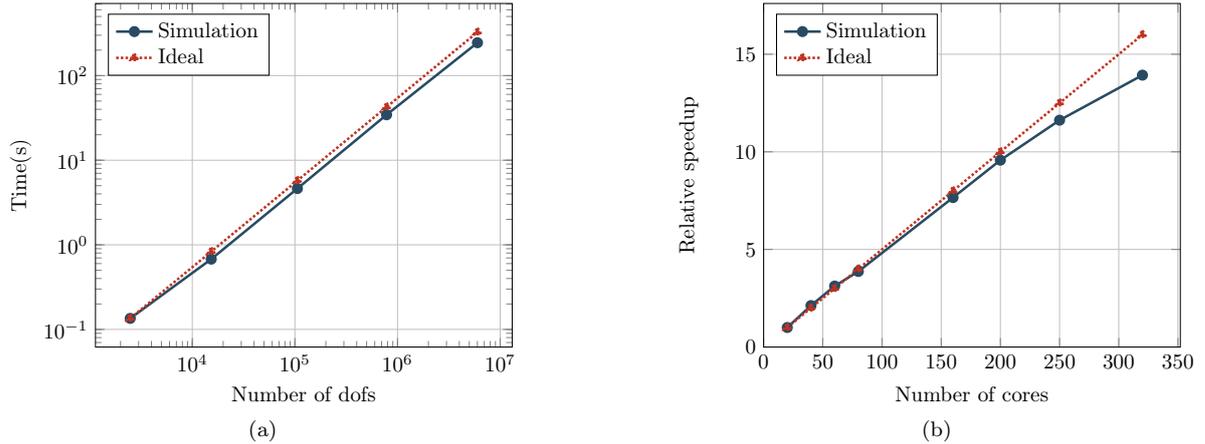

One of the most appealing property of any the multilevel algorithm is its optimal complexity $\mathcal{O}(n)$. 
This means, that the required solution time increases linearly with increased number of dofs. 
We can demonstrate the computational complexity by measuring the execution time of one V-cycle while increasing number of dofs. 
Our experiment considers tension test from \ref{section:fracture_modes}. 
The obtained results, depicted on Figure \ref{fig:perf1}, demonstrate that our implementation of the RMTR method is indeed of the optimal complexity. 
One should note, that the reported times are based on the results obtained on a serial machine.

Furthermore, we examine the strong scalability properties of our {\sc Utopia} implementation of the RMTR method. 
We consider tension test with fixed number of dofs, namely $782, 340$. 
Then, we measure required solution time while increasing number of cores. 
This allows us to compute the relative speedup by using following formula
\begin{align}
\text{relative speedup} = \frac{T_{20}}{T_p}, 
\end{align}
where $T_{20}$ represents computational time required by 20 cores and $T_p$ symbolizes the computational time required by $p$ cores.
Obtained results are reported on Figure \ref{fig:perf2}, which demonstrates that our implementation gives rise to almost ideal scaling up to 250 cores.
For more than $250$ processors, scaling is not ideal anymore, which is caused by an insufficient amount of dofs on the coarsest level. 
Further examination of scaling properties of the RMTR method as well as performing tests with larger examples (with more dofs) is left for the future work.

\begin{remark}
The measured time for both complexity and scalability tests contains also time spent for finite-element assembly. 
\end{remark}

\section{Conclusion}
\label{section:conclusion}
 In this work, we proposed a non-linear multi-level method for solving the monolithic systems arising from the phase-field fracture simulations.
Our method is based on recursive multi-level trust-region methods (RMTR), see \cite{Gratton2008recursive, Gross2009}. 
The proposed variant of RMTR method employs novel level dependent objective functions, which combine coarse level discretizations with a fine level representation of the fracture. 
The usage of this novel formulation improved the performance of the RMTR method substantially.
Furthermore, a comparison with single level trust region method was made and demonstrated a reduction in the number of iterations approximately by order of magnitude.
Performed experiments also indicated  that the our implementation of the RMTR method is of optimal complexity and shows good scaling behavior.

\appendix
\section{Euler-Lagrange equations}
\label{section:weak_form}
The first order necessary conditions for the minimization problem \eqref{eq:min_problem} are defined as:\\
Find pair $(\u, c) \in \bm{V}^t \times H^1(\Omega)$, such that
\begin{align}
\nabla_{\u} \Psi(\u, c; \bm{v}) = 0, \qquad  \forall \bm{v} \ \in \bm{H}^1_0(\Omega) \hspace{2cm}  \nabla_{c} \Psi(\u, c; w) \geq 0, \qquad  \forall w \ \in W^t \cap L^{\infty}(\Omega),
\label{eq:weak_form}
\end{align}
where $\bm{V}^{t} := \{ \u \in \bm{H}^1(\Omega) \ | \  \u = \bm{g}^t    \  \text{on} \ \  \partial \Omega_D  \}$ and $W^t:= \{ c \in H^1(\Omega) \  | \  0 \leq c^{t-1} \leq c \leq 1 \  \text{a.e.\ on} \ \Omega \}$.
The symbol $c^{t-1}$ denotes the solution of the phase-field on the previous time-step. 
The directional derivatives of $\Psi(\u, c)$ with respect to $\u$ and $c$, denoted by $\nabla_{\u} \Psi$ and $\nabla_{c} \Psi$, respectively are defined as
 \begin{equation}
 \begin{aligned}
\nabla_{\u} \Psi(\u, c; \mathbf{v}) &:= \int_{\Omega} \Big[  [g(c)(1-k) + k] \frac{\partial \psi_e^+}{\partial \eps}(\eps(\u)) + \frac{\partial \psi_e^-}{\partial \eps}(\eps(\u)) \ \Big] : \eps(\mathbf{v})  \ d \Omega, \\
\nabla_{c} \Psi(\u, c; w) &:= - \int_{\Omega} 2 (1-c) (1-k) \psi_e^+(\eps(\u)) w  \ d \Omega 
                                     + \int_{\Omega} \pazocal{G}_c \Big[   \frac{1}{ l_s }c w + l_s  \nabla c \nabla w \Big] \  d \Omega, 
\label{eq:weak_form2}
\end{aligned}
\end{equation}
where $k \approx 0$ denotes a regularization parameter.

\section{Additive decomposition of elastic energy}
\label{section:split}
Following \cite{miehe2010thermodynamically}, we additively decompose
the elastic energy density $\psi_e(\eps)$ into $\psi_e^+(\eps)$ due to the tension and $\psi_e^-(\eps)$
due to the compression by using following relation
\begin{align}
\psi_e^{\pm}(\eps) := \frac{1}{2} \lambda \Big( \big< \text{tr}(\eps) \big>^{\pm} \Big)^2 + \mu \ \eps^{\pm} : \eps^{\pm}\, .
\label{eq:energy}
\end{align}
The symbols $\mu$ and $\lambda$ in \eqref{eq:energy} denote the Lam\'e constants. 
The bracket symbol is defined as $\big< x \big>^{\pm}:= \frac{(|x| \pm x)}{2}$.
The strain tensor related to the desirable mode (tension, or compression) is obtained by the spectral decomposition of the principal strains
\begin{align}
\eps^{\pm} := \sum_{i=1}^d \big<\varepsilon_i \big>^{\pm} \bm{n}_i \otimes \bm{n}_i,
\end{align}
where $\bm{n}_i$ are the principal directions. 
The Cauchy stress tensor $\boldsymbol{\sigma} := \frac{\partial \psi}{\partial {\eps}}$ is also decomposed in an additive manner as
\begin{align}
\boldsymbol{\sigma} := [(1-c)^2 (1-k) + k] \boldsymbol{\sigma}^+ - \boldsymbol{\sigma}^-,
\end{align}
with 
\begin{align}
 \boldsymbol{\sigma}^{\pm} :=   \sum_{i=1}^d \big[ \lambda \big< \varepsilon_1 + \varepsilon_2 + \varepsilon_3 \big>^{\pm}  + 2 \mu \big< \epsilon_i \big>^{\pm} \big] \bm{n}_i \otimes \bm{n}_i.
\end{align}

\section{Projection onto feasible set}
\label{section:projection_set}
Let us consider feasible set $\pazocal{F}:= \{ \x \in \R^n \ | \ \l \leq \x \leq \u \}$, 
where the inequalities are meant component-wise. 
The projection $\mathcal{P}(\y)$ of vector $\y$ onto $\pazocal{F}$ is also defined component-wise as
\begin{align}
(\mathcal{P}(\y))_k := 
\begin{cases}
 (\l)_k 	&  \text{if}  \ \ \  (\y)_k \leq  (\l)_k \\
 (\y)_k	&  \text{if}  \ \ \  (\l)_k \leq  (\y)_k \leq  (\u)_k \\
 (\u)_k  	& \text{if}  \ \ \   (\u)_k \leq  (\y)_k
\end{cases}
\label{eq:projection_onto_feas_set}
\end{align}
for $\forall k \in \{ 1, \dots, n \}$.

\section*{Acknowledgments}
The authors thank Maria Nestola, Patrick Zulian and Hardik Kothari for many useful discussions and suggestions.
The authors thank Cyrill von Planta for providing rock geometry. 
This work was supported by Swiss National Science Foundation (SNF) and DFG (Deutsche Forschungsgemeinschaft) under the project SPP1748 "Reliable Simulation Techniques in Solid Mechanics. Development of Non-standard Discretization Methods, Mechanical and Mathematical Analysis".

\bibliography{literature}

\end{document}